\documentclass[10pt]{amsart}
\usepackage{amscd}
\usepackage{amssymb}
\usepackage[all]{xy}
\date{6 April 2006}
\title[Dualizing Complexes and Perverse Sheaves]
{Dualizing Complexes and Perverse Sheaves on 
Noncommutative Ringed Schemes}

\author{Amnon Yekutieli and James J. Zhang }
\address{A. Yekutieli: Department of  Mathematics 
Ben Gurion University, 
Be'er Sheva 84105, 
Israel}
\email{amyekut@math.bgu.ac.il}
\address{J.J. Zhang: Department of Mathematics, Box 354350,
University of Washington, Seattle, Washington 98195, USA}
\email{zhang@math.washington.edu}
\thanks{
Primary: 14A22; Secondary: 14F05, 14J32, 16E30, 16D90, 18E30.}
\thanks{This research was supported by the US-Israel Binational
Science Foundation. The second author was partially supported by 
the US National Science Foundation.}

\newtheorem{thm}[equation]{Theorem}
\newtheorem{cor}[equation]{Corollary}
\newtheorem{prop}[equation]{Proposition}
\newtheorem{lem}[equation]{Lemma}
\theoremstyle{definition}
\newtheorem{dfn}[equation]{Definition}
\newtheorem{rem}[equation]{Remark}
\newtheorem{exa}[equation]{Example}

\newtheorem{conv}[equation]{Convention}

\numberwithin{equation}{section}

\newcommand{\iso}{\xrightarrow{\simeq}}
\newcommand{\inj}{\hookrightarrow}

\newcommand{\xar}{\xrightarrow}
\newcommand{\opn}{\operatorname}
\newcommand{\opnt}[1]{\mathrm{#1}} 
\newcommand{\cat}[1]{\operatorname{\mathsf{#1}}}

\newcommand{\rmitem}[1]{\item[\text{\textup{(#1)}}]}
\newcommand{\mfrak}[1]{\mathfrak{#1}}
\newcommand{\mcal}[1]{\mathcal{#1}}
\newcommand{\msf}[1]{\mathsf{#1}}
\newcommand{\mbf}[1]{\mathbf{#1}}
\newcommand{\mrm}[1]{\mathrm{#1}}
\newcommand{\mbb}[1]{\mathbb{#1}}

\newcommand{\tup}[1]{\textup{#1}}
\newcommand{\bsym}[1]{\boldsymbol{#1}}
\newcommand{\boplus}{\bigoplus\nolimits}
\newcommand{\wtil}[1]{\tilde{#1}}

\newcommand{\bra}[1]{\langle #1 \rangle}
\renewcommand{\k}{\Bbbk}
\newcommand{\bwedge}{\bigwedge\nolimits}

\begin{document}

\begin{abstract}
A quasi-coherent ringed scheme is a pair $(X, \mcal{A})$,
where $X$ is a scheme, and $\mcal{A}$ is a noncommutative
quasi-coherent $\mcal{O}_{X}$-ring. 
We introduce dualizing complexes over 
quasi-coherent ringed schemes and study their properties. 
For a separated differential quasi-coherent ringed scheme
of finite type over a field, we prove existence and uniqueness of 
a rigid dualizing complex. In the proof we use the theory of
perverse coherent sheaves in order to glue local pieces of the 
rigid dualizing complex into a global complex.
\end{abstract}

\maketitle

\setcounter{section}{-1}
\section{Introduction}

The ``classical'' Grothendieck duality theory, dealing with 
dualizing complexes over schemes, was developed in 
the book {\em Residues and Duality} by Hartshorne \cite{RD}. 
Dualizing sheaves and complexes have important roles in several 
areas of algebraic geometry, including moduli spaces, 
resolution of singularities, arithmetic geometry and enumerative 
geometry. 
Various refinements, generalizations and explicit reformulations 
of Grothendieck's theory have appeared since \cite{RD};
a partial list of papers is 
\cite{Kl}, \cite{Li}, \cite{HK}, \cite{Ye2}, \cite{Ne}, 
\cite{AJL} and \cite{Co}. 

A noncommutative affine duality 
theory was introduced in \cite{Ye1}. By 
``affine'' we mean that this theory deals with noncommutative 
algebras over a base field $\k$. In the decade since its 
introduction the theory of noncommutative dualizing complexes has 
progressed in several directions, and it has applications in 
noncommutative algebraic geometry, ring theory, representation 
theory and even mathematical physics. Here is a sample of
papers: \cite{VdB1}, 
\cite{Jo}, \cite{MY}, \cite{WZ}, \cite{Ch2}, \cite{EG}, 
\cite{KKO}, \cite{NV}, \cite{LR} and \cite{AKO}. 
The definition of dualizing complex over a noncommutative algebra 
is recalled in Section \ref{sec2}.

The aim of this paper is to study Grothendieck duality on 
noncommutative spaces. As motivation one should consider the 
role Grothendieck duality plays both in commutative algebraic 
geometry and in noncommutative ring theory. Furthermore, 
recent developments, mainly surrounding homological mirror symmetry 
(cf.\ \cite{BO} and \cite{Do}),
tell us that there is a profound interplay between algebraic 
geometry, noncommutative algebra and derived categories. Duality 
for noncommutative spaces sits right in the middle of these areas 
of research. 

Actually there is one type of non-affine noncommutative space for 
which duality has already been studied. These are the noncommutative 
projective schemes $\opn{Proj} A$ of Artin-Zhang \cite{AZ}. Here 
$A$ is a noetherian connected graded $\k$-algebra satisfying the 
$\chi$-condition. A global duality theory for $\opn{Proj} A$ was 
discussed in \cite{YZ1} and \cite{KKO}. However in the present 
paper we choose to stay closer to the \cite{RD} paradigm, namely 
to develop a theory that has both local and global aspects. 

Here are a few issues one should consider before proposing
a theory of Gro\-then\-dieck duality on noncommutative spaces.
The first is to decide what is meant by a 
noncommutative space -- there are several reasonable choices in 
current literature. The second is to find a suitable formulation of 
duality, which shall include in a natural way the established 
commutative and noncommutative theories. The third issue is 
whether this duality theory applies to a wide enough 
class of spaces. 

The noncommutative spaces we shall concentrate on are the 
noetherian quasi-coherent ringed schemes. A quasi-coherent ringed 
scheme is a pair $(X, \mcal{A})$, where $X$ is a scheme over $\k$, 
and $\mcal{A}$ is a (possibly noncommutative) quasi-coherent 
$\mcal{O}_{X}$-ring (see Definition \ref{dfn3.1}).
This type of geometric object includes schemes 
($\mcal{A} = \mcal{O}_{X}$), 
noncommutative rings ($X = \opn{Spec} \k$) and rings of differential 
operators ($\mcal{A} = \mcal{D}_X$). Our definition of dualizing 
complex extends the established definitions in all these cases. A 
dualizing complex $\mcal{R}$ over 
$(X, \mcal{A})$ is an object of the derived category
$\msf{D}(\cat{Mod} \mcal{A}^{\mrm{e}})$, where
$(X^2, \mcal{A}^{\mrm{e}})$ is the product of $(X, \mcal{A})$
and its opposite $(X, \mcal{A}^{\mrm{op}})$. The duality 
functor 
$\mrm{D} : \msf{D}^{\mrm{b}}_{\mrm{c}}(\cat{Mod} \mcal{A})
\to \msf{D}^{\mrm{b}}_{\mrm{c}}(\cat{Mod} \mcal{A}^{\mrm{op}})$
is a contravariant Fourier-Mukai transform with respect to
$\mcal{R}$. See Definitions \ref{dfn4.1} and \ref{dfn5.8} 
for details.

We shall restrict our attention to dualizing complexes 
that have a ``local'' behavior, namely those that are supported 
on the diagonal $\Delta(X) \subset X^2$. 
Moreover we are interested in dualizing complexes that are 
canonical, or functorial, in a suitable sense -- something 
resembling Grothendieck's dualizing complex 
$\pi^! \k \in \msf{D}^{\mrm{b}}_{\mrm{c}}(\cat{Mod} \mcal{O}_X)$,
where $\pi : X \to \opn{Spec} \k$ is the structural morphism
(see \cite{RD}). Thus we propose to consider {\em rigid dualizing 
complexes}.  

Let $A$ be a noetherian $\k$-algebra. A rigid dualizing 
complex over $A$ is a dualizing complex $R_A$ equipped with 
an isomorphism 
\[ \rho_A : R_A \iso \opn{RHom}_{A^{\mrm{e}}}(A, R_A \otimes R_A) 
\] 
in $\msf{D}(\cat{Mod} A^{\mrm{e}})$, called a 
{\em  rigidifying isomorphism}. This notion is due to Van den Bergh 
\cite{VdB1}. It is known that the pair $(R_A, \rho_A)$ 
is unique up to a unique isomorphism in 
$\msf{D}(\cat{Mod} A^{\mrm{e}})$. 

Passing to sheaves, a rigid dualizing complex over $\mcal{A}$ is 
a dualizing complex $\mcal{R}_{\mcal{A}}$ supported on the diagonal 
in $X^2$, together with a collection $\bsym{\rho} = \{ \rho_U \}$ 
of rigidifying isomorphisms, indexed by the affine open sets 
$U \subset X$. For any such $U$ the complex
$R_A := \mrm{R} \Gamma(U^2, \mcal{R}_{\mcal{A}})$
is a dualizing complex over the ring $A := \Gamma(U, \mcal{A})$, 
and $\rho_A := \rho_U$ is a rigidifying isomorphism for $R_A$.
Moreover, the isomorphisms $\rho_U$ are required to 
satisfy a compatibility condition (see Definition \ref{dfn7.5}).

Our main objective in this paper is to prove existence and 
uniqueness of a rigid dualizing complex over $\mcal{A}$ 
(under suitable assumptions). 

A class of noncommutative spaces for which this can be done
is that of {\em separated differential 
quasi-coherent ringed schemes of finite type over $\k$}.
Suppose $(X, \mcal{A})$ is such a space. 
By definition $X$ is a separated finite type $\k$-scheme, 
and $\mcal{A}$ has a filtration such that 
the graded ring $\opn{gr} \mcal{A}$ is a finite type 
quasi-coherent $\mcal{O}_X$-algebra with a big center. Such a 
filtration is called a {\em differential filtration}. See 
Definition \ref{dfn13.3} for details. The prototypical examples 
are: 
\begin{enumerate}
\item ${\mcal A}$ is a coherent ${\mcal O}_X$-algebra 
(e.g.\  $\mcal{A} = \mcal{O}_X$, or $\mcal{A}$ is 
an Azumaya algebra with center $\mcal{O}_X$);
\item ${\mcal A}$ is the ring ${\mcal D}_X$ 
of differential operators on a smooth scheme $X$ in characteristic
$0$; and
\item ${\mcal A}$ is the universal enveloping algebra 
$\mrm{U}(\mcal{O}_X, \mcal{L})$ of a coherent 
Lie algebroid $\mcal{L}$ on $X$.
\end{enumerate}
In (1) and (3) there are no regularity assumptions on $X$, 
$\mcal{A}$ or $\mcal{L}$. More examples are provided in Section 
\ref{sec6}.

Suppose $(X, \mcal{A})$ is a differential 
quasi-coherent ringed scheme of finite type, $U$ is an affine open 
set and $A := \Gamma(U, \mcal{A})$. The ring $A$ is then a 
differential $\k$-algebra of finite type. In \cite{YZ4} 
we proved that $A$ has a rigid dualizing complex $R_A$, which
is supported on the diagonal in $U^2$.

The fact that $R_A$ is supported on the diagonal implies that it 
sheafifies to a complex 
$\mcal{R}_{\mcal{A}|_U} \in
\msf{D}(\cat{Mod} \mcal{A}^{\mrm{e}}|_{U^2})$,
which is a dualizing complex over $(U, \mcal{A}|_U)$. 
Because of the uniqueness of rigid dualizing complexes we obtain 
isomorphisms 
\begin{equation} \label{eqn0.1}
\mcal{R}_{\mcal{A}|_{U_1}}|_{U_1^2 \cap U_2^2} \cong 
\mcal{R}_{\mcal{A}|_{U_2}}|_{U_1^2 \cap U_2^2}
\end{equation}
in 
$\msf{D}(\cat{Mod} \mcal{A}^{\mrm{e}}|_{U_1^2 \cap U_2^2})$
for pairs of affine open sets $U_1$ and $U_2$, 
and these isomorphisms satisfy the cocycle condition on triple 
intersections.

The next stage is a process of gluing. Here we encounter a 
genuine problem: usually objects in derived categories cannot be 
glued. Grothendieck's solution in the commutative case, in 
\cite{RD}, was to use Cousin complexes. However, as explained 
in \cite{YZ3}, this solution seldom applies in the noncommutative 
context. The discovery at the heart of our present paper is that 
{\em perverse coherent sheaves can be used instead of Cousin 
complexes for gluing dualizing complexes}.

The concepts of t-structure and perverse sheaf are due to Bernstein, 
Beilinson and Deligne \cite{BBD}, and we recall the definitions 
in Section \ref{sec4}. It turns out that rigid dualizing complexes 
give rise to a perverse t-structure
on $\msf{D}^{\mrm{b}}_{\mrm{c}}(\cat{Mod} \mcal{A})$.
Indeed, for any affine open set $U$ we have a rigid dualizing 
complex $R_A$ over $A := \Gamma(U, \mcal{A})$. Let 
$\mrm{D} : \msf{D}^{\mrm{b}}_{\mrm{f}}(\cat{Mod} A) \to
\msf{D}^{\mrm{b}}_{\mrm{f}}(\cat{Mod} A^{\mrm{op}})$
be the duality functor $\opn{RHom}_{A}(-, R_A)$. 
The {\em rigid perverse t-structure} on 
$\msf{D}^{\mrm{b}}_{\mrm{c}}(\cat{Mod} A)$
is defined by
\[ \begin{aligned}
{}^{p}\msf{D}^{\mrm{b}}_{\mrm{c}}(\cat{Mod} A)^{\leq 0} & := & 
\{ M \mid \mrm{H}^i \mrm{D} M = 0 \text{ for all } i < 0 \}, \\
{}^{p}\msf{D}^{\mrm{b}}_{\mrm{c}}(\cat{Mod} A)^{\geq 0} & := & 
\{ M \mid \mrm{H}^i \mrm{D} M = 0 \text{ for all } i > 0 \} .
\end{aligned} \]
The intersection of these subcategories (the heart) is denoted by
${}^{p}\msf{D}^{\mrm{b}}_{\mrm{c}}(\cat{Mod} A)^{0}$.
Next let ``$\star$'' denote either ``$\leq 0$'', ``$\geq 0$'' or 
``$0$'', and define 
\[ \begin{aligned}
{}^{p}\msf{D}^{\mrm{b}}_{\mrm{c}}(\cat{Mod} \mcal{A})^{\star} 
& := \big\{ \mcal{M} \in  
\msf{D}^{\mrm{b}}_{\mrm{c}}(\cat{Mod} \mcal{A}) \mid
\mrm{R} \Gamma(U, \mcal{M}) \in 
{}^{p}\msf{D}^{\mrm{b}}_{\mrm{f}}
\big( \cat{Mod} \Gamma(U, \mcal{A}) \big)^{\star} \\
& \hspace{33ex} \text{ for all affine open sets } U \big\} .
\end{aligned} \]

Here is the first main result of the paper.

\begin{thm} \label{thm0.1}
Let $(X, \mcal{A})$ be a differential quasi-coherent ringed 
$\k$-scheme of finite type.  Then: 
\begin{enumerate}
\item The pair 
\[ \bigl( 
{}^{p}\msf{D}^{\mrm{b}}_{\mrm{c}}(\cat{Mod} \mcal{A})^{\leq 0}, 
{}^{p}\msf{D}^{\mrm{b}}_{\mrm{c}}(\cat{Mod} \mcal{A})^{\geq 0} 
\bigr) \]
is a t-structure on 
$\msf{D}^{\mrm{b}}_{\mrm{c}}(\cat{Mod} \mcal{A})$.
\item The assignment 
$V \mapsto 
{}^{p}\msf{D}^{\mrm{b}}_{\mrm{c}}(\cat{Mod} \mcal{A}|_{V})^{0}$,
for $V \subset X$ open, is a stack of abelian categories on $X$.
\end{enumerate} 
\end{thm}

Part (2) of the theorem says that 
${}^{p}\msf{D}^{\mrm{b}}_{\mrm{c}}(\cat{Mod} \mcal{A})^{0}$
behaves like the category of coherent sheaves $\cat{Coh} \mcal{A}$;
hence its objects are called {\em perverse coherent sheaves}. 
The theorem is restated as Theorem \ref{thm9.1}. 

It is not hard to show that the product 
$(X^2, \mcal{A}^{\mrm{e}})$ exists, and that it too is 
a differential quasi-coherent ringed scheme of finite type.
So by Theorem \ref{thm0.1} the rigid perverse t-structure on 
$\msf{D}^{\mrm{b}}_{\mrm{c}}(\cat{Mod} \mcal{A}^{\mrm{e}})$
exists, and perverse $\mcal{A}^{\mrm{e}}$-modules can be glued. 

According to \cite{YZ4}, for every affine open set $U$
the complex $\mcal{R}_{\mcal{A}|_U}$ is a perverse coherent 
$\mcal{A}^{\mrm{e}}|_{U^2}$-module, i.e.\   
$\mcal{R}_{\mcal{A}|_U} \in 
{}^{p}\msf{D}^{\mrm{b}}_{\mrm{c}}(\cat{Mod} \mcal{A}^{\mrm{e}}
|_{U^2})^{0}$.
Hence the gluing data (\ref{eqn0.1}) gives rise to a global 
dualizing complex 
$\mcal{R}_{\mcal{A}}$. Thus we obtain the second main result of 
our paper (which is repeated as Theorem \ref{thm13.8}):

\begin{thm} \label{thm0.2}
Let $(X, \mcal{A})$ be a separated differential 
quasi-coherent ringed $\k$-scheme of finite type.
Then there exists a rigid dualizing complex 
$(\mcal{R}_{\mcal{A}}, \bsym{\rho})$ over $\mcal{A}$. 
It is unique up to a unique isomorphism in 
$\msf{D}^{\mrm{b}}_{\mrm{c}}(\cat{Mod} \mcal{A}^{\mrm{e}})$.
\end{thm}

We also prove that the rigid trace exists for a finite 
centralizing morphism $f : (Y, \mcal{B}) \to (X, \mcal{A})$;
see Theorem \ref{thm13.9}. 

Here is a section-by-section synopsis of the paper. In Section 1 
we define quasi-coherent ringed schemes, discuss criteria for 
existence of the product 
$(X, \mcal{A}) \times (Y, \mcal{B})$ of two quasi-coherent ringed 
schemes, and consider some properties of the derived category
$\msf{D}(\cat{Mod} \mcal{A})$. In Section 2 we define dualizing 
complexes over quasi-coherent ringed schemes, and present a couple 
of ``exotic'' examples. Section 3 is about rigid dualizing 
complexes: their definition and some properties. In Section 4 we 
study t-structures on ringed spaces. The main result here is that 
a local collection of t-structures can be glued to a global 
t-structure (Theorems \ref{thm4.10} and \ref{thm4.11}).
In Section 5 we put together all the previous results to prove 
Theorem \ref{thm0.2}. Finally, in Section 6 we give examples of
differential quasi-coherent ringed schemes and their rigid 
dualizing complexes. We also take a close look at the commutative 
case $\mcal{A} = \mcal{O}_X$.

\medskip \noindent
\textbf{Acknowledgments.}
The authors wish to thank Eitan Bachmat, Joseph Bernstein, 
Sophie Chemla, Masaki Kashiwara, David Kazhdan, Maxim Kontsevich, 
Thierry Levasseur and Michel Van den Bergh for helpful 
conversations. We also thank the referee for reading the paper 
carefully and suggesting some corrections.

\section{Quasi-Coherent Ringed Schemes}
\numberwithin{equation}{section}

Throughout the paper $\k$ is a base field. All rings are by 
default $\k$-algebras, and all bimodules are central over $\k$.
Let $(X, \mcal{A})$ be a ringed space over $\k$. 
Thus $X$ is a topological space, and $\mcal{A}$ is a sheaf of 
(possibly noncommutative) $\k$-algebras on $X$. By an
$\mcal{A}$-bimodule we mean a sheaf $\mcal{M}$ of $\k$-modules 
on $X$, together with a left $\mcal{A}$-module structure and a 
right $\mcal{A}$-module structure that commute with each other. 
In other words $\mcal{M}$ is a module over the sheaf of rings
$\mcal{A} \otimes_{\k_X} \mcal{A}^{\mrm{op}}$,
where $\k_X$ is the constant sheaf $\k$ on $X$.
An $\mcal{A}$-ring is a sheaf $\mcal{B}$ of rings on $X$
together with a ring homomorphism 
$\mcal{A} \to \mcal{B}$. 
Note that $\mcal{B}$ is an $\mcal{A}$-bimodule.

\begin{dfn} \label{dfn3.1}
\begin{enumerate}
\item A {\em ringed scheme} over $\k$ is a pair 
$(X, \mcal{A})$ consisting of a $\k$-scheme $X$ and an 
$\mcal{O}_{X}$-ring $\mcal{A}$.
\item We say $\mcal{A}$ is a {\em quasi-coherent 
$\mcal{O}_{X}$-ring}, and the pair $(X, \mcal{A})$ is a 
{\em quasi-coherent ringed scheme}, if the $\mcal{O}_{X}$-bimodule
$\mcal{A}$ is a quasi-coherent $\mcal{O}_{X}$-module on both 
sides.
\item A quasi-coherent ringed scheme $(X, \mcal{A})$ is called 
{\em separated} (resp.\ {\em affine})
if $X$ is a separated (resp.\ affine) $\k$-scheme. 
\end{enumerate}
\end{dfn}

Henceforth in this section we consider quasi-coherent ringed 
schemes.

\begin{dfn}
\label{dfn4.7}
A morphism $(f, \psi): (Y, \mcal{B}) \to (X, \mcal{A})$
of quasi-coherent ringed schemes over $\k$ is a morphism of 
schemes $f: Y \to X$, together with a homomorphism
$\psi : \mcal{A} \to f_{*} \mcal{B}$
of $\mcal{O}_{X}$-rings.
\end{dfn}

Often we will use an abbreviated form, and denote by $f^*$ both 
the homomorphism $\mcal{O}_X \to f_{*} \mcal{O}_Y$
and the homomorphism 
$\psi : \mcal{A} \to f_{*} \mcal{B}$.

Suppose we are given two quasi-coherent ringed schemes 
$(X_1, \mcal{A}_1)$ and $(X_2, \mcal{A}_2)$ over $\k$. Let us denote 
by  $X_1 \times X_2 := X_1 \times_{\k} X_2$
the usual product of schemes, and by
$\opnt{p}_{i} : X_1 \times X_2 \to X_i$ the projections. 
We obtain sheaves of $\k_{X \times Y}$-algebras
$\opnt{p}_{i}^{-1} \mcal{O}_{X_i}$
and
$\opnt{p}_{i}^{-1} \mcal{A}_i$,
and there are canonical ring homomorphisms:
\[ \begin{CD}
\opnt{p}_{1}^{-1} \mcal{O}_{X_1} \otimes 
\opnt{p}_{2}^{-1} \mcal{O}_{X_2} @>>> \mcal{O}_{X_1 \times X_2} \\
@VVV \\
\opnt{p}_{1}^{-1} \mcal{A}_1 \otimes \opnt{p}_{2}^{-1} \mcal{A}_2
\end{CD} \]
where 
``$\otimes$'' denotes ``$\otimes_{\k_{X \times Y}}$''.

\begin{dfn} \label{dfn4.1}
Let $(X_1, \mcal{A}_1)$ and $(X_2, \mcal{A}_2)$ be two quasi-coherent 
ringed schemes over $\k$. Their {\em product}
is a quasi-coherent ringed scheme 
$(X_1 \times X_2, \mcal{A}_1 \boxtimes \mcal{A}_2)$, together with a
ring homomorphism
\[ \phi: \opnt{p}_{1}^{-1} \mcal{A}_1 \otimes 
\opnt{p}_{2}^{-1} \mcal{A}_2 \to \mcal{A}_1 \boxtimes \mcal{A}_2 , \]
satisfying the conditions below. 
\begin{enumerate}
\rmitem{i} The diagram
\[  \begin{CD}
\opnt{p}_{1}^{-1} \mcal{O}_{X_1} \otimes \opnt{p}_{2}^{-1} 
\mcal{O}_{X_2} @>>> \mcal{O}_{X_1 \times X_2} \\
@VVV @VVV \\
\opnt{p}_{1}^{-1} \mcal{A}_1 \otimes \opnt{p}_{2}^{-1} 
\mcal{A}_2 @>{\phi}>> \mcal{A}_1 \boxtimes \mcal{A}_2
\end{CD} \]
commutes.
\rmitem{ii} For every pair of affine open sets 
$U_1 \subset X_1$ and $U_2 \subset X_2$ the homomorphism
\[ \Gamma(U_1, \mcal{A}_1) \otimes \Gamma(U_2, \mcal{A}_2) \to 
\Gamma(U_1 \times U_2, \mcal{A}_1 \boxtimes \mcal{A}_2) \]
induced by $\phi$ is bijective.
\end{enumerate}
\end{dfn}

We remind that a denominator set $S$ in a ring $A$ is a 
multiplicatively closed subset satisfying the left and right Ore 
conditions and the left and right torsion conditions. 
Equivalently, it means that the (left and right) ring of fractions 
$A_S = S^{-1} A = A S^{-1}$ exists. See \cite[Section 2.1]{MR} 
for more details.

\begin{lem} \label{lem4.1}
Let $(U, \mcal{A})$ and $(V, \mcal{B})$ be two affine
quasi-coherent ringed schemes, 
$f : V \to U$ a morphism of schemes and 
$\psi : \Gamma(U, \mcal{A}) \to \Gamma(V, \mcal{B})$
a homomorphism of $\Gamma(U, \mcal{O}_U)$-rings. Then there is a 
unique homomorphism of $\mcal{O}_U$-rings
$\wtil{\psi} : \mcal{A} \to f_* \mcal{B}$
such that 
$\Gamma(U, \wtil{\psi}) = \psi$.
\end{lem}

\begin{proof}
Write $C := \Gamma(U, \mcal{O}_U)$, 
$A := \Gamma(U, \mcal{A})$ and
$B := \Gamma(V, \mcal{B})$. Choose an element
$s \in C$, and define the multiplicatively closed sets 
$S := \{ s^i \}_{i \in \mbb{Z}} \subset C$, 
$\bar{S}$ the image of $S$ in $A$ and 
$\bar{T} := \psi(\bar{S}) \subset B$.
Also let $U' := \{ s \neq 0 \} \subset U$ and
$V' := \{ \psi(s) \neq 0 \} \subset V$.
According to \cite[Corollary 5.13]{YZ4}, $\bar{S}$ and
$\bar{T}$ are denominator sets in $A$ and $B$ respectively.  
By the universal property of localization there is a unique
homomorphism of $C_S$-rings
$A_{\bar{S}} \to B_{\bar{T}}$
extending $\psi$. Now by \cite[Proposition 5.6]{YZ4} we have
$A_{\bar{S}} \cong \Gamma(U', \mcal{A})$ and
$B_{\bar{T}} \cong \Gamma(V', \mcal{B}) \cong
\Gamma(U', f_* \mcal{B})$. Since the open sets $U'$, as we change 
$s$, form a basis of the topology of $U$, we obtain an 
$\mcal{O}_U$-ring homomorphism 
$\wtil{\psi} : \mcal{A} \to f_* \mcal{B}$.
\end{proof}

The lemma says that 
$(f, \wtil{\psi}) : (V, \mcal{B}) \to (U, \mcal{A})$
is a morphism of quasi-coherent ringed schemes.

\begin{lem} 
Let $(X_{1}, \mcal{A}_{1})$ and $(X_{2}, \mcal{A}_{2})$
be quasi-coherent ringed 
schemes over $\k$, and assume a product
$(X_{1} \times X_{2}, \mcal{A}_{1} \boxtimes \mcal{A}_{2})$ 
exists. Then the projections 
$\mrm{p}_{i} : X_{1} \times X_{2} \to X_{i}$
extend to morphisms of ringed schemes
\[ \wtil{\mrm{p}}_{i}: 
(X_{1} \times X_{2}, \mcal{A}_{1} \boxtimes \mcal{A}_{2})
\to (X_{i}, \mcal{A}_{i}) . \]
\end{lem}

\begin{proof} 
By symmetry it suffices to show that the first projection 
$\wtil{\opnt{p}}_1$ exists. Let us choose affine open sets 
$U_1 \subset X_1$ and $U_2 \subset X_2$. By Lemma \ref{lem4.1}
the ring homomorphism
\[ \Gamma(U_1, \mcal{A}_{1}) \to
\Gamma(U_1, \mcal{A}_{1}) \otimes \Gamma(U_2, \mcal{A}_{2})
\cong 
\Gamma(U_1 \times U_2, \mcal{A}_{1} \boxtimes \mcal{A}_{2}) \]
extends to a morphism of ringed schemes
\[ \wtil{\mrm{p}}_1|_{U_1 \times U_2} :
(U_1 \times U_2, (\mcal{A}_{1} \boxtimes \mcal{A}_{2})
|_{U_1 \times U_2}) \to
(U_1, \mcal{A}_{1}|_{U_1}) \to (X_1, \mcal{A}_1) . \]
Because the morphisms $\wtil{\mrm{p}}_1|_{U_1 \times U_2}$ 
depend functorially on the affine open rectangles 
$U_1 \times U_2$ we can glue them to obtain a global morphism 
$\wtil{\mrm{p}}_1$. 
\end{proof}

\begin{prop} \label{prop4.1}
Let $(X_{1}, \mcal{A}_{1})$, $(X_{2}, \mcal{A}_{2})$,
$(Y_{1}, \mcal{B}_{1})$ and $(Y_{2}, \mcal{B}_{2})$
be quasi-coherent ringed schemes over $\k$, and let
$\wtil{f}_{i}:  (Y_{i}, \mcal{B}_{i}) 
\to (X_{i}, \mcal{A}_{i})$
be morphisms. Assume products 
$(X_{1} \times X_{2}, \mcal{A}_{1} \boxtimes \mcal{A}_{2})$ 
and
$(Y_{1} \times Y_{2}, \mcal{B}_{1} \boxtimes \mcal{B}_{2})$ 
exist. Then there is a unique morphism 
\[ \wtil{f}_1 \times \wtil{f}_2 : 
(Y_{1} \times Y_{2}, \mcal{B}_{1} \boxtimes \mcal{B}_{2}) \to
(X_{1} \times X_{2}, \mcal{A}_{1} \boxtimes \mcal{A}_{2}) \]
such that
$\wtil{\mrm{p}}_i \circ (\wtil{f}_1 \times \wtil{f}_2) =
\wtil{f}_i \circ \wtil{\mrm{p}}_i$.
\end{prop}

\begin{proof}
Let $U_i \subset X_i$ and $V_i \subset f_1^{-1}(U_i)$
be affine open sets. The ring homomorphisms
$\wtil{f}_i^* : \Gamma(U_i, \mcal{A}_i) \to
\Gamma(V_i, \mcal{B}_i)$
give rise to a ring homomorphism
\[ \begin{aligned}
& \Gamma(U_1 \times U_2, \mcal{A}_{1} \boxtimes \mcal{A}_{2}) 
\cong \Gamma(U_1, \mcal{A}_1) \otimes \Gamma(U_2, \mcal{A}_2) \\
& \quad \xar{\wtil{f}_1^* \otimes \wtil{f}_2^*}
\Gamma(V_1, \mcal{B}_1) \otimes \Gamma(V_2, \mcal{B}_2)
\cong 
\Gamma(V_1 \times V_2, \mcal{B}_{1} \boxtimes \mcal{B}_{2}) . 
\end{aligned} \]
Using Lemma \ref{lem4.1} we obtain a unique morphism of 
ringed schemes
\[ \big( V_1 \times V_2, (\mcal{B}_{1} \boxtimes \mcal{B}_{2})
|_{V_1 \times V_2} \big) \to
\big( U_1 \times U_2, (\mcal{A}_{1} \boxtimes \mcal{A}_{2})
|_{U_1 \times U_2} \big) \]
that's compatible with the projections. By gluing we obtain
$\wtil{f}_1 \times \wtil{f}_2$.
\end{proof}

\begin{cor} \label{cor4.2}
The product 
$(X_{1} \times X_{2}, \mcal{A}_{1} \boxtimes \mcal{A}_{2})$, 
together with the projections $\wtil{\mrm{p}}_{1}$ and
$\wtil{\mrm{p}}_{1}$, is unique up to a unique isomorphism. 
\end{cor}

\begin{proof}
Take $(X_i, \mcal{A}_i) = (Y_i, \mcal{B}_i)$ in Proposition 
\ref{prop4.1}. 
\end{proof}

Observe that given a quasi-coherent $\mcal{O}_X$-ring
$\mcal{A}$ the opposite ring $\mcal{A}^{\mrm{op}}$ is also a  
quasi-coherent $\mcal{O}_X$-ring. 

\begin{dfn} \label{dfn4.2}
Let $(X, \mcal{A})$ be a quasi-coherent ringed scheme over $\k$. 
We denote by
\[ (X^2, \mcal{A}^{\mrm{e}}) := 
(X \times X, \mcal{A} \boxtimes \mcal{A}^{\mrm{op}}) , \]
the product of the quasi-coherent ringed schemes 
$(X, \mcal{A})$ and $(X, \mcal{A}^{\mrm{op}})$.
\end{dfn}

Here is an easy example of a product.

\begin{exa} \label{exa4.1}
Let $(X, \mcal{A})$ be any quasi-coherent ringed scheme. Let 
$B$ be a $\k$-algebra, which we consider as a quasi-coherent 
ringed scheme $(Y, \mcal{B}) := (\opn{Spec} \k, B)$. 
Then the product exists, and it is 
\[ (X \times Y, \mcal{A} \boxtimes \mcal{B}) =
(X, \mcal{A} \otimes B). \]
\end{exa}

The existence of products turns out to be more complicated in 
general, as we see in the next theorem. 

\begin{thm} \label{thm4.2}
Let $(X, \mcal{A})$ and $(Y, \mcal{B})$ be quasi-coherent 
ringed schemes over $\k$. The two conditions below are 
equivalent. 
\begin{enumerate}
\rmitem{i} The product $(X \times Y, \mcal{A} \boxtimes \mcal{B})$
exists. 
\rmitem{ii} Given any pair of affine open sets 
$U \subset X$ and $V \subset Y$, write 
$C := \Gamma(U \times V, \mcal{O}_{X \times Y})$, 
$A := \Gamma(U, \mcal{A})$
and 
$B := \Gamma(V, \mcal{B})$.
Let $S \subset C$ be any multiplicatively closed set. Then the 
image $\bar{S}$ of $S$ in $A \otimes B$ is a denominator set.
\end{enumerate}
\end{thm}

\begin{proof}
(i) $\Rightarrow$ (ii): Choose affine open sets $U \subset X$
and $V \subset Y$. Then 
$(\mcal{A} \boxtimes \mcal{B})|_{U \times V}$ is a quasi-coherent 
$\mcal{O}_{U \times V}$-ring satisfying
\[ \Gamma \big( U \times V, 
(\mcal{A} \boxtimes \mcal{B})|_{U \times V} \big) 
\cong A \otimes B \]
as $C$-rings. According to \cite[Corollary 5.13]{YZ4}
for any multiplicatively closed set $S \subset C$ the set
$\bar{S} \subset A \otimes B$ is a denominator set.

\medskip \noindent
(ii) $\Rightarrow$ (i):
Let $U \subset X$ and $V \subset Y$ be arbitrary affine open 
sets. According to \cite[Corollary 5.13]{YZ4}
there is a quasi-coherent $\mcal{O}_{U \times V}$-ring, call it 
$\mcal{A}|_U \boxtimes \mcal{B}|_V$, such that 
\[ \Gamma(U \times V, \mcal{A}|_U \boxtimes \mcal{B}|_V) 
\cong A \otimes B \]
as $C$-rings. We claim that 
$(U \times V, \mcal{A}|_U \boxtimes \mcal{B}|_V)$
is a product of $(U, \mcal{A}|_U)$ and 
$(V, \mcal{B}|_V)$.

By Lemma \ref{lem4.1} the ring homomorphism 
$A \to A \otimes B$ induces a morphism of ringed schemes
$(U \times V,  \mcal{A}|_U \boxtimes \mcal{B}|_V) \to
(U, \mcal{A}|_U)$, which in turn gives us a homomorphism of 
sheaves of rings
$(\mrm{p}_1^{-1} \mcal{A})|_{U \times V} \to 
\mcal{A}|_U \boxtimes \mcal{B}|_V$. 
Likewise we get
$(\mrm{p}_2^{-1} \mcal{B})|_{U \times V} \to 
\mcal{A}|_U \boxtimes \mcal{B}|_V$. 
Multiplying inside $\mcal{A}|_U \boxtimes \mcal{B}|_V$ we obtain a 
homomorphism of modules
\[ \phi|_{U \times V} : 
(\mrm{p}_1^{-1} \mcal{A})|_{U \times V} \otimes 
(\mrm{p}_2^{-1} \mcal{B})|_{U \times V} \to 
\mcal{A}|_U \boxtimes \mcal{B}|_V . \]
In order to show that this is a homomorphism of sheaves of
rings it suffices to check that for any 
point $z \in U \times V$ the images of 
$(\mrm{p}_1^{-1} \mcal{A})|_{U \times V}$ and 
$(\mrm{p}_2^{-1} \mcal{B})|_{U \times V}$ inside the stalk 
$(\mcal{A}|_U \boxtimes \mcal{B}|_V)_z$ commute. 
Let us write $x := \mrm{p}_1(z)$ and $y := \mrm{p}_2(z)$.
Then 
$(\mrm{p}_1^{-1} \mcal{A})_z \cong \mcal{A}_x$
and $(\mrm{p}_2^{-1} \mcal{B})_z \cong \mcal{B}_y$.
As explained in the proof of Lemma \ref{lem5.6} the ring
$(\mcal{A}|_U \boxtimes \mcal{B}|_V)_z$ is an Ore localization of 
$\mcal{A}_x \otimes \mcal{B}_y$, so indeed there is 
commutation. 

The formation of the homomorphism $\phi|_{U \times V}$ is 
functorial in the affine open rectangle $U \times V$, and hence we 
can glue. By construction conditions (i)-(ii) of Definition 
\ref{dfn4.1} are satisfied.
\end{proof}

\begin{exa}
If $\mcal{A}$ is an $\mcal{O}_X$-algebra (i.e.\ the image of 
$\mcal{O}_X$ is in the center $\mrm{Z}(\mcal{A})$)
and $\mcal{B}$ is an 
$\mcal{O}_Y$-algebra then condition (ii) in the theorem is 
automatically satisfied. Hence the product 
$(X \times Y, \mcal{A} \boxtimes \mcal{B})$
exists. 
\end{exa}

\begin{rem}
Our definition of quasi-coherent ringed scheme is slightly 
different than that of Lunts \cite{Lu}. While in \cite{Lu} the 
existence of a product must be stipulated, here we have the 
criterion in Theorem \ref{thm4.2}. There are even
counterexamples; indeed in \cite[Example 5.14]{YZ4} we exhibit 
an affine quasi-coherent ringed scheme $(X, \mcal{A})$ for which 
the product $(X^2, \mcal{A}^{\mrm{e}})$ does not exist.
\end{rem}

Let $(X, \mcal{A})$ be a quasi-coherent ringed scheme over $\k$. 
We denote by $\cat{Mod} \mcal{A}$ the category of 
left $\mcal{A}$-modules, and by $\msf{D}(\cat{Mod} \mcal{A})$ its 
derived category. 
The category $\cat{Mod} \mcal{A}$ is abelian and it has 
enough injectives. 
Given an injective $\mcal{A}$-module $\mcal{I}$ its restriction 
$\mcal{I}|_{U}$ to an open subset $U$ is an injective 
$\mcal{A}|_{U}$-module. 
Any complex $\mcal{N} \in \msf{D}^{+}(\cat{Mod} \mcal{A})$ has an 
injective resolution $\mcal{N} \to \mcal{I}$, namely a 
quasi-isomorphism to a bounded below complex $\mcal{I}$ of 
injective $\mcal{A}$-modules. This allows us to define the derived 
functor 
\[ \mrm{R} \mcal{H}om_{\mcal{A}}(-, -) :
\msf{D}^{-}(\cat{Mod} \mcal{A})^{\mrm{op}} \times
\msf{D}^{+}(\cat{Mod} \mcal{A}) \to
\msf{D}^+(\cat{Mod} \k_X) , \]
where $\k_X$ is the constant sheaf $\k$ on $X$.
The formula is 
$\mrm{R} \mcal{H}om_{\mcal{A}}(\mcal{M}, \mcal{N}) :=
\mcal{H}om_{\mcal{A}}(\mcal{M}, \mcal{I})$
for any injective resolution $\mcal{N} \to \mcal{I}$.
Since $\mcal{H}om_{\mcal{A}}(\mcal{M}, \mcal{I})$
is a bounded below complex of flasque sheaves it follows that 
\[ \opn{RHom}_{\mcal{A}}(\mcal{M}, \mcal{N}) \cong
\mrm{R} \Gamma(X, \mrm{R} \mcal{H}om_{\mcal{A}}
(\mcal{M}, \mcal{N})) , \]
which is a functor
\[ \opn{RHom}_{\mcal{A}}(-, -) :
\msf{D}^{-}(\cat{Mod} \mcal{A})^{\mrm{op}} \times
\msf{D}^{+}(\cat{Mod} \mcal{A}) \to
\msf{D}^+(\cat{Mod} \k) . \]

For more details regarding derived categories of sheaves 
see \cite{RD}, \cite{KS} or \cite{Bor}.

\begin{rem}
One can of course remove some boundedness restrictions using 
K-injective resolutions, but we are not going to worry about this. 
Already there are enough delicate issues regarding injective 
resolutions of quasi-coherent $\mcal{A}$-modules; see Remark
\ref{rem5.2}.
\end{rem}

\begin{dfn} 
\label{dfn3.2}
Let $(X, \mcal{A})$ be a quasi-coherent ringed scheme over 
$\k$. Suppose $X$ is noetherian, and for every affine 
open set $U \subset X$ the ring $\Gamma(U, \mcal{A})$ is 
left noetherian. Then we call $\mcal{A}$ a {\em left 
noetherian quasi-coherent $\mcal{O}_{X}$-ring}, and the pair 
$(X, \mcal{A})$ is called a {\em left noetherian quasi-coherent 
ringed scheme}. We say $\mcal{A}$ is right noetherian if 
$\mcal{A}^{\mrm{op}}$ is left noetherian. If $\mcal{A}$ is both 
left and right noetherian then the pair $(X, \mcal{A})$ is called a 
{\em noetherian quasi-coherent ringed scheme}.
\end{dfn}

When we speak of a left or right noetherian quasi-coherent 
$\mcal{O}_{X}$-ring $\mcal{A}$ we tacitly assume that $X$ itself is 
noetherian.

An $\mcal{A}$-module $\mcal{M}$ is called quasi-coherent if 
locally, on every sufficiently small open set $U$, the 
$\mcal{A}|_U$-module $\mcal{M}|_U$ has a free presentation
(i.e.\ it is the cokernel of a homomorphism between free 
$\mcal{A}|_U$-modules); 
cf.\ \cite{EGA-I}. Equivalently, $\mcal{M}$ is 
quasi-coherent as $\mcal{O}_X$-module.  If $\mcal{A}$ is 
left noetherian then $\mcal{M}$ is a coherent $\mcal{A}$-module 
if and only if it is quasi-coherent and locally finitely generated. 
We shall denote the categories of quasi-coherent 
(resp.\ coherent) $\mcal{A}$-modules by $\cat{QCoh} \mcal{A}$ 
(resp.\ $\cat{Coh} \mcal{A}$). For a ring $A$ we write 
$\cat{Mod}_{\mrm{f}} A$ for the category of finite (i.e.\ finitely
generated) modules.

\begin{prop}
\label{prop3.16}
Let $(X,\mcal{A})$ be a quasi-coherent ringed scheme, let 
$U \subset X$ be an affine open set and
$A := \Gamma(U, \mcal{A})$.
\begin{enumerate}
\item The functor $\Gamma(U, -)$
is an equivalence of categories
$\cat{QCoh} \mcal{A}|_U \to \cat{Mod} A$.
\item If $\mcal{A}$ is left noetherian then $\Gamma(U, -)$ 
restricts to an equivalence of categories
$\cat{Coh} \mcal{A}|_U \to \cat{Mod}_{\mrm{f}} A$.
\end{enumerate}
\end{prop} 

\begin{proof}
This is a slight generalization of \cite[Corollary 1.4.2 and
Theorem 1.5.1]{EGA-I}. See also \cite[Corollary II.5.5]{Ha}. 
\end{proof}

In the context of the proposition above, given an $A$-module $M$ 
we shall usually denote the corresponding 
quasi-coherent $\mcal{A}|_U$-module by
$\mcal{A}|_U \otimes_{A} M$.

\begin{rem} \label{rem5.2}
By \cite[Remark 3.2]{Ka} any injective object $\mcal{I}$ in
$\cat{QCoh} \mcal{A}$ is a flasque sheaf on $X$. However $\mcal{I}$ 
might 
not be injective in the bigger category $\cat{Mod} \mcal{A}$, and 
there are counterexamples (see \cite[Remark 7.4]{Ka}). Moreover,
we do not know if the restriction to an open set $\mcal{I}|_U$ 
is injective in $\cat{QCoh} \mcal{A}|_{U}$; the results of 
\cite{GJ} seem to indicate otherwise.
This is in contrast 
to the commutative noetherian case $\mcal{A} = \mcal{O}_{X}$,
in which $\mcal{I}|_U$ is always injective in 
$\cat{Mod} \mcal{A}|_U$;  see \cite[Proposition II.7.17]{RD}.
\end{rem}

As is customary we denote by 
$\msf{D}_{\mrm{qc}}(\cat{Mod} \mcal{A})$
(resp. $\msf{D}_{\mrm{c}}(\cat{Mod} \mcal{A})$)
the full subcategory of $\msf{D}(\cat{Mod} \mcal{A})$
consisting of complexes with quasi-coherent (resp.\ coherent)
cohomology modules. These are triangulated subcategories (the 
latter if $\mcal{A}$ is left noetherian). 

We now recall a theorem of Bernstein about equivalences of 
derived categories of $\mcal{A}$-modules. 

\begin{thm}[{\cite[Theorem VI.2.10 and Proposition VI.2.11]{Bor}}]
\label{thm5.3}
Let $(X, \mcal{A})$ be a separated quasi-coherent ringed scheme. 
\begin{enumerate}
\item The inclusion functor 
$\msf{D}^{\mrm{b}}(\cat{QCoh} \mcal{A}) \to
\msf{D}_{\mrm{qc}}^{\mrm{b}}(\cat{Mod} \mcal{A})$
is an equivalence.
\item If in addition $\mcal{A}$ is left noetherian then 
the inclusion functor
$\msf{D}^{\mrm{b}}(\cat{Coh} \mcal{A}) \to
\msf{D}_{\mrm{c}}^{\mrm{b}}(\cat{Mod} \mcal{A})$
is an equivalence.
\end{enumerate}
\end{thm}

\begin{cor} \label{cor5.5} 
Assume $(X, \mcal{A})$ is a quasi-coherent ringed scheme
and $U \subset X$ is an affine open set. Write 
$A := \Gamma(U, \mcal{A})$. Then 
\[ \mrm{R} \Gamma(U, -) : 
\msf{D}_{\mrm{qc}}^{\mrm{b}}(\cat{Mod} \mcal{A}|_{U}) \to 
\msf{D}^{\mrm{b}}(\cat{Mod} A) \]
is an equivalence with inverse 
$M \mapsto \mcal{A}|_U \otimes_{A} M$. 
If $\mcal{A}$ is left noetherian then we get an equivalence
\[ \mrm{R} \Gamma(U, -) : 
\msf{D}_{\mrm{c}}^{\mrm{b}}(\cat{Mod} \mcal{A}|_{U}) \to 
\msf{D}^{\mrm{b}}_{\mrm{f}}(\cat{Mod} A) . \]
\end{cor}

\begin{proof}
By Proposition \ref{prop3.16} any quasi-coherent 
$\mcal{O}_{U}$-module is acyclic for the functor
$\Gamma(U, -)$. Hence if 
$\mcal{M} \in \msf{D}^{\mrm{b}}(\cat{QCoh} \mcal{A}|_U)$, and if 
$\mcal{M} \to \mcal{I}$ is a quasi-isomorphism with $\mcal{I}$ a 
bounded below complex of injective $\mcal{A}|_U$-modules (not 
necessarily quasi-coherent), then 
$\Gamma(U, \mcal{M}) \to \Gamma(U, \mcal{I})$
is a quasi-isomorphism. It follows that the composed functor
\[ \msf{D}^{\mrm{b}}(\cat{QCoh} \mcal{A}|_U) \to
\msf{D}^{\mrm{b}}_{\mrm{qc}}(\cat{Mod} \mcal{A}|_U) 
\xar{\mrm{R} \Gamma(U, -)} 
\msf{D}^{\mrm{b}}(\cat{Mod} A) \]
is an equivalence. Now use Theorem \ref{thm5.3}. 
\end{proof}

\begin{lem} \label{lem5.4}
Let $(X, \mcal{A})$ be a quasi-coherent ringed scheme and let 
$V \subset U$ be two affine open sets. Let 
$A := \Gamma(U, \mcal{A})$ and $A' := \Gamma(V, \mcal{A})$. Then 
the diagram 
\[ \begin{CD}
\msf{D}^{+}_{\mrm{qc}}(\cat{Mod} \mcal{A}|_U)
@>{\mrm{R} \Gamma(U, -)}>> 
\msf{D}^{+}(\cat{Mod} A) \\
@V{\mrm{rest}}VV @V{A' \otimes_A -}VV \\
\msf{D}^{+}_{\mrm{qc}}(\cat{Mod} \mcal{A}|_V)
@>{\mrm{R} \Gamma(V, -)}>> 
\msf{D}^{+}(\cat{Mod} A')
\end{CD} \]
is commutative.
\end{lem}

\begin{proof}
First we note that $A'$ is a flat $A$-module on both sides (see 
\cite[Proposition 5.6]{YZ4}. 
Given $\mcal{M} \in \msf{D}^{+}_{\mrm{qc}}(\cat{Mod} \mcal{A}|_U)$
take a resolution $\mcal{M} \to \mcal{I}$ where $\mcal{I}$ is a 
bounded below complex of injective $\mcal{A}|_U$-modules. 
Then 
$\mcal{M}|_V \to \mcal{I}|_V$ is an injective resolution. We get a 
natural morphism
\[ \mrm{R} \Gamma(U, \mcal{M}) = \Gamma(U, \mcal{I}) \to
\Gamma(V, \mcal{I}) = \mrm{R} \Gamma(V, \mcal{M}) , \]
and hence a morphism
\[ A' \otimes_A \mrm{R} \Gamma(U, \mcal{M}) \to 
\mrm{R} \Gamma(V, \mcal{M}) . \]
To show the latter is an isomorphism it suffices to check for a 
single quasi-coherent $\mcal{A}|_U$-module $\mcal{M}$ -- since 
these are way-out functors, cf.\ \cite[Section I.7]{RD}.
But for such $\mcal{M}$ we have
$\mrm{R} \Gamma(U, \mcal{M}) = \Gamma(U, \mcal{M})$
and 
$\mrm{R} \Gamma(V, \mcal{M}) = \Gamma(V, \mcal{M})$,
so \cite[Proposition 5.6]{YZ4} applies. 
\end{proof}

\section{Dualizing Complexes over Ringed Schemes}
\label{sec2}

Derived functors in the noncommutative situation are inherently 
more complicated than in the commutative situation. Therefore we 
begin this section with some technical results that are needed  
to state our definition of dualizing complexes over quasi-coherent 
ringed schemes.

\begin{lem} \label{lem5.6}
Let $(X, \mcal{A})$ and $(Y, \mcal{B})$ be two quasi-coherent ringed
schemes over $\k$. Assume the product 
$(X \times Y, \mcal{A} \boxtimes \mcal{B})$
exists. Then:
\begin{enumerate}
\item The homomorphisms of sheaves of rings
\[ \opnt{p}_{1}^{-1} \mcal{A} \to 
(\opnt{p}_{1}^{-1} \mcal{A}) \otimes
(\opnt{p}_{2}^{-1} \mcal{B}) \to \mcal{A} \boxtimes \mcal{B} \]
are flat on both sides.
\item An injective 
$\mcal{A} \boxtimes \mcal{B}$ -module is also an 
injective $\opnt{p}_{1}^{-1} \mcal{A}$ -module.
\item Given a $\opnt{p}_{1}^{-1} \mcal{A}$ -module
$\mcal{M}$ and an injective 
$\mcal{A} \boxtimes \mcal{B}$ -module $\mcal{I}$, the 
$\opnt{p}_{2}^{-1} \mcal{B}$ -module
$\mcal{H}om_{\opnt{p}_{1}^{-1} \mcal{A}}(\mcal{M}, \mcal{I})$
is injective.
\end{enumerate}
These statements hold also if we exchange
$\opnt{p}_{1}^{-1} \mcal{A}$ with $\opnt{p}_{2}^{-1} \mcal{B}$, 
or $\opnt{p}_{1}^{-1} \mcal{A}$ with 
$\opnt{p}_{1}^{-1} \mcal{A}^{\mrm{op}}$, etc. 
\end{lem}

\begin{proof} 
(1) Flatness can be checked on stalks at points of $X \times Y$.
Take an arbitrary point $z \in X \times Y$, and let
$x := \opnt{p}_{1}(z) \in X$ and 
$y := \opnt{p}_{2}(z) \in Y$.
The stalks satisfy
$(\opnt{p}_{1}^{-1} \mcal{A})_z \cong \mcal{A}_x$,
$(\opnt{p}_{2}^{-1} \mcal{B})_z \cong \mcal{B}_y$
and
\[ \bigl( (\opnt{p}_{1}^{-1} \mcal{A}) \otimes
(\opnt{p}_{2}^{-1} \mcal{B}) \bigr)_{z} \cong
(\opnt{p}_{1}^{-1} \mcal{A})_{z} \otimes 
(\opnt{p}_{2}^{-1} \mcal{B})_{z}  \cong
\mcal{A}_{x} \otimes \mcal{B}_{y} . \]
Hence the homomorphism
\[ (\opnt{p}_{1}^{-1} \mcal{A})_{z} \to 
\bigl( (\opnt{p}_{1}^{-1} \mcal{A}) \otimes
(\opnt{p}_{2}^{-1} \mcal{B}) \bigr)_{z} \]
is flat on both sides.

Choose affine open sets $U \subset X$ and $V \subset Y$
such that $z \in U \times V$. By definition 
\[ \Gamma(U \times V, \mcal{A} \boxtimes \mcal{B}) 
\cong \Gamma(U, \mcal{A}) \otimes \Gamma(V, \mcal{B}) . \]
Now the stalks $\mcal{A}_{x}$, $\mcal{B}_{y}$ and
$(\mcal{A} \boxtimes \mcal{B})_{z}$ are Ore localizations of 
$\Gamma(U, \mcal{A})$, $\Gamma(V, \mcal{B})$ and
$\Gamma(U \times V, \mcal{A} \boxtimes \mcal{B})$ 
respectively, for suitable denominator sets. Hence 
$(\mcal{A} \boxtimes \mcal{B})_{z}$ is an Ore localization of  
$(\mcal{A}_{x} \otimes \mcal{B}_{y})$, so it's flat.

\medskip \noindent
(2) By part (1) the functor 
$\mcal{M} \mapsto (\mcal{A} \boxtimes \mcal{B})
\otimes_{\opnt{p}_{1}^{-1} \mcal{A}} \mcal{M}$
from $\cat{Mod} \opnt{p}_{1}^{-1} \mcal{A}$ to
$\cat{Mod} \mcal{A} \boxtimes \mcal{B}$
is exact. Hence for an injective $\mcal{A} \boxtimes \mcal{B}$
-module $\mcal{I}$ the functor 
\[ \mcal{M} \mapsto 
\opn{Hom}_{\opnt{p}_{1}^{-1} \mcal{A}}(\mcal{M}, \mcal{I}) \cong
\opn{Hom}_{\mcal{A} \boxtimes \mcal{B}}
\big( (\mcal{A} \boxtimes \mcal{B})
\otimes_{\opnt{p}_{1}^{-1} \mcal{A}} \mcal{M}, \mcal{I} \big) \]
is exact, implying that $\mcal{I}$ is an injective 
$\opnt{p}_{1}^{-1} \mcal{A}$ -module.

\medskip \noindent
(3) Similar to (2).
\end{proof}

\begin{lem} \label{lem5.7}
Let $(X, \mcal{A})$ and $(Y, \mcal{B})$ be quasi-coherent ringed 
schemes over $\k$. Assume that the product 
$(X \times Y, \mcal{A} \boxtimes \mcal{B}^{\mrm{op}})$ exists. Let
$\mcal{R} \in \msf{D}^{+}
(\cat{Mod}\,  \mcal{A} \boxtimes \mcal{B}^{\mrm{op}})$
be some complex. 
\begin{enumerate}
\item There are functors
\[ \begin{aligned}
\mrm{D} & : 
\msf{D}^{\mrm{b}}(\cat{Mod} \mcal{A})^{\mrm{op}} \to 
\msf{D}(\cat{Mod} \mcal{B}^{\mrm{op}}) , \\
\mrm{D}^{\mrm{op}} & : \msf{D}^{\mrm{b}}
(\cat{Mod} \mcal{B}^{\mrm{op}})^{\mrm{op}} \to 
\msf{D}(\cat{Mod} \mcal{A}) 
\end{aligned} \]
defined by the formulas
\[ \begin{aligned}
\mrm{D} \mcal{M} & := \mrm{R} \opnt{p}_{2 *}
\mrm{R} \mcal{H}om_{\opnt{p}_{1}^{-1} \mcal{A}}
(\opnt{p}_{1}^{-1} \mcal{M}, \mcal{R}) , \\
\mrm{D}^{\mrm{op}} \mcal{N} & := 
\mrm{R} \opnt{p}_{1 *}
\mrm{R} \mcal{H}om_{\opnt{p}_{2}^{-1} \mcal{B}^{\mrm{op}}}
(\opnt{p}_{2}^{-1} \mcal{N}, \mcal{R}) .
\end{aligned} \]
\item If $\mrm{H}^{i} \mrm{D} \mcal{M} = 0$ for $i \gg 0$ then 
$\mrm{D}^{\mrm{op}} \mrm{D} \mcal{M}$ is well defined, and there is 
a morphism
$\mcal{M} \to \mrm{D}^{\mrm{op}} \mrm{D} \mcal{M}$
in $\msf{D}(\cat{Mod} \mcal{A})$, which is functorial in 
$\mcal{M}$. Similarly for $\mcal{N}$ and 
$\mrm{D} \mrm{D}^{\mrm{op}} \mcal{N}$.
\end{enumerate}
\end{lem}

\begin{proof}
(1) Choose a quasi-isomorphism $\mcal{R} \to \mcal{I}$ 
with $\mcal{I}$ a bounded below complex of injective 
$\mcal{A} \boxtimes \mcal{B}^{\mrm{op}}$ -modules.
By Lemma \ref{lem5.6}(2) $\mcal{I}$ is a complex of injective 
$\opnt{p}_{1}^{-1} \mcal{A}$ -modules, and we get 
\[ \mrm{R} \mcal{H}om_{\opnt{p}_{1}^{-1} \mcal{A}}
(\opnt{p}_{1}^{-1} \mcal{M}, \mcal{R}) =
\mcal{H}om_{\opnt{p}_{1}^{-1} \mcal{A}}
(\opnt{p}_{1}^{-1} \mcal{M}, \mcal{I}) \in
\msf{D}(\cat{Mod} \opnt{p}_{2}^{-1} \mcal{B}^{\mrm{op}}) . \]
Because 
$\mcal{H}om_{\opnt{p}_{1}^{-1} \mcal{A}}
(\opnt{p}_{1}^{-1} \mcal{M}, \mcal{I})$
is a bounded below complex of injective 
$\opnt{p}_{2}^{-1} \mcal{B}^{\mrm{op}}$ -modules
(see Lemma \ref{lem5.6}(3)) we get
\[ \mrm{R} \opnt{p}_{2 *} \mrm{R}
\mcal{H}om_{\opnt{p}_{1}^{-1} \mcal{A}}
(\opnt{p}_{1}^{-1} \mcal{M}, \mcal{R}) = 
\opnt{p}_{2 *} \mcal{H}om_{\opnt{p}_{1}^{-1} \mcal{A}}
(\opnt{p}_{1}^{-1} \mcal{M}, \mcal{I}) 
\in \msf{D}(\cat{Mod} \mcal{B}^{\mrm{op}}) . \]
Likewise for $\mrm{D}^{\mrm{op}}$. 

\medskip \noindent
(2) Let $\mcal{R} \to \mcal{I}$ be as above.
If $\mrm{H}^{i} \mrm{D} \mcal{M} = 0$ for $i \gg 0$ then 
$\mrm{D} \mcal{M} \cong \mcal{L}$ in 
$\msf{D}(\cat{Mod} \mcal{B}^{\mrm{op}})$,
where 
$\mcal{L} \in \msf{D}^{\mrm{b}}(\cat{Mod} \mcal{B}^{\mrm{op}})$
is a truncation of $\mrm{D} \mcal{M}$. So
$\mrm{D}^{\mrm{op}} \mrm{D} \mcal{M} := \mrm{D}^{\mrm{op}} 
\mcal{L}$
is well defined. Adjunction gives a morphism
\[ \opnt{p}_{2}^{-1} \mcal{L} \cong
\opnt{p}_{2}^{-1} \opnt{p}_{2 *} 
\mcal{H}om_{\opnt{p}_{1}^{-1} \mcal{A}}
(\opnt{p}_{1}^{-1} \mcal{M}, \mcal{I}) \to
\mcal{H}om_{\opnt{p}_{1}^{-1} \mcal{A}}
(\opnt{p}_{1}^{-1} \mcal{M}, \mcal{I}) \]
in $\msf{D}(\cat{Mod}\, \opnt{p}_{1}^{-1} \mcal{A})$.
Applying 
$\mcal{H}om_{\opnt{p}_{2}^{-1} \mcal{B}^{\mrm{op}}}(-, \mcal{I})$
and a second adjunction we get a morphism
\[ \opnt{p}_{1}^{-1} \mcal{M} \to 
\mcal{H}om_{\opnt{p}_{2}^{-1} \mcal{B}^{\mrm{op}}}
\bigl( \mcal{H}om_{\opnt{p}_{1}^{-1} \mcal{A}}
(\opnt{p}_{1}^{-1} \mcal{M}, \mcal{I}), \mcal{I} \bigr) 
\to \mcal{H}om_{\opnt{p}_{2}^{-1} \mcal{B}^{\mrm{op}}}
(\opnt{p}_{2}^{-1} \mcal{L}, \mcal{I}) . \]
Finally a third adjunction gives a morphism
\[ \mcal{M} \to \opnt{p}_{1 *} \opnt{p}_{1}^{-1} \mcal{M}
\to \opnt{p}_{1 *} 
\mcal{H}om_{\opnt{p}_{2}^{-1} \mcal{B}^{\mrm{op}}}
(\opnt{p}_{2}^{-1} \mcal{L}, \mcal{I}) =
\mrm{D}^{\mrm{op}} \mcal{L} . \]
All morphisms occurring above (including the truncation 
$\mrm{D} \mcal{M} \cong \mcal{L}$) are functorial in 
$\mcal{M}$. 
\end{proof}

\begin{dfn} 
\label{dfn5.8}
Let $(X, \mcal{A})$ and $(Y, \mcal{B})$ be separated
quasi-coherent ringed 
schemes over $\k$. Assume $\mcal{A}$ and $\mcal{B}^{\mrm{op}}$ 
are left noetherian, and the product 
$(X \times Y, \mcal{A} \boxtimes \mcal{B}^{\mrm{op}})$
exists. A complex
$\mcal{R} \in \msf{D}^{\mrm{b}}_{\mrm{qc}}
(\cat{Mod}\, \mcal{A} \boxtimes \mcal{B}^{\mrm{op}})$
is called a {\em dualizing complex over $(\mcal{A}, \mcal{B})$} 
if conditions \tup{(i)-(iii)} below hold for the functors 
$\mrm{D}$ and $\mrm{D}^{\mrm{op}}$ from Lemma \ref{lem5.7}(1).
\begin{enumerate}
\rmitem{i} The functors $\mrm{D}$ and $\mrm{D}^{\mrm{op}}$
have finite cohomological dimensions when restricted to 
$\cat{Coh} \mcal{A}$
and
$\cat{Coh} \mcal{B}^{\mrm{op}}$
respectively.
\rmitem{ii} The functor $\mrm{D}$ sends $\cat{Coh} \mcal{A}$
into
$\msf{D}_{\mrm{c}} (\cat{Mod} \mcal{B}^{\mrm{op}})$,
and the functor $\mrm{D}^{\mrm{op}}$ sends 
$\cat{Coh} \mcal{B}^{\mrm{op}}$ into
$\msf{D}_{\mrm{c}} (\cat{Mod} \mcal{A})$.
\rmitem{iii} The adjunction morphisms
$\bsym{1} \to \mrm{D}^{\mrm{op}} \mrm{D}$
in $\msf{D}^{\mrm{b}}_{\mrm{c}} (\cat{Mod} \mcal{A})$, and
$\bsym{1} \to \mrm{D} \mrm{D}^{\mrm{op}}$
in $\msf{D}^{\mrm{b}}_{\mrm{c}} (\cat{Mod} \mcal{B}^{\mrm{op}})$,
are both isomorphisms.
\end{enumerate}
\end{dfn}

Note that conditions (i)-(ii) imply that there are well defined 
functors 
\[ \begin{aligned}
\mrm{D} & : 
\msf{D}^{\mrm{b}}_{\mrm{c}}(\cat{Mod} \mcal{A})^{\mrm{op}} 
\to \msf{D}^{\mrm{b}}_{\mrm{c}} (\cat{Mod} \mcal{B}^{\mrm{op}})
, \\
\mrm{D}^{\mrm{op}} & : 
\msf{D}^{\mrm{b}}_{\mrm{c}}(\cat{Mod} \mcal{B}^{\mrm{op}})^{\mrm{op}}
\to \msf{D}^{\mrm{b}}_{\mrm{c}} (\cat{Mod} \mcal{A}) .
\end{aligned} \]
By Lemma \ref{lem5.7}(2) we have the adjunction morphisms appearing
in condition (iii). 

\begin{conv} \label{conv5.1}
When we talk about a dualizing complex $\mcal{R}$ over
$(\mcal{A}, \mcal{B})$ we tacitly assume that we are given 
separated
quasi-coherent ringed $\k$-schemes $(X, \mcal{A})$ and 
$(Y, \mcal{B})$, that $\mcal{A}$ and $\mcal{B}^{\mrm{op}}$
are left noetherian, and that the product 
$(X \times Y, \mcal{A} \boxtimes \mcal{B}^{\mrm{op}})$
exists. 
\end{conv}

\begin{dfn} \label{dfn7.2}
Let $(X, \mcal{A})$ be a separated noetherian quasi-coherent 
ringed scheme over $\k$. Suppose the product 
$(X^2, \mcal{A}^{\mrm{e}}) = 
(X \times X, \mcal{A} \boxtimes \mcal{A}^{\mrm{op}})$ 
exists. If 
$\mcal{R} \in \msf{D}^{\mrm{b}}_{\mrm{qc}}
(\cat{Mod} \mcal{A}^{\mrm{e}})$
is a dualizing complex over $(\mcal{A}, \mcal{A})$ 
then we say $\mcal{R}$ is a {\em dualizing complex over $\mcal{A}$}.
\end{dfn}

\begin{exa}
Suppose $X$ is a separated finite type $\k$-scheme 
and $\mcal{R}$ is a 
dualizing complex over $X$ in the sense of \cite{RD}. Let 
$\Delta : X \to X^2$ be the diagonal embedding. Then 
$\Delta_* \mcal{R}$ is a dualizing complex over $\mcal{O}_X$
in the noncommutative sense, i.e.\ in the sense of Definition 
\ref{dfn7.2}.
\end{exa}

Let us now recall the definition of dualizing complex over 
a pair of rings from \cite{Ye1} and \cite{YZ2}. 

\begin{dfn} \label{dfn2.1}
Let $A$ be a left noetherian $\k$-algebra and $B$ a right noetherian
$\k$-algebra. A complex 
$R \in \msf{D}^{\mrm{b}}(\cat{Mod}\, A \otimes B^{\mrm{op}})$ 
is called a {\em dualizing complex over $(A,B)$} if 
it satisfies the following three conditions:
\begin{enumerate}
\rmitem{i} $R$ has finite injective dimension over $A$ and over 
$B^{\mrm{op}}$.
\rmitem{ii} $R$ has finitely generated cohomology modules over $A$ 
and over $B^{\mrm{op}}$.
\rmitem{iii} The canonical morphisms 
$B \to \opn{RHom}_A(R,R)$ in 
$\msf{D}(\cat{Mod} B^{\mrm{e}})$, and 
$A \to \opn{RHom}_{B^{\mrm{op}}}(R,R)$ 
in $\msf{D}(\cat{Mod} A^{\mrm{e}})$, are both isomorphisms.
\end{enumerate}
In the case $A = B$ we say $R$ is a {\em dualizing complex over} 
$A$.
\end{dfn}

\begin{exa} \label{exa4.2}
In the situation of Definition \ref{dfn2.1} above let
$R \in \msf{D}^{\mrm{b}}(\cat{Mod}\, A \otimes B^{\mrm{op}})$.
Consider the quasi-coherent ringed schemes
$(\opn{Spec} \k, A)$ and $(\opn{Spec} \k, B)$. 
The duality functors of Lemma \ref{lem5.7} in this case are simply
$\mrm{D} = \opn{RHom}_A(-, R)$ and 
$\mrm{D}^{\mrm{op}} = \opn{RHom}_{A^{\mrm{op}}}(-, R)$.
In view of \cite[Proposition 1.3]{YZ2} we see that 
$R$ is a dualizing complex in the sense of Definition \ref{dfn5.8} 
if and only if it satisfies the conditions 
of Definition \ref{dfn2.1}. 
\end{exa}

In Theorem \ref{thm7.2} 
we show that the above is true for any pair of 
affine quasi-coherent ringed schemes.

We end this section with a couple of examples, that in fact 
digress from the main direction of our paper, yet are quite 
interesting on their own. The first is a contravariant 
Fourier-Mukai transform on elliptic curves.

\begin{prop} \label{prop5.3}
Let $X$ be an elliptic curve over $\k := \mbb{C}$ 
and let $\mcal{R}$ be the Poincar\'e bundle on $X^2$. Then 
$(X^2, \mcal{O}_{X^2})$ is the product of the quasi-coherent 
ringed scheme $(X, \mcal{O}_X)$ with itself, and
$\mcal{R} \in \msf{D}^{\mrm{b}}_{\mrm{c}}(\cat{Mod} \mcal{O}_{X^2})$
is a dualizing complex over $\mcal{O}_X$
in the sense of Definition \tup{\ref{dfn7.2}}.
\end{prop}

\begin{proof}
Conditions (i)-(ii) of Definition \ref{dfn5.8}
are clearly verified in this case, since $X^2$ is smooth and 
proper. To check 
(iii) it suffices to prove that the adjunction morphisms 
$\phi_1: \mcal{O}_{X} \to \mrm{D}^{\mrm{op}} \mrm{D} \mcal{O}_{X}$
and
$\phi_2: \mcal{O}_{X} \to \mrm{D} \mrm{D}^{\mrm{op}} \mcal{O}_{X}$
are isomorphisms. 

We begin by showing that $\phi_1$ and $\phi_2$ are nonzero. 
Going over the adjunctions in the proof of Lemma \ref{lem5.7} we 
see that 
$\mrm{D} \mcal{O}_{X} \cong \mrm{R} \opnt{p}_{2 *} \mcal{R}$
and 
\[ \mrm{Hom}_{\msf{D}(\cat{Mod} \mcal{O}_{X})}
(\mcal{O}_{X}, \mrm{D}^{\mrm{op}} \mrm{D} \mcal{O}_{X}) 
\cong \mrm{Hom}_{\msf{D}(\cat{Mod} \mcal{O}_{X})}
(\mrm{R} \opnt{p}_{2 *} \mcal{R}, 
\mrm{R} \opnt{p}_{2 *} \mcal{R}) . \]
Under the second isomorphism  
$\phi_1$ corresponds to the identity morphism of 
$\mrm{R} \opnt{p}_{2 *} \mcal{R}$, 
which is nonzero. Similarly $\phi_2$ is shown to be nonzero.

Next for any 
$\mcal{M} \in \msf{D}^{\mrm{b}}_{\mrm{c}}(\cat{Mod} \mcal{O}_X)$
define
\[ \begin{aligned}
F \mcal{M} & := \mrm{R} \opnt{p}_{2 *} 
(\opnt{p}_{1}^* \mcal{M} \otimes^{\mrm{L}}_{\mcal{O}_{X^2}}
\mcal{R}) , \\
F^{\mrm{op}} \mcal{M} & := \mrm{R} \opnt{p}_{1 *} 
(\opnt{p}_{2}^* \mcal{M} \otimes^{\mrm{L}}_{\mcal{O}_{X^2}}
\mcal{R}) , \\
E \mcal{M} & := \mrm{R} \mcal{H}om_{\mcal{O}_{X}} 
(\mcal{M}, \mcal{O}_X) . 
\end{aligned} \]
We then have
\[ \begin{aligned}
\mrm{D} \mcal{M} & = 
\mrm{R} \opnt{p}_{2 *} \mrm{R}
\mcal{H}om_{\mcal{O}_{X^2}}
(\opnt{p}_{1}^{*} \mcal{M}, \mcal{R}) \\
& \cong \mrm{R} \opnt{p}_{2 *} 
\bigl( (\opnt{p}_{1}^{*} \mrm{R} \mcal{H}om_{\mcal{O}_X}
(\mcal{M}, \mcal{O}_X)) 
\otimes^{\mrm{L}}_{\mcal{O}_{X^2}} 
\mcal{R} \bigr) = F E \mcal{M} . 
\end{aligned} \]
Likewise we get
$\mrm{D}^{\mrm{op}} \mcal{M} \cong F^{\mrm{op}} E \mcal{M}$.

Let $O \in X$ be the zero for the group structure and
$Z := \{ O \}_{\mrm{red}}$. According to \cite{Mu}
the Fourier-Mukai transform satisfies 
$F \mcal{O}_{X} \cong F^{\mrm{op}} \mcal{O}_{X} \cong 
\mcal{O}_{Z}[-1]$
and 
$F \mcal{O}_{Z} \cong F^{\mrm{op}} \mcal{O}_{Z} \cong 
\mcal{O}_{X}$. 
Because
$E \mcal{O}_{X} \cong \mcal{O}_{X}$ and
$E \mcal{O}_{Z} \cong \mcal{O}_{Z}[-1]$
we obtain 
\[ \mrm{D} \mrm{D}^{\mrm{op}} \mcal{O}_{X} \cong  
\mrm{D}^{\mrm{op}} \mrm{D} \mcal{O}_{X} \cong \mcal{O}_{X} . \]
Finally we use the facts that
$\opn{Hom}_{\msf{D}(\cat{Mod} \mcal{O}_{X})}(\mcal{O}_{X}, 
\mcal{O}_{X}) \cong \mbb{C}$
and that $\phi_1$ and $\phi_2$ are nonzero
to deduce that $\phi_1$ and $\phi_2$ are isomorphisms.
\end{proof}

The previous result can most likely be extended to higher dimensional 
abelian varieties. The next example is based on a result of 
Beilinson \cite{Be}.

\begin{prop} \label{prop5.4}
Consider projective space $X := \mbf{P}^n_{\k}$. Let 
$\mcal{E} := \boplus_{i = 0}^n \mcal{O}_{X}(i)$
and
$B := \opn{End}_{\mcal{O}_{X}}(\mcal{E})^{\mrm{op}}$.
Then 
$\mcal{E} \in \cat{Mod}\, \mcal{O}_X \otimes B^{\mrm{op}}$
is a dualizing complex over $(\mcal{O}_X, B)$.
\end{prop}

\begin{proof}
One has
\[ \mrm{D} \mcal{M} = \mrm{R} \opnt{p}_{2 *} 
\mrm{R} \mcal{H}om_{\mcal{O}_{X}}(\mcal{M}, \mcal{E}) 
\cong \opn{RHom}_{\mcal{O}_{X}}(\mcal{M}, \mcal{E}) \]
for every 
$\mcal{M} \in 
\msf{D}^{\mrm{b}}_{\mrm{c}}(\cat{Mod} \mcal{O}_{X})$.
Since $X$ is smooth and proper over $\k$, this shows that $\mrm{D}$ 
satisfies conditions (i) and (ii) of Definition \ref{dfn5.8}. 

On the other hand  
\[ \mrm{D}^{\mrm{op}} N = \mrm{R} \opnt{p}_{1 *} 
\mrm{R} \mcal{H}om_{\opnt{p}_{2}^{-1} B^{\mrm{op}}}
(\opnt{p}_{2}^{-1} N, \mcal{E}) \]
for every 
$N \in \msf{D}^{\mrm{b}}_{\mrm{f}}(\cat{Mod} B^{\mrm{op}})$.
It is known that $B$ is a finite $\k$-algebra of global 
dimension $n$. So any 
$N \in \cat{Mod}_{\mrm{f}} B^{\mrm{op}}$ 
has a resolution $P \to N$ with each $P^i$ a 
finite projective $B^{\mrm{op}}$-module, and $P^i = 0$ unless
$-n \leq i \leq 0$. Because 
$\mrm{D}^{\mrm{op}} N \cong 
\mcal{H}om_{\opnt{p}_{2}^{-1} B^{\mrm{op}}}
(\opnt{p}_{2}^{-1} P, \mcal{E})$
we see that the functor $\mrm{D}^{\mrm{op}}$ 
satisfies conditions (i) and (ii). 
 
It remains to verify that the adjunctions
$\mcal{M} \to \mrm{D}^{\mrm{op}} \mrm{D} \mcal{M}$
and
$N \to \mrm{D} \mrm{D}^{\mrm{op}} N$
are isomorphisms.

Choose an injective resolution
$\mcal{E} \to \mcal{J}$ over $\mcal{O}_{X} \otimes B^{\mrm{op}}$.
We then have 
$\mrm{D} \mcal{M} \cong \opn{Hom}_{\mcal{O}_{X}}
(\mcal{M}, \mcal{J})$
and
$\mrm{D}^{\mrm{op}} N \cong 
\mcal{H}om_{\opnt{p}_{2}^{-1} B^{\mrm{op}}}
(\opnt{p}_{2}^{-1} N, \mcal{J})$.
Because $B$ and $\mcal{E}$ generate the categories
$\msf{D}^{\mrm{b}}_{\mrm{f}}(\cat{Mod} B^{\mrm{op}})$
and 
$\msf{D}^{\mrm{b}}_{\mrm{c}}(\cat{Mod} \mcal{O}_{X})$
respectively, it suffices to check that 
$\phi_1: \mcal{E} \to \mrm{D}^{\mrm{op}} \mrm{D} \mcal{E}$
and
$\phi_2: B \to \mrm{D} \mrm{D}^{\mrm{op}} B$
are isomorphisms. Now the fact that 
$\opn{Ext}^i_{\mcal{O}_{X}}(\mcal{E}, \mcal{E}) = 0$
for all $i \neq 0$ implies that
\[ B \to \opn{Hom}_{\mcal{O}_{X}}(\mcal{J}, \mcal{J}) \to
\opn{Hom}_{\mcal{O}_{X}}(\mcal{E}, \mcal{J}) \]
are quasi-isomorphisms.  Therefore
\[ \mcal{E} 
\xar{\phi_1} \mcal{H}om_{\opnt{p}_{2}^{-1} B^{\mrm{op}}}
\bigl(\opnt{p}_{2}^{-1} 
\opn{Hom}_{\mcal{O}_{X}}(\mcal{E}, \mcal{J}),
\mcal{J} \bigr) \to \mcal{J} \]
and
\[ B \xar{\phi_2} \opn{Hom}_{\mcal{O}_{X}}
\bigl( \mcal{H}om_{\opnt{p}_{2}^{-1} B^{\mrm{op}}}
(\opnt{p}_{2}^{-1} B, \mcal{J}), \mcal{J} \bigr) 
\to \opn{Hom}_{\mcal{O}_{X}}(\mcal{J}, \mcal{J}) \]
are quasi-isomorphisms. 
\end{proof}

The crucial fact in Proposition \ref{prop5.4} is that 
$\mcal{O}_{X}, \mcal{O}_{X}(1), \ldots, \mcal{O}_{X}(n)$
is an exceptional sequence in the sense of \cite{BO}. 
This result can be extended to other smooth 
complete varieties that admit exceptional sequences.

\section{Rigid Dualizing Complexes}

We are mainly interested in dualizing complexes on ringed
schemes that have a local behavior -- as opposed to, say, the 
dualizing complexes occurring in Propositions \ref{prop5.3} and 
\ref{prop5.4}.

\begin{dfn}
\label{dfn7.1}
Let $(X, \mcal{A})$ and $(Y, \mcal{B})$ be separated 
quasi-coherent ringed schemes over $\k$, and let
$\mcal{R} \in \msf{D}^{\mrm{b}}_{\mrm{qc}}
(\cat{Mod}\, \mcal{A} \boxtimes \mcal{B}^{\mrm{op}})$
be a dualizing complex over $(\mcal{A}, \mcal{B})$. 
If the support of $\mcal{R}$ (i.e.\ the union of the supports 
of the cohomology sheaves $\mrm{H}^i \mcal{R}$)
is contained in the graph of an isomorphism of schemes
$X \iso Y$ then we call $\mcal{R}$ a {\em local dualizing complex}.
\end{dfn}

\begin{lem} \label{lem7.1}
Let $(X, \mcal{A})$ be a left noetherian quasi-coherent ringed 
scheme, $U \subset X$ an open set and $\mcal{M}$ a coherent 
$\mcal{A}|_U$-module. Then $\mcal{M}$ extends to a coherent 
$\mcal{A}$-module.
\end{lem}

\begin{proof}
Let $g: U \to X$ be the inclusion. The sheaf
$g_* \mcal{M}$ is a quasi-coherent $\mcal{O}_X$-module, hence it 
is a quasi-coherent $\mcal{A}$-module. Also  
$g_* \mcal{M} = \bigcup_\alpha \mcal{L}_{\alpha}$
where $\{ \mcal{L}_{\alpha} \}$ is the set of its coherent 
$\mcal{O}_X$-submodules (cf.\ \cite[Corollary 6.9.9]{EGA-I}).
For any $\alpha$ the image 
\[ \mcal{N}_{\alpha} := 
\opn{Im}(\mcal{A} \otimes_{\mcal{O}_{X}} \mcal{L}_{\alpha} \to
g_* \mcal{M}) \]
is a coherent $\mcal{A}$-module. Now
\[ \mcal{M} = (g_* \mcal{M})|_U = 
\bigcup_\alpha\, (\mcal{N}_{\alpha}|_U) , \]
and because $\mcal{M}$ is a noetherian object we get
$\mcal{M} = \mcal{N}_{\alpha}|_U$ for some $\alpha$.
\end{proof}

\begin{prop} 
\label{prop7.3}
Let $(X, \mcal{A})$ and $(Y, \mcal{B})$ be separated 
quasi-coherent ringed 
schemes over $\k$. Assume $\mcal{A}$ and $\mcal{B}^{\mrm{op}}$ 
are left noetherian and the product 
$(X \times Y, \mcal{A} \boxtimes \mcal{B}^{\mrm{op}})$ 
exists. Let $\mcal{R} \in \msf{D}^{\mrm{b}}_{\mrm{qc}}
(\cat{Mod}\, \mcal{A} \boxtimes \mcal{B}^{\mrm{op}})$
be a complex, and let $X = \bigcup U_{i}$ be any open 
covering of $X$. Then the following two conditions are equivalent. 
\begin{enumerate}
\rmitem{i} $\mcal{R}$ is a local dualizing complex over 
$(\mcal{A}, \mcal{B})$.
\rmitem{ii} $\mcal{R}$ is supported on the graph of some isomorphism
$f: X \iso Y$, and for every $i$, defining $V_{i} := f(U_{i})$, 
the restriction
\[ \mcal{R}|_{U_{i} \times V_{i}} \in 
\msf{D}^{\mrm{b}}_{\mrm{qc}}
\big( \cat{Mod}\, (\mcal{A} \boxtimes \mcal{B}^{\mrm{op}})
|_{U_{i} \times V_{i}} \big) \]
is a dualizing complex over 
$(\mcal{A}|_{U_{i}}, \mcal{B}|_{V_{i}})$. 
\end{enumerate}
\end{prop}

\begin{proof}
Assume $\mcal{R}$ has support in the graph of an isomorphism
$f: X \to Y$. We shall write
$\opnt{p}_{(1, i)}: U_i \times V_i \to U_i$
and 
$\opnt{p}_{(2, i)}: U_i \times V_i \to V_i$
for the projections on the open sets. Also we denote by
$\mrm{D}_i$ and $\mrm{D}_i^{\mrm{op}}$ the duality functors 
determined by $\mcal{R}|_{U_i \times V_i}$.
For any 
$\mcal{M} \in \msf{D}^{\mrm{b}}(\cat{Mod} \mcal{A})$
we have
\[ \begin{aligned}
(\mrm{D} \mcal{M})|_{V_i} & = 
(\mrm{R} \opnt{p}_{2 *} 
\mrm{R} \mcal{H}om_{\opnt{p}_{1}^{-1} \mcal{A}}
(\opnt{p}_{1}^{-1} \mcal{M}, \mcal{R}))|_{V_i} \\
& \cong 
\mrm{R} \opnt{p}_{(2, i) *} 
(\mrm{R} \mcal{H}om_{\opnt{p}_{1}^{-1} \mcal{A}}
(\opnt{p}_{1}^{-1} \mcal{M}, \mcal{R})|_{U_i \times V_i}) \\
& \cong 
\mrm{R} \opnt{p}_{(2, i) *} 
(\mrm{R} \mcal{H}om_{\opnt{p}_{(1, i)}^{-1} \mcal{A}}
(\opnt{p}_{(1, i)}^{-1} \mcal{M}|_{U_i}, \mcal{R})
|_{U_i \times V_i}) = \mrm{D}_i(\mcal{M}|_{U_i}) . 
\end{aligned} \]
Likewise 
$(\mrm{D}^{\mrm{op}} \mcal{N})|_{U_i} \cong 
\mrm{D}_i^{\mrm{op}}(\mcal{N}|_{V_i})$
for 
$\mcal{N} \in \msf{D}^{\mrm{b}}(\cat{Mod} \mcal{B}^{\mrm{op}})$.

Any coherent $\mcal{A}$-module restricts to a coherent 
$\mcal{A}|_{U_i}$-module, and the same for coherent 
$\mcal{B}^{\mrm{op}}$-modules. On the other hand by Lemma 
\ref{lem7.1} any coherent $\mcal{A}|_{U_i}$-module extends to a 
coherent $\mcal{A}$-module; and of course the same for 
coherent $\mcal{B}|_{V_i}$-modules. The upshot is that the three 
conditions of Definition \ref{dfn5.8} are satisfied for 
$\mcal{R}$ if and only if they are satisfied for all the complexes
$\mcal{R}|_{U_i \times V_i}$.
\end{proof}

Suppose we are given $\k$-algebras $A$ and $B$, 
and we view them as quasi-coherent 
ringed schemes $(\opn{Spec} \k, A)$ and
$(\opn{Spec} \k, B)$. 
As observed earlier in Example \ref{exa4.2}, for a complex
$R \in \msf{D}^{\mrm{b}}(\cat{Mod}\, A \otimes B^{\mrm{op}})$
Definition \ref{dfn5.8} becomes equivalent 
to the ring-theoretic Definition \ref{dfn2.1}.
As the next theorem shows the ``same'' is true for all affine 
quasi-coherent ringed schemes. 

\begin{thm} \label{thm7.2}
Let $(X, \mcal{A})$ and $(Y, \mcal{B})$ be affine
quasi-coherent ringed schemes over $\k$. 
Assume $\mcal{A}$ and $\mcal{B}^{\mrm{op}}$ 
are left noetherian and the product 
$(X \times Y, \mcal{A} \boxtimes \mcal{B}^{\mrm{op}})$
exists. Let 
$\mcal{R} \in \msf{D}^{\mrm{b}}_{\mrm{qc}}
(\cat{Mod}\, \mcal{A} \boxtimes \mcal{B}^{\mrm{op}})$
be some complex. Define $A := \Gamma(X, \mcal{A})$ and 
$B := \Gamma(Y, \mcal{B})$. 
Then the following two conditions are equivalent. 
\begin{enumerate}
\rmitem{i} $\mcal{R}$ is a dualizing complex over 
$(\mcal{A}, \mcal{B})$ in the sense of Definition 
\tup{\ref{dfn5.8}}.
\rmitem{ii} The complex
\[ R := \mrm{R} \Gamma(X \times Y, \mcal{R}) \in 
\msf{D}(\cat{Mod}\, A \otimes B^{\mrm{op}}) \]
is a dualizing complex over $(A, B)$
in the sense of Definition \tup{\ref{dfn2.1}}.
\end{enumerate}
\end{thm}

\begin{proof}
Let us write
$\mrm{D}_{\mrm{glob}} := \opn{RHom}_{A}(-, R)$ and 
$\mrm{D}_{\mrm{glob}}^{\mrm{op}} := 
\opn{RHom}_{B^{\mrm{op}}}(-, R)$.
According to Corollary \ref{cor5.5} the sheafification functors
\[ \begin{aligned}
\mcal{A} \otimes_{A} - & : 
\msf{D}^{\mrm{b}}_{\mrm{f}}(\cat{Mod} A) \to
\msf{D}^{\mrm{b}}_{\mrm{c}}(\cat{Mod} \mcal{A}) , \\
- \otimes_{B} \mcal{B} & : 
\msf{D}^{\mrm{b}}_{\mrm{f}}(\cat{Mod} B^{\mrm{op}}) \to
\msf{D}^{\mrm{b}}_{\mrm{c}}(\cat{Mod} \mcal{B}^{\mrm{op}})
\end{aligned} \]
are equivalences, with inverses 
$\mrm{R} \Gamma(X, -)$ and $\mrm{R} \Gamma(Y, -)$ respectively. 
Thus it suffices to show that the diagram
\[ \begin{CD}
\msf{D}^{\mrm{b}}_{\mrm{f}}(\cat{Mod} A)^{\mrm{op}}
@>{\mrm{D}_{\mrm{glob}}}>> 
\msf{D}^{\mrm{b}}_{\mrm{f}}(\cat{Mod} B^{\mrm{op}}) \\
@V{\mcal{A} \otimes_A -}VV @A{\mrm{R} \Gamma(Y, -)}AA \\
\msf{D}^{\mrm{b}}_{\mrm{c}}(\cat{Mod} \mcal{A})^{\mrm{op}}
@>{\mrm{D}}>> 
\msf{D}^{\mrm{b}}_{\mrm{c}}(\cat{Mod} \mcal{B}^{\mrm{op}})
\end{CD} \]
and the ``opposite'' diagram (the one involving 
$\mrm{D}^{\mrm{op}}$ and $\mrm{D}_{\mrm{glob}}^{\mrm{op}}$)
are commutative. By symmetry it suffices to check only one of 
them, say the one displayed.

We can assume that 
$\mcal{R} = (\mcal{A} \boxtimes \mcal{B}^{\mrm{op}})
\otimes_{A \otimes B^{\mrm{op}}} R$. 
Choose an injective resolution 
$R \to J$ in $\msf{C}^{+}(\cat{Mod}\, A \otimes B^{\mrm{op}})$,
and let
$\mcal{J} := (\mcal{A} \boxtimes \mcal{B}^{\mrm{op}})
\otimes_{A \otimes B^{\mrm{op}}} J$. 
Then $\mcal{R} \to \mcal{J}$ is a quasi-isomorphism. 
Now $\mcal{J}$ is a complex of injectives in 
$\cat{QCoh}\, \mcal{A} \boxtimes \mcal{B}^{\mrm{op}}$, 
but it might not be a complex of injectives in 
$\cat{Mod}\, \mcal{A} \boxtimes \mcal{B}^{\mrm{op}}$;  
cf.\ Remark \ref{rem5.2}. So we choose an injective resolution 
$\mcal{J} \to \mcal{K}$ in 
$\msf{C}^{+}(\cat{Mod}\, \mcal{A} \boxtimes \mcal{B}^{\mrm{op}})$. 
Let 
$\mcal{M} := \mcal{A} \otimes_A M 
\in \msf{D}^{\mrm{b}}_{\mrm{c}}(\cat{Mod} \mcal{A})$. 
Then
\[ \mrm{D} \mcal{M}  = 
\opnt{p}_{2 *} 
\mcal{H}om_{\opnt{p}_{1}^{-1} \mcal{A}}
(\opnt{p}_{1}^{-1} \mcal{M}, \mcal{K}) . \]
Since the latter is a complex of flasque 
$\mcal{B}^{\mrm{op}}$-modules on $Y$, we get
\[ \begin{aligned}
\mrm{R} \Gamma(Y, \mrm{D} \mcal{M}) & = 
\Gamma \bigl( Y, \opnt{p}_{2 *} 
\mcal{H}om_{\opnt{p}_{1}^{-1} \mcal{A}}
(\opnt{p}_{1}^{-1} \mcal{M}, \mcal{K}) \bigr) \\
& \cong \Gamma \bigl( X \times Y, 
\mcal{H}om_{\opnt{p}_{1}^{-1} \mcal{A}}
(\opnt{p}_{1}^{-1} \mcal{M}, \mcal{K}) \bigr) \\
& \cong \opn{Hom}_{\opnt{p}_{1}^{-1} \mcal{A}}
(\opnt{p}_{1}^{-1} \mcal{M}, \mcal{K}) .
\end{aligned} \]

Now choose a bounded above resolution $P \to M$ by finitely
generated free $A$-modules, and let 
$\mcal{P} := \mcal{A} \otimes_A P$. Then
$\opnt{p}_{1}^{-1} \mcal{P} \to \opnt{p}_{1}^{-1} \mcal{M}$
is a quasi-isomorphism of $\opnt{p}_{1}^{-1} \mcal{A}$ -modules, 
and so
\[ \opn{Hom}_{\opnt{p}_{1}^{-1} \mcal{A}}
(\opnt{p}_{1}^{-1} \mcal{M}, \mcal{K}) \to
\opn{Hom}_{\opnt{p}_{1}^{-1} \mcal{A}}
(\opnt{p}_{1}^{-1} \mcal{P}, \mcal{K}) \]
is a quasi-isomorphism of $B^{\mrm{op}}$-modules. 
Each $\mcal{J}^p$ and $\mcal{K}^p$ is acyclic for the functor
$\Gamma(X \times Y, -)$, since $\mcal{J}^p$ is quasi-coherent and 
$\mcal{K}^p$ is injective. Therefore 
$\Gamma(X \times Y, \mcal{J}) \to \Gamma(X \times Y, \mcal{K})$
is a quasi-isomorphism. Now 
$\opnt{p}_{1}^{-1} \mcal{P}$
is a bounded above complex of finitely generated free 
$\opnt{p}_{1}^{-1} \mcal{A}$ -modules, and thus  
\[ \opn{Hom}_{\opnt{p}_{1}^{-1} \mcal{A}}
(\opnt{p}_{1}^{-1} \mcal{P}, \mcal{J}) \to
\opn{Hom}_{\opnt{p}_{1}^{-1} \mcal{A}}
(\opnt{p}_{1}^{-1} \mcal{P}, \mcal{K}) \]
is a quasi-isomorphism of $B^{\mrm{op}}$-modules. But
\[ \opn{Hom}_{\opnt{p}_{1}^{-1} \mcal{A}}
(\opnt{p}_{1}^{-1} \mcal{P}, \mcal{J}) \cong
\opn{Hom}_{A}(P, J) = \opn{RHom}_{A}(M, R) = 
\mrm{D}_{\mrm{glob}} M. \]
\end{proof}

Next we wish to recall the definition of rigid dualizing complex 
over a $\k$-algebra, which is due to Van den Bergh \cite{VdB1}.

\begin{dfn} \label{dfn2.3}
Let $R$ be a dualizing complex over $A$. If there is an isomorphism
\[ \rho: R \to \opn{RHom}_{A^{\mrm{e}}}(A, R \otimes R) \]
in $\msf{D}(\cat{Mod} A^{\mrm{e}})$ then we call $(R, \rho)$, 
or just $R$, a {\em rigid} dualizing complex. 
The isomorphism $\rho$ is 
called a {\em rigidifying isomorphism}.
\end{dfn}

For a detailed explanation of this definition see 
\cite[Section 3]{YZ2}.
According to \cite[Proposition 8.2]{VdB1} and 
\cite[Corollary 3.4]{YZ2} a rigid dualizing complex $(R, \rho)$ is 
unique up to a unique isomorphism in 
$\msf{D}(\cat{Mod} A^{\mrm{e}})$. It is important 
to note that rigidity is a relative notion (relative to the base 
field $\k$). 

A ring homomorphism $A \to A'$ is called a {\em localization}
(in the sense of Silver \cite{Si}) if $A'$ is a flat
$A$-module on both sides and $A' \otimes_A A' \cong A'$.

\begin{dfn} \label{dfn6.1}
Let $A \to A'$ be a localization homomorphism between two noetherian
$\k$-algebras. Suppose the rigid dualizing complexes $(R, \rho)$ 
and $(R', \rho')$ of $A$ and $A'$ respectively exist. A {\em rigid 
localization morphism} is a morphism
\[ \mrm{q}_{A' / A} : R \to R' \]
in $\msf{D}(\cat{Mod} A^{\mrm{e}})$ satisfying the conditions
below. 
\begin{enumerate}
\rmitem{i} The morphisms $A' \otimes_A R \to R'$ 
and $R \otimes_A A' \to R'$ in $\msf{D}(\cat{Mod} A^{\mrm{e}})$
induced by $\mrm{q}_{A' / A}$ are isomorphisms.
\rmitem{ii} The diagram
\[ \begin{CD}
R @>{\rho}>> \opn{RHom}_{A^{\mrm{e}}}(A, R \otimes R) \\
@V \opn{q} VV @VV{\opn{q} \otimes \opn{q}}V \\
R' @>{\rho'}>> \opn{RHom}_{(A')^{\mrm{e}}}(A', R' \otimes R')
\end{CD} \]
in $\msf{D}(\cat{Mod} A^{\mrm{e}})$, where
$\mrm{q} := \mrm{q}_{A' / A}$, is commutative.
\end{enumerate}
We shall sometimes express this by saying that
$\mrm{q}_{A' / A} : (R, \rho) \to (R', \rho')$
is a rigid localization morphism.
\end{dfn}

According to \cite[Theorem 6.2]{YZ4} a rigid localization morphism 
$\mrm{q}_{A' / A} : (R, \rho) \to (R', \rho')$
is unique (if it exists). 

Let us denote by $\cat{Aff} X$ the set of affine open subsets of 
$X$. Suppose $\mcal{R}$ is a local dualizing complex over 
$\mcal{A}$ supported on the diagonal
$\Delta(X) \subset X^2$. By Proposition \ref{prop7.3} and
Theorem \ref{thm7.2}, 
for every $U \in \cat{Aff} X$ the complex
$\mrm{R} \Gamma(U^2, \mcal{R})$ is a dualizing complex over 
$A := \Gamma(U, \mcal{A})$. If $V \subset U$ is another affine 
open set and $A' := \Gamma(V, \mcal{A})$ then both 
$A \to A'$ and $A^{\mrm{e}} \to {A'}^{\mrm{e}}$
are localizations. Moreover the restriction 
$\msf{D}^{\mrm{b}}_{\mrm{c}}(\cat{Mod} \mcal{A}^{\mrm{e}}|_{U^2})
\to 
\msf{D}^{\mrm{b}}_{\mrm{c}}(\cat{Mod} \mcal{A}^{\mrm{e}}|_{V^2})$
induces a morphism $\mrm{q}_{A' / A} : R \to R'$
in $\msf{D}(\cat{Mod} A^{\mrm{e}})$; cf.\ Lemma \ref{lem5.4}.

\begin{dfn} \label{dfn7.5}
Let $(X, \mcal{A})$ be a separated noetherian quasi-coherent 
ringed scheme over $\k$. Assume the product 
$(X^2, \mcal{A}^{\mrm{e}})$ exists, and is also noetherian. A 
{\em rigid dualizing complex} over $\mcal{A}$ is a pair
$(\mcal{R}, \bsym{\rho})$, where:
\begin{enumerate}
\item 
$\mcal{R} \in \msf{D}^{\mrm{b}}_{\mrm{c}}
(\cat{Mod} \mcal{A}^{\mrm{e}})$ 
is a local dualizing complex over $\mcal{A}$ supported
on the diagonal $\Delta(X) \subset X^2$. 
\item $\bsym{\rho} = \{ \rho_U \}_{U \in \cat{Aff} X}$ 
is a collection of rigidifying isomorphisms, namely for each 
$U \in \cat{Aff} X$, letting $A := \Gamma(U, \mcal{A})$ and
$R:= \mrm{R} \Gamma(U^2, \mcal{R})$, the pair
$(R, \rho_U)$ is a rigid dualizing complex over $A$.
\end{enumerate}
The following compatibility condition is required of the data
$(\mcal{R}, \bsym{\rho})$:
\begin{enumerate}
\item[($*$)] Given $V \subset U$ in $\cat{Aff} X$, write
$A' := \Gamma(V, \mcal{A})$ and 
$R' := \mrm{R} \Gamma(V^2, \mcal{R})$.
Let $\mrm{q}_{A' / A} : R \to R'$
be the morphism in $\msf{D}(\cat{Mod} A^{\mrm{e}})$ coming from  
restriction. Then 
\[ \mrm{q}_{A' / A} : (R, \rho_U) \to (R', \rho_V) \]
is a rigid localization morphism. 
\end{enumerate}
\end{dfn}

\begin{exa} \label{exa7.4} 
Let $X$ be a separated finite type $\k$-scheme, smooth of 
dimension $n$. There is a canonical isomorphism 
\[ \rho : \Delta_* \Omega^n_{X / \k} \iso
\mcal{E}xt^n_{\mcal{O}_{X^2}}
(\Delta_* \mcal{O}_{X}, \Omega^{2n}_{X^2 / \k}) , \]
and 
$\mcal{E}xt^i_{\mcal{O}_{X^2}}
(\Delta_* \mcal{O}_{X}, \Omega^{2n}_{X^2 / \k}) = 0$
for $i \neq n$ (see \cite[Proposition III.7.2]{RD}). 
Therefore on any affine open set $U = \opn{Spec} A \subset X$
we get an isomorphism 
\[ \rho_U : \Omega^n_{A / \k}[n] \iso
\opn{RHom}_{A^{\mrm{e}}} \bigl( A, \Omega^n_{A / \k}[n] \otimes 
\Omega^n_{A / \k}[n] \bigr) \]
in $\msf{D}(\cat{Mod} A^{\mrm{e}})$, and the collection
$\bsym{\rho} := \{ \rho_U \}$
is compatible with localization. We see that 
$(\Delta_* \Omega^n_{A / \k}[n], \bsym{\rho})$
is a rigid dualizing complex over $X$ in the sense of 
Definition \ref{dfn7.5}. 
\end{exa}

\begin{dfn}
\label{dfn7.6}
A morphism $f: (Y, \mcal{B}) \to (X, \mcal{A})$ between noetherian 
quasi-coherent ringed schemes is called {\em finite} if 
$f: Y \to X$ is finite and 
$f_{*} \mcal{B}$ is a coherent $\mcal{A}$-module on both sides.
\end{dfn}

Given a morphism $f:(Y, \mcal{B}) \to (X, \mcal{A})$ of ringed 
schemes then 
$f^{\mrm{op}} : (Y, \mcal{B}^{\mrm{op}}) \to 
(X, \mcal{A}^{\mrm{op}})$ 
is also a morphism.
According to Proposition \ref{prop4.1}, if the products exist then
there is a morphism
\[ f^{\mrm{e}} := f \times f^{\mrm{op}} : 
(Y^{2}, \mcal{B}^{\mrm{e}}) \to 
(X^{2}, \mcal{A}^{\mrm{e}}) . \]

\begin{dfn} \label{dfn2.4}
Let $A \to B$ be a finite homomorphism of $\k$-algebras. Assume the
rigid dualizing complexes $(R_A, \rho_A)$ and $(R_B, \rho_B)$ exist. 
Let $\opn{Tr}_{B / A}: R_B \to R_A$ be a morphism in 
$\msf{D}(\cat{Mod} A^{\mrm{e}})$.
We say $\opn{Tr}_{B / A}$ is a {\em rigid trace} if it satisfies  
the following two conditions:
\begin{enumerate}
\rmitem{i} $\opn{Tr}_{B / A}$ induces isomorphisms
\[ R_B \cong \opn{RHom}_A(B, R_A) \cong 
\opn{RHom}_{A^{\mrm{op}}}(B, R_A) \]
in $\msf{D}(\cat{Mod} A^{\mrm{e}})$.
\rmitem{ii} The diagram
\[ \begin{CD}
R_B @>{\rho_B}>> \opn{RHom}_{B^{\mrm{e}}}(B, R_B \otimes R_B) \\
@V \opn{Tr} VV @VV{\opn{Tr} \otimes \opn{Tr}}V \\
R_A @>{\rho_A}>> \opn{RHom}_{A^{\mrm{e}}}(A, R_A \otimes R_A)
\end{CD} \]
in $\msf{D}(\cat{Mod} A^{\mrm{e}})$, where
$\opn{Tr} := \opn{Tr}_{B / A}$, is commutative.
\end{enumerate}
Often we shall say that
$\opn{Tr}_{B / A} : (R_B, \rho_B) \to (R_A, \rho_A)$
is a rigid trace morphism.
\end{dfn}

By \cite[Theorem 3.2]{YZ2}, a rigid trace $\opn{Tr}_{B / A}$ 
is unique (if it exists). 

\begin{dfn} \label{dfn7.7}
Let $f: (Y, \mcal{B}) \to (X, \mcal{A})$ be a finite morphism
between noetherian separated quasi-coherent ringed $\k$-schemes.
Assume both $(X, \mcal{A})$ and $(Y, \mcal{B})$ have rigid 
dualizing complexes $(\mcal{R}_{\mcal{A}}, \bsym{\rho}_{\mcal{A}})$ 
and 
$(\mcal{R}_{\mcal{B}}, \bsym{\rho}_{\mcal{B}})$ respectively.
A {\em rigid trace} is a morphism 
\[ \opn{Tr}_{f}: \mrm{R} f^{\mrm{e}}_{*} \mcal{R}_{\mcal{B}} \to
\mcal{R}_{\mcal{A}} \]
in $\msf{D}(\cat{Mod} \mcal{A}^{\mrm{e}})$ satisfying the 
following condition. 
\begin{enumerate}
\item[($\dag$)] Let $U \subset X$ be any affine open set,
$V := f^{-1}(U)$, $A:= \Gamma(U, \mcal{A})$ and 
$B:= \Gamma(V, \mcal{B})$. 
Let
$R_{A}:= \mrm{R} \Gamma(U^2, \mcal{R}_{\mcal{A}})$
and 
$R_{B} := \mrm{R} \Gamma(V^2, \mcal{R}_{\mcal{B}})$
be the rigid dualizing complexes, with their respective 
rigidifying isomorphisms $\rho_U$ and $\rho_V$. Then 
\[ \mrm{R} \Gamma(U^2, \opn{Tr}_{f}):
(R_{B}, \rho_V) \to (R_{A}, \rho_U) \]
is a rigid trace morphism.  
\end{enumerate}
\end{dfn}

\begin{exa}
Suppose $X$ and $Y$ are separated $\k$-schemes, smooth of 
dimensions $m$ and $n$ respectively. Let $f : X \to Y$ be a finite 
morphism. According to \cite[Theorem III.10.5]{RD} 
there is a trace morphism
$\opn{Tr}_f : f_* \Omega^m_{X / \k}[m] \to 
\Omega^n_{Y / \k}[n]$
in $\msf{D}(\cat{Mod} \mcal{O}_Y)$. We know that 
$\Delta_* \Omega^m_{X / \k}[m]$ and 
$\Delta_* \Omega^n_{Y / \k}[n]$ are the rigid 
dualizing complexes of these schemes (see Example \ref{exa7.4}). 
Condition TRA2 of \cite[Theorem III.10.2]{RD} 
implies that 
$\Delta_*(\opn{Tr}_f)$ it is a rigid trace 
in the sense of Definition \ref{dfn7.7}. 
\end{exa}

\begin{rem}
For $\k$-algebras one has the notion of {\em Auslander dualizing 
complex}, see \cite{YZ2}. In \cite{YZ4} it was shown that the 
Auslander condition is closely related to the rigid perverse 
t-structure. It should be interesting to extend the 
Auslander condition to  the geometric context, i.e.\ to dualizing 
complexes over quasi-coherent ringed schemes.
\end{rem}

\section{Perverse Sheaves on Ringed Spaces}
\label{sec4}

This section deals with gluing t-structures, in a rather general 
context. 
Let us begin by recalling the following basic definition due to 
Beilinson, Bernstein and Deligne \cite{BBD}. We shall 
follow the exposition in \cite[Chapter X]{KS}.  

\begin{dfn} \label{dfn6.2} 
Suppose $\cat{D}$ is a triangulated category and
$\cat{D}^{\leq 0}, \cat{D}^{\geq 0}$ are two full subcategories.
Let
$\cat{D}^{\leq n} := \cat{D}^{\leq 0}[-n]$
and
$\cat{D}^{\geq n} := \cat{D}^{\geq 0}[-n]$.
We say $(\cat{D}^{\leq 0}, \cat{D}^{\geq 0})$ is a
{\em t-structure on $\cat{D}$} if:
\begin{enumerate}
\rmitem{i} $\cat{D}^{\leq -1} \subset \cat{D}^{\leq 0}$
and
$\cat{D}^{\geq 1} \subset \cat{D}^{\geq 0}$.
\rmitem{ii} $\opn{Hom}_{\cat{D}}(M, N) = 0$ for
$M \in \cat{D}^{\leq 0}$ and $N \in \cat{D}^{\geq 1}$.
\rmitem{iii} For any $M \in \cat{D}$ there is a distinguished
triangle
\[ M' \to M \to M'' \to M'[1] \]
in $\cat{D}$ with $M'\in \cat{D}^{\leq 0}$ and
$M'' \in \cat{D}^{\geq 1}$.
\end{enumerate}
When these conditions are satisfied we define
the {\em heart of $\cat{D}$} to be the full subcategory
$\cat{D}^{0} := \cat{D}^{\leq 0} \cap \cat{D}^{\geq 0}$.
\end{dfn}

It is known that the heart $\cat{D}^{0}$ is an abelian category, 
in which short exact sequences are distinguished triangles in 
$\cat{D}$ with vertices in $\cat{D}^{0}$.

Let $(X, \mcal{A})$ be a ringed space, i.e.\ a topological space 
$X$ endowed with a sheaf of (not necessarily commutative) rings 
$\mcal{A}$. We denote by $\cat{Mod} \mcal{A}$ the category of 
sheaves of left $\mcal{A}$-modules, and by
$\msf{D}(\cat{Mod} \mcal{A})$
the derived category. 

The triangulated category
$\msf{D}(\cat{Mod} \mcal{A})$
has the {\em standard t-structure}, in which 
\[ \begin{aligned}
\msf{D}(\cat{Mod} \mcal{A})^{\leq 0}
& := \{ \mcal{M} \in 
\msf{D}(\cat{Mod} \mcal{A}) \mid
\mrm{H}^i \mcal{M} = 0 \text{ for all } i > 0 \} , \\
\msf{D}(\cat{Mod} \mcal{A})^{\geq 0}
& := \{ \mcal{M} \in 
\msf{D}(\cat{Mod} \mcal{A}) \mid
\mrm{H}^i \mcal{M} = 0 \text{ for all } i < 0 \} .
\end{aligned} \]
The heart 
$\msf{D}(\cat{Mod} \mcal{A})^{0}$
is equivalent to $\cat{Mod} \mcal{A}$. Other t-structures on 
$\msf{D}(\cat{Mod} \mcal{A})$, or on some triangulated full 
subcategory 
$\msf{D} \subset \msf{D}(\cat{Mod} \mcal{A})$,
will be called {\em perverse t-structures}, and the notation 
${}^{p}\msf{D}^{\star}$ shall be used. 

A stack on $X$ is a ``sheaf of categories.'' The general definition 
(cf.\ \cite{LMB}) is quite forbidding; 
but we shall only need the following special instance
(cf.\ \cite[Section X.10]{KS}). Given two open sets 
$V \subset U$ in $X$, the restriction functor 
$\msf{D}(\cat{Mod} \mcal{A}|_U) \to \msf{D}(\cat{Mod} \mcal{A}|_V)$ 
is denoted by
$\mcal{M} \mapsto \mcal{M}|_V$.

\begin{dfn} \label{dfn4.4}
Let $(X, \mcal{A})$ be a ringed space.
Suppose that for every open set $U \subset X$ 
we are given a full subcategory 
$\cat{C}(U) \subset \msf{D}(\cat{Mod} \mcal{A}|_U)$.
The collection of categories 
$\cat{C} = \{ \cat{C}(U) \}$ is called a {\em stack of 
subcategories of $\msf{D}(\cat{Mod} \mcal{A})$} 
if the following axioms hold.
\begin{enumerate}
\rmitem{a} Let $V \subset U$ be open sets in $X$ and 
$\mcal{M} \in \cat{C}(U)$. Then 
$\mcal{M}|_V \in \cat{C}(V)$. 
\rmitem{b} Descent for objects: 
given an open covering $U = \bigcup V_{i}$, objects
$\mcal{M}_{i} \in \cat{C}(V_{i})$ and isomorphisms
$\phi_{i, j} : \mcal{M}_{i}|_{V_{i} \cap V_{j}} \iso
\mcal{M}_{j}|_{V_{i} \cap V_{j}}$
satisfying the cocycle condition
$\phi_{i, k} = \phi_{j, k} \circ \phi_{i, j}$
on triple intersections, there exists an object 
$\mcal{M} \in \cat{C}(U)$ and isomorphisms
$\phi_{i}: \mcal{M}|_{V_{i}} \iso \mcal{M}_{i}$
such that 
$\phi_{i, j}\circ \phi_{i} = \phi_{j}$.
\rmitem{c} Descent for morphisms: 
given two objects $\mcal{M}, \mcal{N} \in \cat{C}(U)$,
an open covering $U = \bigcup V_{i}$ and morphisms 
$\psi_i : \mcal{M}|_{V_i} \to \mcal{N}|_{V_i}$
such that 
$\psi_i|_{V_{i} \cap V_{j}} = \psi_j|_{V_{i} \cap V_{j}}$,
there is a unique morphism $\psi : \mcal{M} \to \mcal{N}$
such that $\psi|_{V_i} = \psi_i$.
\end{enumerate}
\end{dfn}

Unlike \cite{LMB}, our stacks do not consist of groupoids. On the 
contrary, we will work with stacks of abelian categories. Here is 
an example.

\begin{exa}
Take 
$\cat{C}(U) := 
\msf{D}(\cat{Mod} \mcal{A}|_U)^{0}$,
the heart for the standard t-structure. Since
$\cat{C}(U)$ is canonically equivalent to $\cat{Mod} \mcal{A}|_U$
it follows that $\cat{C} = \{ \cat{C}(U) \}$
is a stack of subcategories of $\msf{D}(\cat{Mod} \mcal{A})$.
\end{exa}

By abuse of notation we shall denote by 
$\cat{Mod} \mcal{A}$ the stack 
$U \mapsto \msf{D}(\cat{Mod} \mcal{A}|_U)^{0}$.

\begin{dfn}
Suppose $\msf{C}$ is a stack of subcategories 
of $\msf{D}(\cat{Mod} \mcal{A})$ such that for any open set $U$ the 
subcategory $\msf{C}(U)$ is a thick abelian subcategory
of $\msf{D}(\cat{Mod} \mcal{A}|_U)^{0}$. Then we call $\msf{C}$ a 
{\em thick abelian substack of $\cat{Mod} \mcal{A}$}.
\end{dfn}

\begin{exa}
Consider a noetherian scheme $(X, \mcal{O}_X)$. 
Define $\msf{C}(U) := \cat{Coh} \mcal{O}_U$, 
the category of coherent sheaves on $U$. 
This is a thick abelian substack of 
$\cat{Mod} \mcal{O}_X$, which we denote (by abuse of notation) 
$\cat{Coh} \mcal{O}_X$.
\end{exa}

 From now until the end of this section we fix a ringed space
$(X, \mcal{A})$, a thick abelian substack
$\msf{C} \subset \cat{Mod} \mcal{A}$ and a basis
$\mfrak{B}$ of the topology of $X$.
For any open set $U \subset X$ the full subcategory
$\msf{D}^{\mrm{b}}_{\mrm{c}}(\cat{Mod} \mcal{A}|_U)$,
whose objects are the bounded complexes $\mcal{M}$ 
such that $\mrm{H}^i \mcal{M} \in \msf{C}(U)$
for all $i$, is a triangulated subcategory of
$\msf{D}^{\mrm{b}}(\cat{Mod} \mcal{A}|_U)$. 

\begin{dfn} \label{dfn6.3} 
Suppose that for every open set $U \in \mfrak{B}$ 
we are given a t-structure
\[ \mfrak{T}_U =
\bigl( 
{}^{p}\msf{D}^{\mrm{b}}_{\mrm{c}}(\cat{Mod} \mcal{A}|_U)^{\leq 0},
{}^{p}\msf{D}^{\mrm{b}}_{\mrm{c}}(\cat{Mod} \mcal{A}|_U)^{\geq 0}
\bigr) \]
on 
$\msf{D}^{\mrm{b}}_{\mrm{c}}(\cat{Mod} \mcal{A}|_U)$. 
Furthermore suppose this collection of t-structures 
$\{ \mfrak{T}_U \}_{U \in \mfrak{B}}$
satisfies the following condition. 
\begin{enumerate}
\item[($\diamondsuit$)] 
Let $U \in \mfrak{B}$, let $U = \bigcup U_i$ be any covering 
with $U_i \in \mfrak{B}$, and let ``$\star$'' denote either 
``$\leq 0$'' or ``$\geq 0$''. Then for any 
$\mcal{M} \in \msf{D}^{\mrm{b}}_{\mrm{c}}(\cat{Mod} \mcal{A}|_U)$
the following are equivalent:
\begin{enumerate}
\rmitem{i} 
$\mcal{M} \in 
{}^{p}\msf{D}^{\mrm{b}}_{\mrm{c}}(\cat{Mod} \mcal{A}|_U)^{\star}$.
\rmitem{ii} 
$\mcal{M}|_{U_i} \in 
{}^{p}\msf{D}^{\mrm{b}}_{\mrm{c}}
(\cat{Mod} \mcal{A}|_{U_i})^{\star}$
for all $i$.
\end{enumerate}
\end{enumerate}
Then we call $\{ \mfrak{T}_U \}_{U \in \mfrak{B}}$
a {\em local collection of t-structures} on $(X, \mcal{A})$.
\end{dfn}

\begin{dfn} \label{dfn4.3}
Let $\{ \mfrak{T}_U \}_{U \in \mfrak{B}}$ be a local 
collection of t-structures on $(X, \mcal{A})$. 
Define full subcategories
\[ \begin{aligned}
& {}^{p}\msf{D}^{\mrm{b}}_{\mrm{c}}(\cat{Mod} \mcal{A})^{\star} 
:= \\
& \qquad \{ \mcal{M} \in 
\msf{D}^{\mrm{b}}_{\mrm{c}}(\cat{Mod} \mcal{A}) & \mid \
\mcal{M}|_U \in
{}^{p}\msf{D}^{\mrm{b}}_{\mrm{c}}
(\cat{Mod} \mcal{A}|_U)^{\star} 
\text{ for every } U \in \mfrak{B} \} ,
\end{aligned} \]
where ``$\star$'' is either ``$\leq$'' or ``$\geq$''.
\end{dfn}

The next lemmas are modifications of material in 
\cite[Section 10.2]{KS}. We assume a local collection of 
t-structures 
$\{ \mfrak{T}_U \}_{U \in \mfrak{B}}$
is given.

\begin{lem} \label{lem4.4}
Let
$\mcal{M} \in 
{}^{p}\msf{D}^{\mrm{b}}_{\mrm{c}}(\cat{Mod} \mcal{A})^{\leq 0}$
and
$\mcal{N} \in 
{}^{p}\msf{D}^{\mrm{b}}_{\mrm{c}}(\cat{Mod} \mcal{A})^{\geq 0}$.
\begin{enumerate}
\item The sheaf 
$\mrm{H}^{i} \mrm{R} \mcal{H}om_{\mcal{A}}(\mcal{M}, \mcal{N})$
vanishes for all $i \leq -1$.
\item The assignment 
$U \mapsto \opn{Hom}_{\msf{D}(\cat{Mod} \mcal{A}|_{U})}
(\mcal{M}|_{U}, \mcal{N}|_{U})$,
for open sets $U \subset X$, is a sheaf on $X$.
\item 
$\opn{Hom}_{\msf{D}(\cat{Mod} \mcal{A})}
(\mcal{M}, \mcal{N}[i]) = 0$
for all $i \leq -1$.
\end{enumerate}
\end{lem}

\begin{proof}
(1) We note that 
\[ \mrm{H}^{i} \mrm{R} \mcal{H}om_{\mcal{A}}(\mcal{M}, \mcal{N}) 
\cong 
\mrm{H}^{0} \mrm{R} \mcal{H}om_{\mcal{A}}(\mcal{M}, \mcal{N}[i]) 
, \]
and 
$\mcal{N}[i] \in 
{}^{p}\msf{D}^{\mrm{b}}_{\mrm{c}}(\cat{Mod} \mcal{A})^{\geq 1}$,
because $i \leq -1$. 

Now
$\mrm{H}^{0} \mrm{R} \mcal{H}om_{\mcal{A}}(\mcal{M}, \mcal{N}[i])$
is isomorphic to the sheaf associated to the presheaf 
\[ U \mapsto \mrm{H}^{0} \opn{RHom}_{\mcal{A}|_{U}}
(\mcal{M}|_{U}, \mcal{N}[i]|_{U}) \cong
\opn{Hom}_{\msf{D}(\cat{Mod} \mcal{A}|_{U})}
(\mcal{M}|_{U}, \mcal{N}[i]|_{U}) . \]
Thus it suffices to prove that
\begin{equation} \label{eqn4.1}
\opn{Hom}_{\msf{D}(\cat{Mod} \mcal{A}|_{U})}
(\mcal{M}|_{U}, \mcal{N}[i]|_{U}) = 0 
\end{equation} 
for all open sets $U \in \mfrak{B}$. 
Since by Definition \ref{dfn4.3}
we have $\mcal{M}|_{U} \in 
{}^{p}\msf{D}^{\mrm{b}}_{\mrm{c}}
(\cat{Mod} \mcal{A}|_{U})^{\leq 0}$
and
$\mcal{N}[i]|_{U} \in 
{}^{p}\msf{D}^{\mrm{b}}_{\mrm{c}}
(\cat{Mod} \mcal{A}|_{U})^{\geq 1}$,
the assertion follows from condition (ii) in Definition 
\ref{dfn6.2}.

\medskip \noindent
(2) Let us write
$\mcal{L} := \mrm{R} \mcal{H}om_{\mcal{A}}(\mcal{M}, \mcal{N})
\in \msf{D}(\cat{Mod} \mbb{Z}_X)$. 
By part (1) we know that
$\mrm{H}^i \mcal{L} = 0$ for all $i < 0$. So after truncation
we can assume $\mcal{L}^i = 0$ for all $i < 0$. Hence 
$\mrm{H}^{0} \mrm{R} \Gamma(U, \mcal{L})
\cong \Gamma(U, \mrm{H}^{0} \mcal{L}).$
But 
\[ \mrm{H}^{0} \mrm{R} \Gamma(U, \mcal{L})
\cong \opn{Hom}_{\msf{D}(\cat{Mod} \mcal{A}|_{U})}
(\mcal{M}|_{U}, \mcal{N}|_{U}) . \]
We see that the presheaf
$U \mapsto \opn{Hom}_{\msf{D}(\cat{Mod} \mcal{A}|_{U})}
(\mcal{M}|_{U}, \mcal{N}|_{U})$
is actually a sheaf, namely the sheaf 
$\mrm{H}^{0} \mcal{L}$.

\medskip \noindent
(3) Here 
$\mcal{N}[i] \in 
{}^{p}\msf{D}^{\mrm{b}}_{\mrm{c}}(\cat{Mod} \mcal{A})^{\geq 1}
\subset
{}^{p}\msf{D}^{\mrm{b}}_{\mrm{c}}(\cat{Mod} \mcal{A})^{\geq 0}$,
so equation (\ref{eqn4.1}) holds for all open sets in 
$\mfrak{B}$. Applying part (2) of the lemma 
to $\mcal{M}$ and $\mcal{N}[i]$ we get 
$\opn{Hom}_{\msf{D}(\cat{Mod} \mcal{A})}
(\mcal{M}, \mcal{N}[i]) = 0$. 
\end{proof}

\begin{lem} \label{lem4.5}
Suppose the two rows in the diagram below are distinguished
triangles with
$\mcal{M}'_{1} \in {}^{p}\msf{D}^{\mrm{b}}_{\mrm{c}}
(\cat{Mod} \mcal{A})^{\leq 0}$
and
$\mcal{M}''_{2} \in 
{}^{p}\msf{D}^{\mrm{b}}_{\mrm{c}} 
(\cat{Mod} \mcal{A})^{\geq 1}$. 
Given any morphism 
$g: \mcal{M}_{1} \to \mcal{M}_{2}$
there exist unique morphisms
$f: \mcal{M}'_{1} \to \mcal{M}'_{2}$
and
$h: \mcal{M}''_{1} \to \mcal{M}''_{2}$
making the diagram commutative. If $g$ is an isomorphism then so 
are $f$ and $h$.
\[ \begin{CD}
\mcal{M}'_{1} @>{\alpha_{1}}>> \mcal{M}_{1} @>{\beta_{1}}>> 
\mcal{M}''_{1} @>{\gamma_{1}}>> \mcal{M}'_{1}[1] \\
@V{f}VV @V{g}VV @V{h}VV @V{f[1]}VV  \\
\mcal{M}'_{2} @>{\alpha_{2}}>> \mcal{M}_{2} @>{\beta_{2}}>> 
\mcal{M}''_{2} @>{\gamma_{2}}>> \mcal{M}'_{2}[1]
\end{CD} \]
\end{lem}

\begin{proof}
By assumption we have
$\mcal{M}'_{1} \in {}^{p}\msf{D}^{\mrm{b}}_{\mrm{c}}
(\cat{Mod} \mcal{A})^{\leq 0}$
and
$\mcal{M}''_{2}[i] \in 
{}^{p}\msf{D}^{\mrm{b}}_{\mrm{c}} 
(\cat{Mod} \mcal{A})^{\geq 1}$
for all $i \leq 0$. So by to Lemma \ref{lem4.4}(3) we get
$\opn{Hom}_{\msf{D}(\cat{Mod} \mcal{A})}
(\mcal{M}'_1, \mcal{M}''_2[i]) = 0$
for $i \leq 0$. 
According to \cite[Proposition 1.1.9]{BBD} there exist unique
morphisms $f$ and $h$ making the diagram commutative. 
The assertion about isomorphisms is also in 
\cite[Proposition 1.1.9]{BBD}.
\end{proof}

If $V \subset X$ is any open subset define
$\mfrak{B}|_V := \{ U \in \mfrak{B} \mid U \subset V \}$.
Then 
$\{ \mfrak{T}_U \}_{U \in \mfrak{B}|_V}$ is a local 
collection of t-structures on the ringed space
$(V, \mcal{A}|_V)$; 
hence using Definition \ref{dfn4.3} we get subcategories
${}^{p}\msf{D}^{\mrm{b}}_{\mrm{c}}
(\cat{Mod} \mcal{A}|_V)^{\star}$
of 
$\msf{D}^{\mrm{b}}_{\mrm{c}}(\cat{Mod} \mcal{A}|_V)$.

\begin{lem} \label{lem6.1} 
Let $V \subset X$ be an open set,
and let $V = \bigcup_{j} V_{j}$ be some open covering 
\tup{(}where $V$ and $V_j$ are not necessarily in
$\mfrak{B}$\tup{)}. 
Write ``$\star$'' for either ``$\leq 0$'' 
or ``$\geq 0$''. Then the following are equivalent for
$\mcal{M} \in 
\msf{D}^{\mrm{b}}_{\mrm{c}}(\cat{Mod} \mcal{A}|_V)$.
\begin{enumerate}
\rmitem{i} 
$\mcal{M} \in 
{}^{p}\msf{D}^{\mrm{b}}_{\mrm{c}}(\cat{Mod} \mcal{A}|_V)^{\star}$.
\rmitem{ii} 
$\mcal{M}|_{V_j} \in
{}^{p}\msf{D}^{\mrm{b}}_{\mrm{f}}
(\cat{Mod} \mcal{A}|_{V_j})^{\star}$
for all $j$.
\end{enumerate}
\end{lem}

\begin{proof}
The implication (i) $\Rightarrow$ (ii) is immediate from 
Definition \ref{dfn4.3}. For the reverse implication 
use condition ($\diamondsuit$) of Definition 
\ref{dfn6.3} twice with suitable $\mfrak{B}$-coverings of $V$ 
and of $V_j$.
\end{proof}

\begin{lem} \label{lem4.8}
Let $U \subset X$ be an open set, and let
$U = \bigcup_{i=1}^n U_i$ be some covering by open sets 
\tup{(}$U$ and $U_i$ not necessarily in $\mfrak{B}$\tup{)}. Let 
$\mcal{M} \in \msf{D}^{\mrm{b}}_{\mrm{c}}(\cat{Mod} \mcal{A}|_U)$,
and suppose there are distinguished triangles
\begin{equation}
\label{equ4.9}
T_{i}= \Bigl(
\mcal{M}_{i}' \xar{\alpha_{i}} \mcal{M}|_{U_{i}} \xar{\beta_{i}}
\mcal{M}_{i}'' \xar{\gamma_{i}} \mcal{M}_{i}'[1]
\Bigr)
\end{equation}
in $\msf{D}^{\mrm{b}}_{\mrm{c}}(\cat{Mod} \mcal{A}|_{U_{i}})$
with 
$\mcal{M}_{i}' \in 
{}^{p}\msf{D}^{\mrm{b}}_{\mrm{c}}
(\cat{Mod} \mcal{A}|_{U_{i}})^{\leq 0}$
and 
$\mcal{M}_{i}'' \in 
{}^{p}\msf{D}^{\mrm{b}}_{\mrm{c}}
(\cat{Mod} \mcal{A}|_{U_{i}})^{\geq 1}$.
Then there exists a distinguished triangle
\[ \mcal{M}' \xar{\alpha} \mcal{M} \xar{\beta}
\mcal{M}'' \xar{\gamma} \mcal{M}'[1] \]
in $\msf{D}^{\mrm{b}}_{\mrm{c}}(\cat{Mod} \mcal{A}|_U)$
whose restriction to $U_{i}$ is isomorphic to $T_{i}$.
\end{lem}

\begin{proof} First assume $n=2$. 
Denote by $g_{i} : U_{i} \to U$ the inclusions, for $i = 1, 2$.
Also write
$U_{(1, 2)} := U_{1} \cap U_{2}$ and
$g_{(1, 2)} : U_{(1, 2)} \to U$. For any open immersion $g$ let
$g_{!}$ be extension by zero, which is an exact functor.

The restriction of the triangles $T_{1}$ and $T_{2}$ to 
$U_{(1, 2)}$ and Lemma \ref{lem4.5} give rise to an isomorphism
$f: \mcal{M}'_{1}|_{U_{(1, 2)}} \iso \mcal{M}'_{2}|_{U_{(1, 2)}}$
in 
$\msf{D}^{\mrm{b}}_{\mrm{c}}(\cat{Mod} \mcal{A}|_{U_{(1, 2)}})$
satisfying 
$\alpha_{2} \circ f = \alpha_{1}$.
Therefore we get a morphism
\[ \delta: g_{(1, 2)!} (\mcal{M}'_{1}|_{U_{(1, 2)}}) \to
g_{1 !} \mcal{M}'_{1} \oplus g_{2 !} \mcal{M}'_{2} \]
in $\msf{D}^{\mrm{b}}(\cat{Mod} \mcal{A}|_{U})$, whose components 
are extensions by zero of the identity and $f$ respectively. 
(Note that the complexes $g_{1 !} \mcal{M}'_{1}$ etc.\ might 
not have cohomologies in $\msf{C}$.)

Define $\mcal{M}'$ to be the cone of $\delta$. So there is a  
distinguished triangle
\[ g_{(1, 2)!} (\mcal{M}'_{1}|_{U_{(1, 2)}}) \xar{\delta}
g_{1 !} \mcal{M}'_{1} \oplus g_{2 !} \mcal{M}'_{2}
\to \mcal{M}' \to g_{(1, 2)!} (\mcal{M}'_{1}|_{U_{(1, 2)}})[1] \]
in $\msf{D}(\cat{Mod} \mcal{A}|_{U})$. 
Applying the cohomological functor
$\opn{Hom}_{\msf{D}(\cat{Mod} \mcal{A}|_{U})}(-, \mcal{M}|_{U})$
to this triangle we get an exact sequence
\[ \begin{aligned}
\opn{Hom}_{\msf{D}(\cat{Mod} \mcal{A}|_{U})}
(\mcal{M}', \mcal{M}|_{U})
& \to
\opn{Hom}_{\msf{D}(\cat{Mod} \mcal{A}|_{U})}
(g_{1 !} \mcal{M}'_{1} \oplus g_{2 !} \mcal{M}'_{2}, 
\mcal{M}|_{U}) \\
& \to
\opn{Hom}_{\msf{D}(\cat{Mod} \mcal{A}|_{U})}
\bigl( g_{(1, 2)!} (\mcal{M}'_{1}|_{U_{(1, 2)}}), \mcal{M}|_{U} 
\bigr) .
\end{aligned} \]
The pair $(\alpha_{1}, -\alpha_{2})$ in the middle term goes to 
zero, and hence it comes from some morphism
$\alpha: \mcal{M}' \to \mcal{M}|_{U}$.
By construction the restriction 
$\mcal{M}'|_{U_{i}} \cong \mcal{M}'_{i}$,
so with the help of Lemma \ref{lem6.1} we deduce that 
$\mcal{M}' \in 
{}^{p}\msf{D}^{\mrm{b}}_{\mrm{c}}
(\cat{Mod} \mcal{A}|_{U})^{\leq 0}$.
Also the restriction of 
$\mcal{M}' \xar{\alpha} \mcal{M}|_{U}$
to $U_{i}$ is 
$\mcal{M}'_{i} \xar{\alpha_{i}} \mcal{M}|_{U_{i}}$.

Define $\mcal{M}''$ to be the cone of $\alpha$. So we have a 
distinguished triangle
\[ T = \Bigl( \mcal{M}' \xar{\alpha} \mcal{M}|_{U} \xar{\beta}
\mcal{M}'' \xar{\gamma} \mcal{M}'[1] \Bigr) \]
in $\msf{D}^{\mrm{b}}(\cat{Mod} \mcal{A}|_{U})$. By
Lemma \ref{lem4.5} there is a (unique) isomorphism
$\mcal{M}''|_{U_{i}} \iso \mcal{M}''_{i}$
such that $T|_{U_{i}} \iso T_{i}$
is an isomorphism of triangles. Therefore
$\mcal{M}''|_{U_{i}} \in 
{}^{p}\msf{D}^{\mrm{b}}_{\mrm{c}}
(\cat{Mod} \mcal{A}|_{U_i})^{\geq 1}$.
Using Lemma \ref{lem6.1} we see that 
$\mcal{M}'' \in 
{}^{p}\msf{D}^{\mrm{b}}_{\mrm{c}}
(\cat{Mod} \mcal{A}|_{U})^{\geq 1}$.

When $n>2$ the statement follows from induction and the case $n=2$.
\end{proof}

\begin{thm} \label{thm4.10}
Let $(X, \mcal{A})$ be a ringed space with a local collection 
of t-structures \linebreak $\{ \mfrak{T}_U \}_{U \in \mfrak{B}}$.
Assume $X$ is a noetherian topological space. 
Then the pair 
\[ \bigl( 
{}^{p}\msf{D}^{\mrm{b}}_{\mrm{c}}(\cat{Mod} \mcal{A})^{\leq 0}, 
{}^{p}\msf{D}^{\mrm{b}}_{\mrm{c}}(\cat{Mod} \mcal{A})^{\geq 0} 
\bigr) \]
from Definition \tup{\ref{dfn4.3}} is a t-structure on 
$\msf{D}^{\mrm{b}}_{\mrm{c}}(\cat{Mod} \mcal{A})$.
\end{thm}

\begin{proof}
Condition (i) of Definition \ref{dfn6.2} is trivially verified.

Let
$\mcal{M} \in {}^{p}\msf{D}^{\mrm{b}}_{\mrm{c}}
(\cat{Mod} \mcal{A})^{\leq 0}$
and
$\mcal{N} \in {}^{p}\msf{D}^{\mrm{b}}_{\mrm{c}}
(\cat{Mod} \mcal{A})^{\geq 1}$.
So 
$\mcal{N}[1] \in {}^{p}\msf{D}^{\mrm{b}}_{\mrm{c}}
(\cat{Mod} \mcal{A})^{\geq 0}$,
and by Lemma \ref{lem4.4}(3) with $i = -1$ we have
$\opn{Hom}_{\msf{D}(\cat{Mod} \mcal{A})}
(\mcal{M}, \mcal{N})  = 0$.
This verifies condition (ii).

It remains to prove condition (iii). Let a complex
$\mcal{M} \in \msf{D}^{\mrm{b}}_{\mrm{c}}(\cat{Mod} \mcal{A})$
be given. Choose a covering
$X = \bigcup_{i = 1}^{n} U_{i}$ by open sets in $\mfrak{B}$. 
For every $i$ the t-structure $\mfrak{T}_{U_i}$ gives rise to a 
distinguished triangle
\[ S_{i} = \bigl( 
\mcal{M}_{i}' \xar{\alpha_{i}} \mcal{M}|_{U_{i}} \xar{\beta_{i}} 
\mcal{M}_{i}'' \xar{\gamma_{i}} \mcal{M}_{i}'[1]
\bigr) \]
in $\msf{D}^{\mrm{b}}_{\mrm{c}}(\cat{Mod} \mcal{A}|_{U_{i}})$
with 
$\mcal{M}_{i}' \in 
{}^{p}\msf{D}^{\mrm{b}}_{\mrm{c}}
(\cat{Mod} \mcal{A}|_{U_{i}})^{\leq 0}$
and 
$\mcal{M}_{i}'' \in 
{}^{p}\msf{D}^{\mrm{b}}_{\mrm{c}}
(\cat{Mod} \mcal{A}|_{U_{i}})^{\geq 1}$.
By Lemma \ref{lem4.8} there is a triangle
\[ \mcal{M}' \xar{\alpha} \mcal{M} \xar{\beta} 
\mcal{M}'' \xar{\gamma} \mcal{M}[1] \]
in 
$\msf{D}^{\mrm{b}}_{\mrm{c}}(\cat{Mod} \mcal{A})$
whose restriction to each $U_i$ is isomorphic to $S_i$. Therefore 
by Lemma \ref{lem6.1} one has
$\mcal{M}' \in 
{}^{p}\msf{D}^{\mrm{b}}_{\mrm{c}}(\cat{Mod} \mcal{A})^{\leq 0}$
and 
$\mcal{M}'' \in 
{}^{p}\msf{D}^{\mrm{b}}_{\mrm{c}}(\cat{Mod} \mcal{A})^{\geq 1}$.
\end{proof}

\begin{thm} \label{thm4.11}
Let $(X, \mcal{A})$ be a ringed space with a local collection 
of t-structures \linebreak $\{ \mfrak{T}_U \}_{U \in \mfrak{B}}$.
Assume $X$ is a noetherian topological space. 
For any open set $V \subset X$ let
${}^{p}\msf{D}^{\mrm{b}}_{\mrm{c}}(\cat{Mod} \mcal{A}|_{V})^{0}$
be the heart of the t-structure from Theorem \tup{\ref{thm4.10}}. 
Then  
$V \mapsto$ \linebreak
${}^{p}\msf{D}^{\mrm{b}}_{\mrm{c}}(\cat{Mod} \mcal{A}|_{V})^{0}$
is a stack of subcategories of 
$\msf{D}^{\mrm{b}}_{\mrm{c}}(\cat{Mod} \mcal{A})$.
\end{thm}

\begin{proof}
Axiom (a) follows from Lemma \ref{lem6.1}. 
Axiom (c) is Lemma \ref{lem4.4}(2). Let us prove axiom (b). 
Suppose we are given open sets
$V = \bigcup_{i \in I} V_{i} \subset X$, 
complexes 
$\mcal{M}_{i} \in 
{}^{p}\msf{D}^{\mrm{b}}_{\mrm{c}}(\cat{Mod} \mcal{A}|_{V_{i}})^{0}$
and isomorphisms
$\phi_{i, j}: \mcal{M}_{i}|_{V_{i} \cap V_{j}} \to 
\mcal{M}_{j}|_{V_{i} \cap V_{j}}$
satisfying the cocycle condition. 
Since $X$ is noetherian, and in view of axiom (c), we may assume 
$I = \{ 1, \ldots, n \}$. Let us define
$W_i := \bigcup_{j = 1}^i V_j$. By induction on $i$ we will 
construct an object
$\mcal{N}_i \in 
{}^{p}\msf{D}^{\mrm{b}}_{\mrm{c}}(\cat{Mod} \mcal{A}|_{W_i})^{0}$
with isomorphisms
$\psi_{i, j} : \mcal{N}_i |_{V_j} \iso \mcal{M}_j$
for all $j \leq i$ that are compatible with the $\phi_{j, k}$. 
Then $\mcal{M} := \mcal{N}_n$ will be the desired global object on 
$V = W_n$. 

So assume $i < n$ and $\mcal{N}_i$ has already been defined. 
For any $j \leq i$ we have an isomorphism
\[ \phi_{j, i+1} \circ \psi_{i, j} : 
\mcal{N}_i|_{V_j \cap V_{i+1}} \iso 
\mcal{M}_j|_{V_j \cap V_{i+1}} \iso 
\mcal{M}_{i+1}|_{V_j \cap V_{i+1}} , \]
and these satisfy the cocycle condition. According to 
Lemma \ref{lem4.4}(2) there is an isomorphism
\[ \psi_{i, i+1} : 
\mcal{N}_i|_{W_i \cap V_{i+1}} \iso 
\mcal{M}_{i+1}|_{W_i \cap V_{i+1}} \] 
in 
${}^{p}\msf{D}^{\mrm{b}}_{\mrm{c}}
(\cat{Mod} \mcal{A}|_{W_i \cap V_{i+1}})^{0}$.
Denote by $f_{i+1} : W_i \to W_{i+1}$,
$g_{i+1} : V_{i+1} \to W_{i+1}$ and
$h_{i+1} : W_i \cap V_{i+1} \to W_{i+1}$ the inclusions. Define
$\mcal{N}_{i+1} \in 
\msf{D}(\cat{Mod} \mcal{A}|_{W_{i+1}})$
to be the cone of the morphism
\[ h_{(i+1) !}(\mcal{N}_{i}|_{W_i \cap V_{i+1}})
\xar{(\gamma, \psi_{i, i+1})} 
f_{(i+1) !} \mcal{N}_{i} \oplus g_{(i+1) !} \mcal{M}_{i+1} \] 
where $\gamma$ is the canonical morphism. We obtain a 
distinguished triangle
\[ h_{(i+1) !}(\mcal{N}_{i}|_{W_i \cap V_{i+1}}) \to
f_{(i+1) !} \mcal{N}_{i} \oplus g_{(i+1) !} \mcal{M}_{i+1} 
\to \mcal{N}_{i+1} \to
h_{(i+1) !}(\mcal{N}_{i}|_{W_i \cap V_{i+1}})[1] \]
in $\msf{D}(\cat{Mod} \mcal{A}|_{W_{i+1}})$.
Upon restriction to $W_i$ we get an isomorphism
$\mcal{N}_{i} \cong \mcal{N}_{i+1}|_{W_i}$; and upon restriction 
to $V_{i+1}$ we get an isomorphism
$\mcal{N}_{i+1}|_{V_{i+1}} \iso \mcal{M}_{i+1}$ which we call 
$\psi_{i+1, i+1}$. Finally from Lemma \ref{lem6.1} we see that
$\mcal{N}_{i+1} \in 
{}^{p}\msf{D}^{\mrm{b}}_{\mrm{c}}(\cat{Mod} 
\mcal{A}|_{W_{i+1}})^{0}$.
\end{proof}

\begin{rem}
In \cite{BBD} the authors considered t-structures on certain
subcategories of $\msf{D}(\cat{Mod} \k_X)$, for a topological 
space $X$ and a constant sheaf of rings $\k_X$. 
Perverse t-structures on 
$\msf{D}^{\mrm{b}}_{\mrm{c}}(\cat{Mod} \mcal{A})$
as above where $\mcal{A}$ is ``quasi-coherent''
have only been considered recently; see \cite{Bz}, 
\cite{Br}, \cite{Ka} and \cite{YZ4}.
\end{rem}

\section{Differential Quasi-Coherent Ringed Schemes
of Finite Type} \label{sec5}

In this section all the pieces of our puzzle come together, and we 
prove the main result Theorem \ref{thm0.2} -- repeated here as 
Theorem \ref{thm13.8}. As before $\k$ denotes the base field. All 
schemes are over $\k$, all rings are $\k$-algebras, all bimodules 
are central over $\k$, and all homomorphisms are over $\k$.

\begin{dfn} \label{dfn13.1}
Let $X$ be a finite type $\k$-scheme and let 
$\mcal{A}$ be a quasi-coherent $\mcal{O}_{X}$-ring. A 
{\em differential quasi-coherent $\mcal{O}_{X}$-filtration
of finite type} on $\mcal{A}$ is an ascending filtration 
$F = \{ F_i \mcal{A} \}_{i \in \mbb{Z}}$ by subsheaves with 
the following properties:
\begin{enumerate}
\rmitem{i} Each $F_i \mcal{A}$ is an $\mcal{O}_{X}$-sub-bimodule
of $\mcal{A}$, quasi-coherent on both sides.
\rmitem{ii} $F_{-1} \mcal{A} = 0$ and 
$\mcal{A} = \bigcup F_i \mcal{A}$.
\rmitem{iii} $1 \in F_0 \mcal{A}$ and
$(F_i \mcal{A}) \cdot (F_j \mcal{A}) \subset F_{i+j} \mcal{A}$.
\rmitem{iv} The graded sheaf of rings 
$\opnt{gr}^F \mcal{A}$ is an $\mcal{O}_{X}$-algebra.
\rmitem{v} The center $\mrm{Z}(\opnt{gr}^F \mcal{A})$ is a finite 
type quasi-coherent $\mcal{O}_{X}$-algebra.
\rmitem{vi} $\opnt{gr}^F \mcal{A}$ is a coherent
$\mrm{Z}(\opnt{gr}^F \mcal{A})$-module.
\end{enumerate}
\end{dfn}

By properties (i) and (iii) we get a ring homomorphism
$\mcal{O}_{X} \to \opnt{gr}^F \mcal{A}$. Property (iv) tells us 
that the image of $\mcal{O}_{X}$ is inside 
$\mrm{Z}(\opnt{gr}^F \mcal{A})$. By (v-vi) we see that
$\Gamma(U, \opnt{gr}^F \mcal{A})$ is a noetherian ring for any 
affine open set $U$, so $(X, \mcal{A})$ is a noetherian 
quasi-coherent ringed scheme. 

Observe that the definition is left-right symmetric: if $\mcal{A}$ 
is a differential quasi-coherent $\mcal{O}_{X}$-ring of finite 
type then so is $\mcal{A}^{\mrm{op}}$. 
The name ``differential filtration'' signifies the resemblance to 
Grothendieck's definition of differential operators in 
\cite{EGA-IV}. 

\begin{dfn} \label{dfn13.3}
Let $X$ be a finite type $\k$-scheme. A 
{\em differential quasi-coherent $\mcal{O}_{X}$-ring of finite 
type} is an $\mcal{O}_{X}$-ring $\mcal{A}$ that admits some 
differential quasi-coherent $\mcal{O}_{X}$-filtration of 
finite type. The pair $(X, \mcal{A})$ is then called a 
{\em differential quasi-coherent ringed $\k$-scheme finite type}. 
\end{dfn}

Let us recall a couple of definition from \cite{YZ4}.

\begin{dfn}[{\cite[Definition 2.1]{YZ4}}]
Suppose $C$ is a commutative $\k$-algebra and $A$ 
is a $C$-ring. A {\em differential $C$-filtration} on
$A$ is a filtration $F = \{ F_{i} A \}_{i \in \mbb{Z}}$ with the 
following properties:
\begin{enumerate}
\rmitem{i} Each $F_i A$ is a $C$-sub-bimodule. 
\rmitem{ii} $F_{-1} A = 0$ and $A = \bigcup F_{i} A$.
\rmitem{iii} $1 \in F_0 A$ and 
$F_i A \cdot F_j A \subset F_{i + j} A$.
\rmitem{iv} The graded ring $\opnt{gr}^{F} A$ 
is a $C$-algebra.
\end{enumerate}
A is called a {\em differential $C$-ring} if it admits 
some differential $C$-filtration.
\end{dfn}

\begin{dfn}[{\cite[Definition 2.2]{YZ4}}]
Let $C$ be a commutative noetherian $\k$-algebra, 
and let $A$ be a $C$-ring. 
\begin{enumerate}
\item A {\em differential $C$-filtration of finite type} on $A$ is 
a differential $C$-filtration $F = \{ F_{i} A \}$ such that the 
graded $C$-algebra $\opnt{gr}^{F} A$ is a finitely generated
module over its center $\mrm{Z}(\opnt{gr}^{F} A)$, and 
$\mrm{Z}(\opnt{gr}^{F} A)$  is a finitely generated $C$-algebra.
\item We say $A$ is a {\em differential $C$-ring of finite type} 
if it admits some differential $C$-filtration of finite type. 
\item If $A$ is a differential $\k$-ring of finite type then 
we also call $A$ a {\em differential $\k$-algebra of finite type}.
\end{enumerate}
\end{dfn}

The next proposition says that a differential $C$-ring of finite 
type is just the ring theoretic counterpart of a differential 
quasi-coherent ringed $\k$-scheme finite type.

\begin{prop} \label{prop13.6}
Let $C$ be a noetherian $\k$-algebra and $U := \opn{Spec} C$.
\begin{enumerate}
\item Given a differential quasi-coherent
$\mcal{O}_{U}$-ring of finite type $\mcal{A}$ the ring
$A := \Gamma(U, \mcal{A})$ is a differential $C$-ring of finite 
type. 
\item Given a differential $C$-ring of finite type $A$
there is a differential quasi-coherent
$\mcal{O}_{U}$-ring of finite type
$\mcal{A}$, unique up to isomorphism, such that 
$\Gamma(U, \mcal{A}) \cong A$ as $C$-rings.
\end{enumerate}
\end{prop}

\begin{proof}
The proof of (1) is straightforward. For (2) 
use \cite[Proposition 5.17 and Corollary 5.13]{YZ4}. 
As filtration on $\mcal{A}$ we may take 
$F_{i} \mcal{A} := \mcal{O}_{U} \otimes_{C} F_{i} A$
where $\{ F_{i} A \}$ is any differential $C$-filtration 
of finite type on $A$.
\end{proof}

\begin{prop} \label{prop13.7}
Suppose $(X, \mcal{A})$ and $(Y, \mcal{B})$ are two differential 
quasi-coherent ringed schemes of finite type over $\k$. Then 
the product 
$(X \times Y, \mcal{A} \boxtimes \mcal{B})$ exists 
\tup{(}cf.\ Definition \tup{\ref{dfn6.2})}, and it too is a 
differential quasi-coherent ringed scheme of finite type 
over $\k$.
\end{prop}

\begin{proof}
Let $A, B, C$ be as in condition (ii) of Theorem \ref{thm4.2}. 
Then \cite[Proposition 2.9]{YZ4} tells us that $A \otimes B$
is a differential $C$-ring of finite type. Using
\cite[Proposition 5.14 and Corollary 5.13]{YZ4}
we see that condition (ii) of Theorem \ref{thm4.2} is 
satisfied. Therefore the product 
$(X \times Y, \mcal{A} \boxtimes \mcal{B})$ exists. 

Given differential filtrations of $\mcal{A}$ and $\mcal{B}$
one can construct a filtration on 
$\mcal{A} \boxtimes \mcal{B}$ by gluing 
together the affine filtrations described in 
the proof of \cite[Proposition 2.9]{YZ4}.
\end{proof}

\begin{dfn} \label{dfn8.3}
Let $A$ be a noetherian $\k$-algebra with rigid dualizing 
complex $R_A$, and let
\[ \mrm{D} := \opn{RHom}_A(-, R_A) : 
\msf{D}^{\mrm{b}}_{\mrm{f}}(\cat{Mod} A) \to
\msf{D}^{\mrm{b}}_{\mrm{f}}(\cat{Mod} A^{\mrm{op}}) \]
be the induced duality. The {\em rigid perverse t-structure} on
$\msf{D}^{\mrm{b}}_{\mrm{f}}(\cat{Mod} A)$ is defined by
\[ \begin{aligned}
{}^{p}\msf{D}^{\mrm{b}}_{\mrm{f}}(\cat{Mod} A)^{\leq 0} & :=
\{ M \in \msf{D}^{\mrm{b}}_{\mrm{f}}(\cat{Mod} A) \mid
\mrm{H}^i \mrm{D} M = 0 \text{ for all } i < 0 \} , \\
{}^{p}\msf{D}^{\mrm{b}}_{\mrm{f}}(\cat{Mod} A)^{\geq 0} & :=
\{ M \in \msf{D}^{\mrm{b}}_{\mrm{f}}(\cat{Mod} A) \mid
\mrm{H}^i \mrm{D} M = 0 \text{ for all } i > 0 \} .
\end{aligned} \]
An object 
$M \in {}^{p}\msf{D}^{\mrm{b}}_{\mrm{f}}(\cat{Mod} A)^{0}$
is called a {\em perverse $A$-module}.
\end{dfn}

The rigid perverse t-structure on 
$\msf{D}^{\mrm{b}}_{\mrm{f}}(\cat{Mod} A)$
is of course dual to the standard t-structure on 
$\msf{D}^{\mrm{b}}_{\mrm{f}}(\cat{Mod} A^{\mrm{op}})$, via the 
duality $\mrm{D}$. 

Suppose $(X, \mcal{A})$ is a differential quasi-coherent 
ringed scheme of finite type over $\k$. By Proposition 
\ref{prop13.6}(1), for every affine open set
$U \subset X$ the $\k$-algebra $\Gamma(U, \mcal{A})$ is a 
differential $\k$-algebra of finite type. Hence according to 
\cite[Theorem 8.1]{YZ4} the algebra $\Gamma(U, \mcal{A})$ 
has a rigid dualizing complex, and the rigid perverse t-structure 
on
$\msf{D}^{\mrm{b}}_{\mrm{f}}(\cat{Mod} \mcal{A})$
exists. 

\begin{dfn}
Let $(X, \mcal{A})$ be a differential quasi-coherent 
ringed scheme of finite type over $\k$. Taking ``$\star$'' to be 
either ``$\leq 0$'', ``$\geq 0$'' or ``$0$'', we define 
classes of objects
\[ \begin{aligned}
{}^{p}\msf{D}^{\mrm{b}}_{\mrm{c}}(\cat{Mod} \mcal{A})^{\star} 
& := \big\{ \mcal{M} \in  
\msf{D}^{\mrm{b}}_{\mrm{c}}(\cat{Mod} \mcal{A}) \mid
\mrm{R} \Gamma(U, \mcal{M}) \in 
{}^{p}\msf{D}^{\mrm{b}}_{\mrm{f}}
\big( \cat{Mod} \Gamma(U, \mcal{A}) \big)^{\star} \\
& \hspace{33ex} \text{ for all affine open sets } U \big\} .
\end{aligned} \]
\end{dfn}

\begin{thm} \label{thm9.1}
Let $(X, \mcal{A})$ be a differential quasi-coherent 
ringed sch\-eme of finite type over $\k$. Then: 
\begin{enumerate}
\item The pair 
\[ \bigl( 
{}^{p}\msf{D}^{\mrm{b}}_{\mrm{c}}(\cat{Mod} \mcal{A})^{\leq 0}, 
{}^{p}\msf{D}^{\mrm{b}}_{\mrm{c}}(\cat{Mod} \mcal{A})^{\geq 0} 
\bigr) \]
is a t-structure on 
$\msf{D}^{\mrm{b}}_{\mrm{c}}(\cat{Mod} \mcal{A})$.
\item The assignment 
$V \mapsto 
{}^{p}\msf{D}^{\mrm{b}}_{\mrm{c}}(\cat{Mod} \mcal{A}|_{V})^{0}$,
for $V \subset X$ open, is a stack of abelian categories on $X$.
\end{enumerate} 

\end{thm}

\begin{proof}
For an affine open set $U \subset X$ let
\[ \mfrak{T}_U := \bigl( 
{}^{p}\msf{D}^{\mrm{b}}_{\mrm{c}}(\cat{Mod} \mcal{A}|_U)^{\leq 0}, 
{}^{p}\msf{D}^{\mrm{b}}_{\mrm{c}}(\cat{Mod} \mcal{A}|_U)^{\geq 0} 
\bigr) . \]
By Theorems \ref{thm4.10} and \ref{thm4.11} 
it suffices to prove that 
$\{ \mfrak{T}_U \}$ is a local collection of t-structures on 
$(X, \mcal{A})$. Namely, 
given an affine open set $U$, a complex
$\mcal{M} \in 
\msf{D}^{\mrm{b}}_{\mrm{c}}(\cat{Mod} \mcal{A}|_U)$,
and an affine open covering $U = \bigcup U_k$, one has
$\mcal{M} \in 
{}^p\msf{D}^{\mrm{b}}_{\mrm{c}}(\cat{Mod} \mcal{A}|_U)^{\star}$
if and only if
$\mcal{M}|_{U_k} \in 
{}^p\msf{D}^{\mrm{b}}_{\mrm{c}}(\cat{Mod} \mcal{A}|_{U_k})^{\star}$
for all $k$. 

Write 
$A := \Gamma(U, \mcal{A})$, 
$M := \mrm{R} \Gamma(U, \mcal{M})$,
$A_{k} := \Gamma(U_{k}, \mcal{A})$ and
$M_{k} := \mrm{R} \Gamma(U_{k}, \mcal{M})$.
We must show that 
$M \in 
{}^{p}\msf{D}^{\mrm{b}}_{\mrm{f}}(\cat{Mod} A)^{\star}$
if and only if 
$M_k \in 
{}^{p}\msf{D}^{\mrm{b}}_{\mrm{f}}(\cat{Mod} A_k)^{\star}$
for all $k$.

Denote by $R_A$ and $R_{A_k}$ the rigid dualizing 
complexes of the rings $A$ and $A_k$ respectively. By 
\cite[Corollary 5.20]{YZ4} $A_k$ is a localization of $A$. 
According to 
\cite[Theorem 8.14 and Propositions 5.21 and 5.17]{YZ4}
each $A$-bimodule 
$\mrm{H}^i R_A$ is evenly localizable to $A_k$. Hence 
\cite[Theorem 6.2]{YZ4} tells us that 
$R_{A_k} \cong R_A \otimes_{A} A_k$
in $\msf{D}(\cat{Mod}\, A \otimes A_k^{{\mrm{op}}})$, 
and likewise 
$R_{A_k} \cong A_k \otimes_{A} R_A$
in $\msf{D}(\cat{Mod}\, A_k \otimes A^{{\mrm{op}}})$.

Define complexes 
$N := \opn{RHom}_{A}(M, R_{A})$
and
$N_k := \opn{RHom}_{A_k}(M_k, R_{A_k})$. 
According to Lemma \ref{lem5.4} we have 
$M_k \cong A_k \otimes_{A} M$. By \cite[Lemma 3.7]{YZ3}
\[ \begin{aligned}
N \otimes_{A} A_k 
& \cong \opn{RHom}_{A}(M, R_{A} \otimes_{A} A_k) 
\cong \opn{RHom}_{A}(M, R_{A_k}) \\
& \cong \opn{RHom}_{A_k}(A_k \otimes_{A} M, R_{A_k}) 
\cong \opn{RHom}_{A_k}(M_k, R_{A_k}) = N_k . 
\end{aligned} \]

Let us consider the case where ``$\star$'' is ``$\geq$''.
Suppose 
$M \in
{}^{p}\msf{D}^{\mrm{b}}_{\mrm{f}}(\cat{Mod} A)^{\geq 0}$.
By definition of the perverse rigid t-structure we have
$N \in 
\msf{D}^{\mrm{b}}_{\mrm{f}}(\cat{Mod} A^{\mrm{op}})^{\leq 0}$,
i.e.\ $\mrm{H}^j N = 0$ for all $j > 0$. 
Hence
$\mrm{H}^j N_k \cong (\mrm{H}^j N) \otimes_{A} A_k = 0$,
implying that 
$M_k \in 
\msf{D}^{\mrm{b}}_{\mrm{f}}(\cat{Mod} A_k)^{\geq 0}$. 

Conversely suppose 
$M_k \in
{}^{p}\msf{D}^{\mrm{b}}_{\mrm{f}}(\cat{Mod} A_k)^{\geq 0}$
for all $k$. We may assume that 
$U = \bigcup_{k = 1}^n U_k$. Let 
$C := \Gamma(U, \mcal{O}_{X})$
and
$C_k := \Gamma(U_k, \mcal{O}_{X})$.
The ring homomorphism 
$C \to \prod_{k = 1}^{n} C_{k}$
is faithfully flat (cf.\ \cite[Proposition 5.6]{YZ4}). 
Applying $A \otimes_C -$ it follows that 
$A \to \prod_{k = 1}^{n} A_{k}$
is faithfully flat (on both sides), so we get an injection 
\[ \mrm{H}^{j} N \inj \boplus_{k} (\mrm{H}^{j} N) 
\otimes_{A} A_{k} . \]
As above we conclude that 
$N \in \msf{D}^{\mrm{b}}_{\mrm{f}}(\cat{Mod} 
A^{\mrm{op}})^{\leq 0}$, and hence
$M \in {}^{p}\msf{D}^{\mrm{b}}_{\mrm{f}}(\cat{Mod} A)^{\geq 0}$. 

The case where ``$\star$'' is ``$\leq 0$'' is handled similarly. 
The only difference is that we have to verify the vanishing of 
$\mrm{H}^{j} N$ and $\mrm{H}^{j} N_k$ for $j < 0$. 
\end{proof}

\begin{dfn} \label{dfn9.4}
The t-structure on 
$\msf{D}^{\mrm{b}}_{\mrm{c}}(\cat{Mod} \mcal{A})$
in Theorem \ref{thm9.1} is called the 
{\em rigid perverse t-structure}. An object
$\mcal{M} \in 
{}^{p}\msf{D}^{\mrm{b}}_{\mrm{c}}(\cat{Mod} \mcal{A})^{0}$
is called a {\em perverse coherent $\mcal{A}$-module}. 
\end{dfn}

\begin{thm} \label{thm13.8}
Let $(X, \mcal{A})$ be a separated differential 
quasi-coherent ringed \linebreak
scheme of finite type over $\k$.
Then $\mcal{A}$ has a rigid dualizing complex 
$(\mcal{R}_{\mcal{A}}, \bsym{\rho})$. It is unique up to a unique 
isomorphism in 
$\msf{D}^{\mrm{b}}_{\mrm{c}}(\cat{Mod} \mcal{A}^{\mrm{e}})$.
\end{thm}

\begin{proof}
By Proposition \ref{prop13.7} the product 
$(X^2, \mcal{A}^{\mrm{e}})$
is also a separated differential quasi-coherent ringed scheme 
of finite type over $\k$. Hence by Theorem \ref{thm9.1}
we obtain the rigid perverse t-structure on
$\msf{D}^{\mrm{b}}_{\mrm{c}} (\cat{Mod} \mcal{A}^{\mrm{e}})$,
and there is a stack of abelian categories
$W \mapsto 
{}^p\msf{D}^{\mrm{b}}_{\mrm{c}}
(\cat{Mod} \mcal{A}^{\mrm{e}}|_{W})^0$
on $X^2$.

Let us choose, for ease of notation, 
an indexing $\{ U_i \}$ of the set
$\cat{Aff} X$ of affine open sets of $X$. 
For any index $i$ let
$A_{i}:= \Gamma(U_{i}, \mcal{A})$, which is a differential 
$\k$-algebra of finite type. Then $A_i$ has a rigid dualizing 
complex $(R_i, \rho_i)$. By \cite[Theorem 8.9]{YZ4} the 
complex $R_i$ is a perverse bimodule, i.e.\
$R_i \in {}^{p}\msf{D}^{\mrm{b}}_{\mrm{f}}
(\cat{Mod} A_{i}^{\mrm{e}})^{0}$. 

For a pair of indices $i, j$ let
$A_{(i, j)}:= \Gamma(U_i \cap U_j, \mcal{A})$. 
Define 
\[ R_{i \to (i, j)} := 
A_{(i, j)} \otimes_{A_i} R_i \otimes_{A_i} A_{(i, j)} 
\in \msf{D}(\cat{Mod} A_{(i, j)}^{\mrm{e}}) , \]
and let
$\mrm{q}_{A_{(i, j)} / A_i} : R_{A_i} \to R_{i \to (i, j)}$
be the morphism with formula $r \mapsto 1 \otimes r \otimes 1$.
According to \cite[Theorem 6.2]{YZ4} the complex 
$R_{i \to (i, j)}$ is dualizing over the ring $A_{(i, j)}$, and it 
has a unique rigidifying isomorphism 
$\rho_{i \to (i, j)}$ such that 
$\mrm{q}_{A_{(i, j)} / A_i}$ is a rigid localization morphism. 
Likewise we obtain a rigid dualizing complex 
$(R_{j \to (i, j)}, \rho_{j \to (i, j)})$.
By \cite[Corollary 4.3]{YZ2} there is a unique isomorphism
\[ \phi_{(i, j)} : R_{i \to (i, j)} \iso
R_{j \to (i, j)} \]
in $\msf{D}(\cat{Mod} A_{(i, j)}^{\mrm{e}})$ that's a rigid trace. 
The isomorphisms $\phi_{(i, j)}$ will then satisfy the cocycle 
condition in
$\msf{D}(\cat{Mod} A_{(i, j, k)}^{\mrm{e}})$, where
$A_{(i, j, k)} := \Gamma(U_i \cap U_j \cap U_k, \mcal{A})$ 
for a triple intersection.

Now consider the affine ringed scheme
$(U_{i}^2, \mcal{A}^{\mrm{e}}|_{U_{i}^2})$.
Let us denote by 
$\mcal{R}_{i} \in \msf{D}^{\mrm{b}}_{\mrm{c}}
(\cat{Mod} \mcal{A}^{\mrm{e}}|_{U_{i}^2})$ 
the sheafification of $R_i$. By definition of the rigid 
perverse t-structure on 
$\msf{D}^{\mrm{b}}_{\mrm{c}}
(\cat{Mod} \mcal{A}^{{\mrm{e}}}|_{U_{i}^2})$
we have 
$\mcal{R}_{i} \in {}^{p}\msf{D}^{\mrm{b}}_{\mrm{c}}
(\cat{Mod} \mcal{A}^{\mrm{e}}|_{U_{i}^2})^{0}$.
By Lemma \ref{lem5.4} we obtain induced isomorphisms
\[ \phi_{(i, j)} : 
\mcal{R}_{i}|_{U_{i}^2 \cap U_{j}^2} 
\iso \mcal{R}_{j}|_{U_{i}^2 \cap U_{j}^2} \]
in 
${}^{p}\msf{D}^{\mrm{b}}_{\mrm{c}}
(\cat{Mod} \mcal{A}^{\mrm{e}}
|_{U_{i}^2 \cap U_{j}^2})^{0}$,
and these satisfy the cocycle condition on triple intersections.
Let $V := X^2 - \Delta(X)$. Since
$\mcal{R}_{i}|_{U_{i}^2 \cap V} = 0$
we have sufficient gluing data corresponding to the open covering
$X^2 = \bigl( \bigcup\nolimits_i U_{i}^2 \bigr) \cup V$,
and by Theorem \ref{thm9.1} we deduce the existence and uniqueness
of a global complex 
$\mcal{R}_{\mcal{A}} \in 
{}^{p}\msf{D}^{\mrm{b}}_{\mrm{c}}
(\cat{Mod} \mcal{A}^{\mrm{e}})^{0}$
together with isomorphisms
$\phi_i : \mcal{R}_{\mcal{A}}|_{U_i} \iso \mcal{R}_i$.
By Theorem \ref{thm7.2} and Proposition \ref{prop7.3} 
$\mcal{R}_{\mcal{A}}$ is a local dualizing complex over 
$\mcal{A}$. And by construction $\mcal{R}_{\mcal{A}}$ comes 
equipped with a collection 
$\bsym{\rho} = \{ \rho_i \}$ of rigidifying isomorphisms 
that is compatible with the sheaf structure. 
\end{proof}

Finite morphisms between ringed schemes were defined in 
Definition \ref{dfn7.6}. A ring homomorphism $A \to B$ is called 
{\em centralizing} if $B = \sum A b_i$ for elements
$b_i \in B$ that commute with $A$.

\begin{dfn} \label{dfn13.4}
Let $f: (Y, \mcal{B}) \to (X, \mcal{A})$ be a finite
morphism between two noetherian quasi-coherent ringed 
schemes over $\k$. We say $f$ is {\em finite centralizing} 
if for any affine open set $U \subset X$ the finite ring 
homomorphism 
$f^* : \Gamma(U, \mcal{A}) \to \Gamma(U, f_* \mcal{B})$
is centralizing. 
\end{dfn}

\begin{exa}
Suppose $\mcal{I} \subset \mcal{A}$ is a coherent sheaf of 
(two-sided) ideals, and define 
$\mcal{B} :=  \mcal{A} /  \mcal{I}$.
Then $(X,  \mcal{B})$ is a differential 
quasi-coherent ringed scheme of finite type over $\k$, and 
$(X,  \mcal{B}) \to (X,  \mcal{A})$
is a finite centralizing morphism.
\end{exa}

\begin{prop} \label{prop13.1}
Let $f : (Y, \mcal{B}) \to (X, \mcal{A})$
be a finite centralizing morphism between two separated differential 
quasi-coherent ringed schemes of finite type over $\k$. Consider 
the the rigid perverse t-structures on these ringed schemes.
\begin{enumerate}
\item Let 
$\mcal{M} \in \msf{D}^{\mrm{b}}_{\mrm{c}}(\cat{Mod} \mcal{B})$. 
Then 
$\mcal{M} \in 
{}^{p}\msf{D}^{\mrm{b}}_{\mrm{c}}(\cat{Mod} \mcal{B})^0$
if and only if
$\mrm{R} f_* \mcal{M} \in 
{}^{p}\msf{D}^{\mrm{b}}_{\mrm{c}}(\cat{Mod} \mcal{A})^0$.
\item Assume $f^* : \mcal{A} \to f_* \mcal{B}$ is surjective. Then 
the functor
\[ \mrm{R} f_* : 
{}^{p}\msf{D}^{\mrm{b}}_{\mrm{c}}(\cat{Mod} \mcal{B})^0
\to {}^{p}\msf{D}^{\mrm{b}}_{\mrm{c}}(\cat{Mod} \mcal{A})^0 \]
is fully faithful. 
\end{enumerate}
\end{prop}

\begin{proof}
(1) Let $U \subset X$ be an affine open set and $V := f^{-1}(U)$. 
Define $A := \Gamma(U, \mcal{A})$, $B := \Gamma(V, \mcal{B})$ and
$M := \mrm{R} \Gamma(V, \mcal{M})$. 
By Lemma \cite[Lemma 6.12]{YZ5} it is enough to show that 
$M \in {}^{p}\msf{D}^{\mrm{b}}_{\mrm{c}}(\cat{Mod} B)^0$
if and only if
$\opn{rest}_{B / A} M \in 
{}^{p}\msf{D}^{\mrm{b}}_{\mrm{c}}(\cat{Mod} A)^0$. 
This is done in \cite[Proposition 8.3(1)]{YZ4}.

\medskip \noindent
(2) This follows from Theorem \ref{thm9.1} and 
\cite[Proposition 5.3(2)]{YZ4}. 
\end{proof}

Rigid trace morphisms between rigid dualizing complexes
were defined in Definition \ref{dfn7.7}.

\begin{thm}
\label{thm13.9}
Let $(X, \mcal{A})$ and $(Y, \mcal{B})$ be two separated
differential quasi-coherent ringed schemes of finite type over $\k$, 
and let $f: (Y, \mcal{B}) \to (X, \mcal{A})$
be a finite centralizing morphism. Then there exists a unique
rigid trace 
\[ \opn{Tr}_{f} : \mrm{R} f^{\mrm{e}}_{*} \mcal{R}_{\mcal{B}}
\to \mcal{R}_{\mcal{A}} . \]
\end{thm}

\begin{proof}
The morphism
$f^{\mrm{e}} : (Y^2, \mcal{B}^{\mrm{e}}) \to 
(X^2, \mcal{A}^{\mrm{e}})$
is also finite centralizing, so by Proposition \ref{prop13.1} 
we get
$\mrm{R} f^{\mrm{e}}_{*} \mcal{R}_{\mcal{B}} \in
{}^{p}\msf{D}^{\mrm{b}}_{\mrm{c}}(\cat{Mod} \mcal{A}^{\mrm{e}})^0$.
Hence by Theorem \ref{thm9.1} any morphism
$\mrm{R} f^{\mrm{e}}_{*} \mcal{R}_{\mcal{B}}
\to \mcal{R}_{\mcal{A}}$
is determined locally. 

Choose an affine open covering
$X = \bigcup_{i = 1}^{n} U_i$ and let
$A_i := \Gamma(U_i, \mcal{A})$ and
$B_i := \Gamma(f^{-1}(U_i), \mcal{B}) \cong
\Gamma(U_i, f_* \mcal{B})$.
Since $A_i \to B_i$ is a finite centralizing homomorphism, 
\cite[Proposition 8.2]{YZ4} asserts the existence of the trace
$\opn{Tr}_{A_i / B_i}: R_{B_i} \to R_{A_i}$
in $\msf{D}^{\mrm{b}}_{\mrm{f}}(\cat{Mod} A^{\mrm{e}})$.
The uniqueness of $\opn{Tr}_{A_i / B_i}$ is always true. By
Corollary \ref{cor5.5} we get a morphism
\[ \opn{Tr}_i: 
(\mrm{R} f^{\mrm{e}}_{*} \mcal{R}_{\mcal{B}})|_{U_i^2}
\to \mcal{R}_{\mcal{A}}|_{U_i^2} \]
in 
${}^{p}\msf{D}^{\mrm{b}}_{\mrm{c}}
(\cat{Mod}\, (\mcal{A}^{\mrm{e}})|_{U_i^2})^0$. 

By \cite[Proposition 6.3]{YZ4} the rigid trace localizes.
Therefore the morphisms $\opn{Tr}_i$ coincide on intersections
$U_i^2 \cap U_j^2$. Both 
$\mrm{R} f^{\mrm{e}}_{*} \mcal{R}_{\mcal{B}}$
and 
$\mcal{R}_{\mcal{A}}$
are supported on the diagonal $\Delta(X)$. Therefore we have 
gluing data for a global morphism 
\[ \opn{Tr}_{f} : \mrm{R} f^{\mrm{e}}_{*} \mcal{R}_{\mcal{B}}
\to \mcal{R}_{\mcal{A}}  \]
in
${}^p\msf{D}^{\mrm{b}}_{\mrm{c}}(\cat{Mod} \mcal{A}^{\mrm{e}})^0$
as required. 

The uniqueness of $\opn{Tr}_{f}$ is a consequence of the fact that 
${}^p\msf{D}^{\mrm{b}}_{\mrm{c}}(\cat{Mod} \mcal{A}^{\mrm{e}})^0$
is a stack and the uniqueness of the traces $\opn{Tr}_i$.
\end{proof}

\begin{exa}
If $X$ is a finite type $\k$-scheme and 
$\mcal{A} = \mcal{O}_{X}$ then conditions (i)-(iii) in Theorem 
\ref{thm9.1} hold. Hence
$\msf{D}^{\mrm{b}}_{\mrm{c}}(\cat{Mod} \mcal{O}_X)$
has the rigid perverse t-structure. In Section \ref{sec6} we show 
that the perverse coherent $\mcal{O}_X$-modules are nothing but 
the Cohen-Macaulay complexes in 
$\msf{D}^{\mrm{b}}_{\mrm{c}}(\cat{Mod} \mcal{A})$.
\end{exa}

\begin{rem}
Suppose $(X, \mcal{A})$ is a noetherian quasi-coherent ringed 
scheme. In his recent paper \cite{Ka} Kashiwara proves that 
any bounded filtration of $X$ by families of 
supports induces a perverse t-structure on 
$\msf{D}^{\mrm{b}}_{\mrm{qc}}(\cat{Mod} \mcal{A})$. In particular 
this is true for the filtration by codimension 
(coniveau) $\mfrak{S} = \{ \mfrak{S}^d \}$.

When $\k := \mbb{C}$, $X$ is smooth and $\mcal{A} := \mcal{D}_X$, 
Kashiwara proves that the t-structure induced by 
$\mfrak{S}$ on the category 
$\msf{D}^{\mrm{b}}_{\mrm{rh}}(\cat{Mod} \mcal{D}_X)$
of regular holonomic complexes is dual to the standard 
t-structure on 
$\msf{D}^{\mrm{b}}_{\mrm{c}}(\cat{Mod} \mbb{C}_{X_{\mrm{an}}})$,
via the Riemann-Hilbert correspondence 
$\mrm{R} \mcal{H}om_{\mcal{D}_X}(-, \mcal{O}_{X_{\mrm{an}}})$.

On the other hand if $\mcal{A}$ is a coherent 
$\mcal{O}_X$-algebra the t-structure induced by $\mfrak{S}$ 
on $\msf{D}^{\mrm{b}}_{\mrm{c}}(\cat{Mod} \mcal{A})$
is the rigid perverse t-structure, as can be seen from 
\cite[Theorem 0.6]{YZ4}. 
\end{rem}

\section{Examples and Complements} \label{sec6}

In this section we present several examples of differential 
quasi-coherent ringed schemes of finite type over $\k$ 
and their rigid dualizing complexes,
including the the commutative case 
$\mcal{A} = \mcal{O}_X$. Finally we discuss in detail the ringed 
scheme $(\mbf{A}^1, \mcal{D}_{\mbf{A}^1})$ and some of its 
partial compactifications. 

Let $X$ be a separated $\k$-scheme, and 
consider the diagonal embedding
$\Delta :  X \to X^2$. A {\em central 
rigid dualizing complex} over $X$ is a pair
$(\mcal{R}_X, \bsym{\rho}_X)$, where 
$\mcal{R}_X \in \msf{D}^{\mrm{b}}_{\mrm{c}}(\cat{Mod} \mcal{O}_X)$
and 
$(\mrm{R} \Delta_* \mcal{R}_X, \bsym{\rho}_X)$
is a rigid dualizing complex in the sense of Definition 
\ref{dfn7.5}. 

\begin{prop}
The scheme $X$ admits a 
central rigid dualizing complex, which is unique up to a unique 
isomorphism in $\msf{D}(\cat{Mod} \mcal{O}_X)$.
\end{prop}

\begin{proof}
See \cite[Theorem 0.1]{YZ5}. 
This can also be proved with the methods of Section \ref{sec5}.
\end{proof}

In case the structural morphism $\pi : X \to \opn{Spec} \k$ is 
embeddable then there is an isomorphism 
$\mcal{R}_X \cong \pi^! \k$, where $\pi^!$ is the twisted inverse 
image from \cite{RD}; see \cite[Proposition 6.18]{YZ5}.

Recall that a complex
$\mcal{M} \in \msf{D}^{\mrm{b}}(\cat{Mod} \mcal{O}_X)$ 
is called {\em Cohen-Macaulay} (with respect to the filtration by 
dimension) if for every point $x$ the local 
cohomologies $\mrm{H}^i_x \mcal{M}$ all vanish except for
$i = -\opn{dim} \overline{\{ x \}}$. Equivalently, $\mcal{M}$ is 
Cohen-Macaulay if there is an isomorphism
$\mcal{M} \cong \mrm{E} \mcal{M}$ in 
$\msf{D}(\cat{Mod} \mcal{O}_X)$, where $\mrm{E} \mcal{M}$ is the 
Cousin complex. See \cite[Section IV.3]{RD} or
\cite[Theorem 2.11]{YZ3}.

\begin{thm} \label{thm13.4}
Let $X$ be a separated finite type scheme over $\k$, let 
$\mcal{R}_X$ be the central rigid dualizing complex of $X$, and 
let $\mrm{D}_X$ be the auto-duality functor
$\mrm{R} \mcal{H}om_{\mcal{O}_X}(-, \mcal{R}_X)$. 
Then the following conditions are equivalent for 
$\mcal{M} \in
\msf{D}^{\mrm{b}}_{\mrm{c}}(\cat{Mod} \mcal{O}_X)$.
\begin{enumerate}
\rmitem{i} $\mcal{M}$ is a perverse coherent sheaf
\tup{(}for the rigid perverse t-structure\tup{)}. 
\rmitem{ii} $\mrm{D} \mcal{M}$ is a coherent sheaf, i.e.\ 
$\mrm{H}^i \mrm{D} \mcal{M} = 0$ for all $i \neq 0$.
\rmitem{iii} $\mcal{M}$ is a Cohen-Macaulay complex.
\end{enumerate} 
\end{thm}

\begin{proof}
All three conditions can be checked locally, so we may assume $X$ 
is affine. Then (i) and (ii) are equivalent by definition. 

Since the dualizing complex $\mcal{R}_X$ is adjusted to Krull 
dimension it follows that the Cousin complex
$\mcal{K}_{X} := \mrm{E} \mcal{R}_X$ is a residual complex.
Thus for every $i$ we have
$\mcal{K}_{X}^{-i} = \bigoplus \mcal{K}_{X}(x)$,
where $x$ runs over the points with 
$\opn{dim} \overline{\{ x \}} = i$. Each
$\mcal{K}_{X}(x)$ is a quasi-coherent 
$\mcal{O}_X$-module with support the closed set
$\overline{\{ x \}}$, and as module it is an injective hull of 
the residue field $\bsym{k}(x)$ as $\mcal{O}_{X, x}$-module. 
Since each $\mcal{K}_{X}^{-i}$ is an injective $\mcal{O}_X$-module 
it follows that 
$\mrm{D} \mcal{M} \cong
\mcal{H}om_{\mcal{O}_X}(\mcal{M}, \mcal{K}_{X})$.

Let us prove that (ii) implies (iii). Let
$\mcal{N}$ be the coherent sheaf $\mrm{H}^0 \mrm{D} \mcal{M}$;
so $\mcal{M} \cong \mrm{D} \mcal{N}$
in $\msf{D}(\cat{Mod} \mcal{O}_X)$. 
For every $i$ the sheaf
\[ \mcal{H}om_{\mcal{O}_X}(\mcal{N}, \mcal{K}_X^{-i}) \cong
\bigoplus_{\opn{dim} \overline{\{ x \}}\, =\, i}
 \mcal{H}om_{\mcal{O}_X}(\mcal{N}, \mcal{K}_{X}(x)) \]
is flasque and pure of dimension $i$ (or it is $0$).
Therefore 
\[ \mrm{E} \mcal{H}om_{\mcal{O}_X}(\mcal{N}, \mcal{K}_X)
= \mcal{H}om_{\mcal{O}_X}(\mcal{N}, \mcal{K}_X) . \]
We conclude that
$\mrm{E} \mcal{M} \cong \mrm{E} \mrm{D} \mcal{N}
\cong \mrm{D} \mcal{N} \cong \mcal{M}$
in $\msf{D}(\cat{Mod} \mcal{O}_X)$. 

Now let us prove that (iii) implies (i). We are assuming that
$\mrm{E} \mcal{M} \cong \mcal{M}$
in $\msf{D}(\cat{Mod} \mcal{O}_X)$. 
It suffices to prove that 
$M \in {}^{p}\msf{D}^{\mrm{b}}_{\mrm{f}}(\cat{Mod} C)^{0}$,
where $C := \Gamma(X, \mcal{O}_X)$
and $M := \mrm{R} \Gamma(X, \mcal{M})$. 
By \cite[Lemma 7.8]{YZ4} it is enough to show that 
$\opn{dim} \mrm{H}^{-i} M \leq i$
and $\mrm{H}^j_{\msf{M}_i} M = 0$
for all $i$ and all $j < -i$.
Because $\mrm{R} \Gamma(X, -)$ can be computed using flasque 
resolutions we get
$M \cong \Gamma(X, \mrm{E} \mcal{M})$. 
But 
$\opn{dim} \Gamma(X, \mrm{E} \mcal{M})^{-i} \leq i$,
and hence $\opn{dim} \mrm{H}^{-i} M \leq i$. 
For the vanishing we use the fact that 
$\mrm{R} \Gamma_{\msf{M}_i}$ can be computed by flasque 
resolutions (cf.\ \cite[Definition 1.13]{YZ3}), 
and that each $C$-module
$\Gamma(X, \mrm{E} \mcal{M})^{-i}$
is a flasque $C$-module and is pure of $\opn{dim} = i$
(unless it is zero). This gives
\[ \mrm{R} \Gamma_{\msf{M}_i} M \cong
\Gamma_{\msf{M}_i} \Gamma(X, \mrm{E} \mcal{M}) =
\Gamma(X, \mrm{E} \mcal{M})^{\geq -i} . \]
Thus $\mrm{H}^j_{\msf{M}_i} M = 0$ for all $j < -i$.
\end{proof}

\begin{rem}
Theorem \ref{thm13.4} implies that the category of Cohen-Macaulay 
complexes is an abelian subcategory of 
$\msf{D}^{\mrm{b}}_{\mrm{c}}(\cat{Mod} \mcal{O}_X)$. 
This fact that seems to have eluded Grothendieck. 
Sastry has communicated to us another 
proof of the equivalence (iii) $\Leftrightarrow$ (ii) in 
Theorem \ref{thm13.4}, using local duality at closed points of $X$. 
\end{rem}

We now go to some noncommutative examples. 

\begin{exa} \label{exa13.4}
Let $X$ be any separated $\k$-scheme of finite type and $\mcal{A}$ 
a coherent $\mcal{O}_{X}$-algebra. Then $\mcal{A}$ is a 
differential quasi-coherent $\mcal{O}_{X}$-ring of finite type;
as filtration $F$ we may take the trivial filtration where
$F_{-1} \mcal{A} := 0$ and $F_{0} \mcal{A} := \mcal{A}$. By Theorem 
\ref{thm13.8} the rigid dualizing complex $\mcal{R}_{\mcal{A}}$ 
exists. In this case $\mcal{R}_{\mcal{A}}$ can 
be chosen to be 
$\mrm{R} \mcal{H}om_{\mcal{O}_X}(\mcal{A}, \mcal{R}_X)$
where $\mcal{R}_X 
\in \msf{D}^{\mrm{b}}_{\mrm{c}}(\cat{Mod} \mcal{O}_X)$ 
is the central rigid dualizing complex of $X$.
If $\mcal{A}$ happens to be an Azumaya $\mcal{O}_{X}$-algebra 
then 
$\mcal{R}_{\mcal{A}} \cong \mcal{A} \otimes_{\mcal{O}_{X}}
\mcal{R}_{X}$
(see \cite[Theorem 6.2]{YZ3}).
\end{exa}

\begin{exa} \label{exa13.3}
Suppose $\opn{char} \k = 0$ and $X$ is a smooth separated
scheme. Then
the ring $\mcal{D}_{X}$ of differential operators on $X$ is 
a differential quasi-coherent $\mcal{O}_{X}$-ring of finite type. 
As filtration $F$ we may take the order filtration, in which
$F_{0} \mcal{D}_{X} := \mcal{O}_{X}$,
$F_{1} \mcal{D}_{X} := \mcal{O}_{X} \oplus \mcal{T}_{X}$,
and
$F_{i + 1} \mcal{D}_{X} := F_{i} \mcal{D}_{X} \cdot 
F_{1} \mcal{D}_{X}$.
Here 
$\mcal{T}_{X} := \mcal{H}om_{\mcal{O}_{X}}
(\Omega^1_{X / \k}, \mcal{O}_{X})$ 
is the tangent sheaf.
The product turns out to be
\[ (X^2, \mcal{D}_{X}^{\, \mrm{e}}) \cong
(X^2, \mcal{D}_{X^2}(\mcal{O}_{X} \boxtimes \bsym{\omega}_X)) \]
where $n := \opn{dim} X$,
$\bsym{\omega}_X := \Omega^n_{X / \k}$,
$\mcal{O}_{X} \boxtimes \bsym{\omega}_X :=
\mcal{O}_{X^2} \otimes_{\mrm{p}_2^{-1} \mcal{O}_{X}} 
\mrm{p}_2^{-1} \bsym{\omega}_X$
and
$\mcal{D}_{X^2}(-)$
is the ring of twisted differential operators. By 
\cite[Proposition 2.6]{Ye4} the rigid dualizing complex 
of $\mcal{D}_{X}$ is 
$\mcal{D}_{X}[2n] \in 
\msf{D}^{\mrm{b}}_{\mrm{c}}(\cat{Mod} \mcal{D}_{X}^{\mrm{\, e}})$.
\end{exa}

\begin{exa}
Let $X$ be any finite type separated $\k$-scheme, $\mcal{M}$ a 
coherent $\mcal{O}_X$-module and 
$\mcal{D}_X(\mcal{M}) := \mcal{D}iff_{\mcal{O}_X / \k}
(\mcal{M}, \mcal{M})$
the ring of differential operators from $\mcal{M}$ to itself
(cf.\ \cite{EGA-IV}). The order filtration $F$ makes 
$\mcal{D}_X(\mcal{M})$ into a differential quasi-coherent 
$\mcal{O}_X$-ring, but unless $\opn{char} \k =0$, 
$X$ is smooth and $\mcal{M}$ is locally free, this is usually not
of finite type. However if we take 
$\mcal{N} := F_1 \mcal{D}_X(\mcal{M})$, then 
the subring $\mcal{O}_X \bra{\mcal{N}} \subset \mcal{D}_X(\mcal{M})$ 
is a differential quasi-coherent $\mcal{O}_X$-ring of finite type.
\end{exa}

\begin{exa} \label{exa13.5}
Generalizing Example \ref{exa13.3}, 
suppose $X$ is a separated finite type $\k$-scheme 
and $\mcal{L}$ is a 
coherent $\mcal{O}_{X}$-module endowed with a $\k$-linear 
Lie bracket $[-, -]$ and an $\mcal{O}_{X}$-linear Lie homomorphism
$\alpha : \mcal{L} \to \mcal{T}_{X}$, satisfying the conditions 
stated in \cite[Example 2.6]{YZ4} on affine open sets. Such 
$\mcal{L}$ is called a Lie algebroid on $X$. 
By \cite[Proposition 5.17 and Corollary 5.13]{YZ4}
the universal enveloping algebras 
$\mrm{U}(C; L)$ sheafify to a differential quasi-coherent 
$\mcal{O}_{X}$-ring of finite type
$\mrm{U}(\mcal{O}_{X}; \mcal{L})$.

Now suppose $\opn{char} \k = 0$ and $X$ is a smooth scheme of 
dimension $n$. When $\mcal{L} = \mcal{T}_{X}$ then 
$\mrm{U}(\mcal{O}_{X}; \mcal{L}) = \mcal{D}_X$ 
as in Example \ref{exa13.3}. 
More generally when $\mcal{L}$ is a locally free 
$\mcal{O}_{X}$-module of rank $r$
the rigid dualizing complex 
$\mcal{R}_{\mrm{U}(\mcal{O}_{X}; \mcal{L})}$ of
$\mrm{U}(\mcal{O}_{X}; \mcal{L})$ was computed by 
Chemla \cite{Ch2} for $X$ affine. Since the rigidifying 
isomorphism used there was canonical it glues, and we obtain
\begin{equation} \label{eqn13.2}
\mcal{R}_{\mrm{U}(\mcal{O}_{X}; \mcal{L})} =
\mrm{U}(\mcal{O}_{X}; \mcal{L}) \otimes_{\mcal{O}_{X}} 
\Omega^n_{X / \k} \otimes_{\mcal{O}_{X}} 
\bigl( \bwedge^r_{\mcal{O}_{X}} \mcal{L} \bigr) [n + r] .
\end{equation}
\end{exa}

\begin{exa}
This is a special case of Example \ref{exa13.5}.
Suppose $\opn{char} \k = p > 0$ and $X$ is smooth of dimension 
$n$. The sheaf of rings 
$\mcal{A} := \mrm{U}(\mcal{O}_{X}; \mcal{T}_X)$ 
is called the sheaf of {\em crystalline differential operators}
(see \cite{BFG}). There is a ring homomorphism 
$\mcal{A} \to \mcal{D}_X$, which is neither injective nor 
surjective. Also $\mcal{A}$ is finite over its center.
Equation \ref{eqn13.2} shows that
the rigid dualizing complex of $\mcal{A}$ is $\mcal{A}[2n]$, 
just like in Example \ref{exa13.3}. 
\end{exa}

\begin{exa} \label{exa13.6}
Assume $\opn{char} \k = 0$. Let $C := \k[x]$ with $x$ an 
indeterminate, and 
$A := \mcal{D}(C)$ the first Weyl algebra. Writing 
$y := \frac{\partial}{\partial x} \in A$ we have
$A \cong \k \bra{x, y}/([y, x] - 1)$. Consider the filtration $F$ 
on $A$ in which 
$\opnt{deg}^F(x) = \opnt{deg}^F(y) = 1$, and let
$\wtil{A} := \opnt{Rees}^F A$. Then $\wtil{A}$ is the homogenized 
Weyl algebra, with linear generators 
$u := x t$, $v := y t$ and $t$, such that $t$ is central and 
$[v, u] = t^2$. Inside $\wtil{A}$ we have the commutative subring
$\wtil{C} := \opnt{Rees}^F C \cong \k[u, t]$. 

Since $\wtil{A}$ is a differential $\wtil{C}$-ring the 
localizations $\wtil{A}_t$ and $\wtil{A}_u$ exist. Their degree 
$0$ components glue to a quasi-coherent sheaf of rings 
$\mcal{A}$ on $\opn{Proj} \wtil{C} \cong \mbf{P}^1_{\k}$. We 
obtain a differential quasi-coherent ringed $\k$-scheme of finite 
type $(\mbf{P}^1, \mcal{A})$. The restriction to the open set
$\{ t \neq 0 \}$ recovers the ringed scheme
$(\mbf{A}^1, \mcal{D}_{\mbf{A}^1})$, and in particular
$\Gamma(\{ t \neq 0 \}, \mcal{A}) = A$. Thus 
$(\mbf{P}^1, \mcal{A})$ can be viewed as a partial 
compactification of $(\mbf{A}^1, \mcal{D}_{\mbf{A}^1})$. 
Note that the ringed 
scheme $(\mbf{P}^1, \mcal{A})$ is not isomorphic to 
$(\mbf{P}^1, \mcal{D}_{\mbf{P}^1})$.
Indeed, consider the affine open set 
$W := \{ u \neq 0 \} \subset \mbf{P}^1$.
Let $z := u^{-1} t$, so 
$\Gamma(W, \mcal{O}_{\mbf{P}^1}) = \k[z]$,
$\Gamma(W, \mcal{A}) = \k \bra{z, u^{-1} v}$ and
$\Gamma(W, \mcal{D}_{\mbf{P}^1}) = 
\k \bra{z, \frac{\partial}{\partial z}}$.
A calculation shows that 
$u^{-1} v = -z^3 \frac{\partial}{\partial z}$, and that the
subring 
$\k \bra{z, -z^3 \frac{\partial}{\partial z}}$
is not isomorphic as $\k[z]$-ring to the Weyl algebra 
$\k \bra{z, \frac{\partial}{\partial z}}$. 
\end{exa}

\begin{exa}
Consider another geometric object associated to the situation
of Example \ref{exa13.6}: 
the projective spectrum $\opn{Proj} \wtil{A}$ in the sense of 
Artin-Zhang \cite{AZ}. We claim that $(\mbf{P}^1, \mcal{A})$ is an 
``open subscheme'' of the ``complete surface'' 
$\opn{Proj} \wtil{A}$, whose complement consists of one point. 
To make this statement precise we have to pass to abelian 
categories. Recall that 
$\cat{QGr} \wtil{A}$ is the quotient category
$\cat{GrMod} \wtil{A} / \{ \mfrak{m}\text{-torsion} \}$, where
$\mfrak{m}$ is the augmentation ideal of $\wtil{A}$; and
$\opn{Proj} \wtil{A}$ is the geometric object such that 
``$\cat{QCoh} \opn{Proj} \wtil{A} = \cat{QGr} \wtil{A}$''.
Because $\cat{QCoh} \mcal{A}$ is obtained by gluing the categories
$\cat{QCoh} \mcal{A}|_{\{ t \neq 0 \}}$ and 
$\cat{QCoh} \mcal{A}|_{\{ u \neq 0 \}}$, we see that 
$\cat{QCoh} \mcal{A}$ is equivalent 
to the quotient of $\cat{QGr} \wtil{A}$ by the localizing 
subcategory $\{ \mfrak{a}\text{-torsion} \}$, 
where $\mfrak{a}$ is the two-sided ideal
$(t, u) \subset  \wtil{A}$.
\end{exa}

Let's continue with the setup of the last two examples. 
The subcategory of noetherian objects in $\cat{QGr} \wtil{A}$
is denoted by $\cat{qgr} \wtil{A}$; it is equivalent to the 
category 
$\cat{GrMod}_{\mrm{f}} \wtil{A} / \{ \mfrak{m}\text{-torsion} \}$.
In \cite{KKO} there is a duality 
\[ \mrm{D}_{\cat{qgr}} : 
\msf{D}^{\mrm{b}}(\cat{qgr} \wtil{A}) \to
\msf{D}^{\mrm{b}}(\cat{qgr} \wtil{A}^{\mrm{op}}) . \]
The formula (with a slight adjustment) is given below. Let us 
denote by 
$\pi : \cat{GrMod} \wtil{A} \to \cat{QGr} \wtil{A}$ the 
localization functor and by 
$\wtil{\Gamma} : \cat{QGr} \wtil{A} \to \cat{GrMod} \wtil{A}$
its adjoint
\[ \wtil{\Gamma} \mcal{M} := \boplus_{i \geq 0} 
\opn{Hom}_{\cat{QGr} \wtil{A}}(\pi \wtil{A}(-i), \mcal{M}) . \]
Denote by $\pi^{\mrm{op}}$ and $\wtil{\Gamma}^{\mrm{op}}$ the 
corresponding right versions. Let $R_{\wtil{A}}$ be the balanced 
dualizing complex of $\wtil{A}$; so 
$R_{\wtil{A}} \cong \bsym{\omega}_{\wtil{A}}[3]$
where $\bsym{\omega}_{\wtil{A}}$ is the dualizing bimodule of 
$\wtil{A}$. Then
\begin{equation} \label{eqn5.9}
\mrm{D}_{\cat{qgr}} \mcal{M} := \pi^{\mrm{op}}
\opn{RHom}^{\mrm{gr}}_{\wtil{A}}(\mrm{R} \wtil{\Gamma} \mcal{M},
R_{\wtil{A}}[-1]) . 
\end{equation}
With this normalization the duality $\mrm{D}_{\cat{qgr}}$
is compatible with the global $\k$-linear duality, namely
\[ \opn{RHom}_{\cat{QGr} \wtil{A}}(\pi \wtil{A}, \mcal{M})^*
\cong
\opn{RHom}_{\cat{QGr} \wtil{A}^{\mrm{op}}}
(\pi^{\mrm{op}} \wtil{A}^{\mrm{op}}, 
\mrm{D}_{\cat{qgr}} \mcal{M}) . \]
(Actually the formula (\ref{eqn5.9}) gives a duality for any 
noetherian connected graded $\k$-algebra $\wtil{A}$ 
admitting a balanced dualizing complex, be it regular or not; 
cf.\ \cite{YZ1}.)

Kazhdan has recently asked us whether the duality 
$\mrm{D}_{\cat{qgr}}$
is compatible with the duality on $(\mbf{P}^1, \mcal{A})$
via the ``open embedding'' 
$g: (\mbf{P}^1, \mcal{A}) \inj \opn{Proj} \wtil{A}$. 
We can answer positively:

\begin{prop} \label{prop5.5}
The diagram
\[ \begin{CD}
\msf{D}^{\mrm{b}}(\cat{qgr} \wtil{A})
@>{ \mrm{D}_{\cat{qgr}} }>>
\msf{D}^{\mrm{b}}(\cat{qgr} \wtil{A}^{\mrm{op}}) \\
@V{ g^* }VV @V{ (g^{\mrm{op}})^* }VV \\
\msf{D}^{\mrm{b}}(\cat{Coh} \mcal{A})
@>{ \mrm{D}_{\mcal{A}} }>>
\msf{D}^{\mrm{b}}(\cat{Coh} \mcal{A}^{\mrm{op}})
\end{CD} \]
where $\mrm{D}_{\mcal{A}}$ is the duality 
determined by the rigid dualizing complex $\mcal{R}_{\mcal{A}}$,
is commutative.
\end{prop}

\begin{proof}
Consider the graded bimodule $\wtil{A}(i)$ for any integer $i$. It 
is generated by the element $1$ which is in degree $-i$. 
When we invert $t$ we get a generator 
$t^{i} \in (\wtil{A}(i)_t)_0$. Likewise when 
we invert $u$ we get a generator 
$u^{i} \in (\wtil{A}(i)_u)_0$.
Therefore we obtain a sheaf $\mcal{A}(i)$
of $\mcal{A}$-bimodules on $\mbf{P}^1$, such that
\[ \mcal{A}(i)|_{\{ t \neq 0 \}} = 
\mcal{A}|_{\{ t \neq 0 \}} \cdot t^i
\cong \mcal{A}|_{\{ t \neq 0 \}} \]
and 
\[ \mcal{A}(i)|_{\{ u \neq 0 \}} 
= \mcal{A}|_{\{ u \neq 0 \}} \cdot u^i
= u^i \cdot \mcal{A}|_{\{ u \neq 0 \}} , \]
and the gluing on the open set $\{ u t \neq 0 \}$ is 
multiplication by $t^{-i} u^i = x^i$. 
Note that the $\mcal{A}$-bimodule $\mcal{A}(i)$ is not 
locally centrally generated if $i \neq 0$.

Since $\mcal{A}(i)$ is a differential 
$\mcal{O}_{\mbf{P}^1}$-bimodule we can 
view it as a coherent $\mcal{A}^{\mrm{e}}$-module supported on the 
diagonal $\Delta(\mbf{P}^1) \subset (\mbf{P}^1)^2$.
As left $\mcal{A}$-module we have 
$\mcal{A}(i) \cong g^* \pi \wtil{A}(i)$; and likewise on the 
right.

Let us compute the balanced dualizing complex $R_{\wtil{A}}$. 
As mentioned above, since $\wtil{A}$ is an AS-regular 
algebra it follows that 
$R_{\wtil{A}} \cong \bsym{\omega}_{\wtil{A}}[3]$, 
and 
$\bsym{\omega}_{\wtil{A}} = \wtil{A}^{\sigma}(-3)$ 
for some graded algebra automorphism $\sigma$.
See \cite[Proposition 1.1]{Ye4}. By \cite[Corollary 3.6]{Ye4} 
the automorphism $\sigma$ is trivial on 
the center, so it is $\k[t]$-linear.
Now $(\wtil{A}_t)_0 \cong A$, and we know that 
$R_A \cong A[2]$ by \cite[Theorem 2.6]{Ye4}. Hence $\sigma$ is 
in fact trivial and
$R_{\wtil{A}} \cong \wtil{A}(-3)[3]$.

Let us denote 
$B := \Gamma(\{ u \neq 0 \}, \mcal{A}) = 
(\wtil{A}_u)_0$. 
The localization of rigid dualizing complexes (see 
\cite[Theorem 6.2]{YZ4}) 
shows that the rigid dualizing complex of $B$ is 
\[ R_B \cong ((R_{\wtil{A}})_u)_0[-1] \cong
(\wtil{A}(-3)[2]_u)_0 \cong B u^{-3} [2] \]
and hence
\[ \mcal{R}_{\mcal{A}}|_{\{ u \neq 0 \}} \cong
\mcal{A}(-3)[2]|_{\{ u \neq 0 \}} . \]
On the open set $\{ t \neq 0 \}$ we have
\[ \mcal{R}_{\mcal{A}}|_{\{ t \neq 0 \}} \cong
\mcal{A}[2]|_{\{ t \neq 0 \}} \cong
\mcal{A}(-3)[2]|_{\{ t \neq 0 \}} . \]
The bimodule 
$\mcal{A}|_{\{ t u \neq 0 \}} \cong 
\mcal{A}(-3)|_{\{ t u \neq 0 \}}$ 
has only constant automorphisms (because $\mrm{Z}(A_x) = \k$), 
so it follows that $\mcal{R}_{\mcal{A}} \cong \mcal{A}(-3)[2]$
in $\msf{D}(\cat{Mod} \mcal{A}^{\mrm{e}})$.

The category $\msf{D}^{\mrm{b}}(\cat{qgr} \wtil{A})$
is a quotient of 
$\msf{D}^{\mrm{b}}(\cat{GrMod}_{\mrm{f}} \wtil{A})$. 
We know that $\wtil{A}$ is regular, and hence any object in 
$\msf{D}^{\mrm{b}}(\cat{GrMod}_{\mrm{f}} \wtil{A})$
is isomorphic to a bounded complex of finite 
free graded modules. Now for any $i$ we get
\[ (g^{\mrm{op}})^*\, \mrm{D}_{\cat{qgr}}\, \pi\, \wtil{A}(i) \cong 
\mrm{D}_{\mcal{A}}\, g^*\, \pi\, \wtil{A}(i)  \cong
\mcal{A}(-3 - i)[2] \]
in $\msf{D}^{\mrm{b}}(\cat{Coh} \mcal{A}^{\mrm{op}})$. Maps
$\wtil{A}(i) \to \wtil{A}(j)$ in 
$\cat{GrMod}_{\mrm{f}} \wtil{A}$ 
are right multiplication by elements 
$a \in \wtil{A}(j - i)$, and these are sent under both dualities 
to left multiplication by the same $a$.
Hence
$(g^{\mrm{op}})^*\, \mrm{D}_{\cat{qgr}}\, \mcal{P} \cong
\mrm{D}_{\mcal{A}}\, g^*\, \mcal{P}$
for any object 
$\mcal{P} \in \msf{D}^{\mrm{b}}(\cat{qgr} \wtil{A})$.
Likewise for morphisms.
\end{proof}

If $A$ is a finite $\k$-algebra of finite global dimension then 
the rigid dualizing complex of $A$ is
$A^* := \opn{Hom}_{\k}(A, \k)$. It is known that 
$M \mapsto A^* \otimes^{\mrm{L}}_{A} M$
is the Serre functor of $\msf{D}^{\mrm{b}}_{\mrm{f}}(\cat{Mod} A)$. 
Namely $A^* \otimes^{\mrm{L}}_{A} -$ is an auto-equivalence of 
$\msf{D}^{\mrm{b}}_{\mrm{f}}(\cat{Mod} A)$, and 
there is a bifunctorial nondegenerate pairing
\[ \opn{Hom}_{\msf{D}(\cat{Mod} A)}(M, N) \times
\opn{Hom}_{\msf{D}(\cat{Mod} A)}(N, A^* \otimes^{\mrm{L}}_{A} M)
\to \k \]
for $M, N \in \msf{D}^{\mrm{b}}_{\mrm{f}}(\cat{Mod} A)$.
Cf.\ \cite{MY} and \cite{BO}. 
Likewise, for a smooth $n$-dimensional
projective scheme $X$ the rigid dualizing complex is 
$\Omega^n_{X / \k}[n]$ (see Example \ref{exa7.4}), 
and the Serre functor is 
$\mcal{M} \mapsto \Omega^n_{X / \k}[n] \otimes_{\mcal{O}_X}
\mcal{M}$. 
Kontsevich has asked us whether this is true in greater generality. 
We have the following partial answer.

Let $\mcal{A}$ be a coherent $\mcal{O}_X$-ring.
We say that $\mcal{A}$ is {\em regular} 
if for any affine open set $U$ the ring
$\Gamma(U, \mcal{A})$ has finite global dimension. As explained 
in Example \ref{exa13.4}, the rigid dualizing complex 
$\mcal{R}_{\mcal{A}}$ can be assumed to live in 
$\msf{D}^{\mrm{b}}_{\mrm{c}}(\cat{Mod}\,
\mcal{A} \otimes_{\mcal{O}_{X}} \mcal{A}^{\mrm{op}})$.

\begin{prop}
Let $X$ be a projective $\k$-scheme, and let $\mcal{A}$ be
a regular coherent $\mcal{O}_X$-algebra, with 
rigid dualizing complex 
$\mcal{R}_{\mcal{A}} \in \msf{D}^{\mrm{b}}_{\mrm{c}}(\cat{Mod}\,
\mcal{A} \otimes_{\mcal{O}_{X}} \mcal{A}^{\mrm{op}})$.
Then 
$\mcal{M} \mapsto \mcal{R}_{\mcal{A}} \otimes_{\mcal{A}}^{\mrm{L}} 
\mcal{M}$
is a Serre functor of 
$\msf{D}^{\mrm{b}}_{\mrm{c}}(\cat{Mod} \mcal{A})$.
\end{prop}

\begin{proof}
Choose an ample $\mcal{O}_X$-module $\mcal{L}$, and let
$\mcal{M}(1) := \mcal{M} \otimes_{\mcal{O}_X} \mcal{L}$.
Define 
$\wtil{A} := \k \oplus \boplus_{i \geq 1} 
\Gamma(X, \mcal{A}(i))$,
which is a connected graded noetherian $\k$-algebra, finite over 
its center. Then
$\cat{Coh} \mcal{A} \cong \cat{qgr} \wtil{A}$,
and this category has finite global dimension. 
Let $R_{\wtil{A}}$ be the balanced dualizing complex of 
$\wtil{A}$. Then 
$\mcal{R}_{\mcal{A}} \cong \pi R_{\wtil{A}}[-1]$,
and the claim follows from \cite[Theorem A.4]{NV}.
\end{proof}


\end{document}